\documentclass[a4paper,12pt]{amsart}
\usepackage[frenchb,english]{babel}
\usepackage[latin1]{inputenc}
\usepackage{pdfsync}
\usepackage[T1]{fontenc}   
\usepackage{amsmath}
\usepackage{amssymb}
\usepackage{amsxtra}
\usepackage{amscd}
\usepackage{amsthm}
\usepackage{amsfonts}
\usepackage{eucal}
\usepackage[all]{xy}
\usepackage{graphicx}
\usepackage{comment}
\usepackage{epsfig}
\usepackage{psfrag}
\usepackage{mathrsfs}
\usepackage{amscd}
\usepackage{rotating}
\usepackage{lscape}
\usepackage{amsbsy}
\usepackage{verbatim}
\usepackage{moreverb}
\usepackage{url}
\usepackage{color}
\usepackage{xcolor}
\usepackage{stmaryrd}
\usepackage{enumitem}

\newcommand{\nc}{\newcommand}
\nc{\renc}{\renewcommand}

\nc\restr[2]{{ 
  \left.\kern-\nulldelimiterspace    #1  
  \vphantom{\big|}  
  \right|_{#2}  
  }}

\newtheorem{lem}[subsubsection]{Lemma}
\newtheorem{prop}[subsubsection]{Proposition}

\newtheorem{thm}[subsubsection]{Theorem}

\newtheorem{rem}[subsubsection]{Remark}

\newtheorem{notation}[subsubsection]{Notation}

\renc{\sec}{\section}
\nc{\ssec}{\subsection}
\nc{\sssec}{\subsubsection}

\theoremstyle{definition}
\newtheorem{defi}[subsubsection]{Definition}
\newtheorem{example}[subsubsection]{Example}

\numberwithin{equation}{section}

\nc{\on}{\operatorname}

\nc\wt{\widetilde}
\nc\wh{\widehat}
\nc\ol{\ov}
\nc{\oc}[1]{{\overset{\circ}{#1}}}
\nc{\ov}[1]{{\overline{#1}}}
\nc{\isor}[1]{{\xrightarrow[\raisebox{0.25 em}{\smash{\ensuremath{\sim}}}]{#1}}}
\nc{\modmod}{/ \! \! /}

\nc{\mc}{\mathcal}
\nc{\mf}{\mathfrak}
\nc{\mr}{\mathrm}
\nc{\mb}{\mathbb}
\nc{\mbf}{\mathbf}
\nc{\ms}{\mathscr}

\nc{\R}{{\mathbb R}}
\nc{\Z}{{\mathbb Z}}
\nc{\N}{{\mathbb N}}
\nc{\C}{{\mathbb C}}
\nc{\Q}{{\mathbb Q}}

\nc{\Fq}{{\mathbb F}_q}
\nc{\Fl}{{\mathbb F}_\ell}
\nc{\Fqbar}{\ol{{\mathbb F}_q}}
\nc{\Flbar}{\ol{{\mathbb F}_\ell}}
\nc{\Zl}{{\mathbb Z}_\ell}
\nc{\Zlbar}{\ol{{\mathbb Z}_\ell}}
\nc{\Ql}{E}
\nc{\Qlbar}{\ol{{\mathbb Q}_\ell}}
\nc{\hl}{\overset{\leftarrow}h{}}
\nc{\hr}{\overset{\rightarrow}h{}}
\nc{\Gr}{{\on{Gr}}}
\nc{\Hecke}{\on{Hecke}}
 \nc{\Hom}{\on{Hom}}
 \nc{\Coker}{\on{Coker}}
 \nc{\Ker}{\on{Ker}}
 \nc{\Lie}{\on{Lie}}
\nc{\Loc}{\on{Loc}}
\nc{\Pic}{\on{Pic}}
\nc{\Bun}{\on{Bun}}
\nc{\IC}{\on{IC}}
\nc{\Aut}{\on{Aut}}
\nc{\Perv}{\on{Perv}}
\nc{\pos}{{\on{pos}}}
\nc{\Sym}{\on{Sym}}

\nc{\ta} {{}^\tau}
\nc {\tu}[1]{{}^{\tau^{#1}}\!}

\nc{\Id}{\on{Id}}
\nc{\Fil}{\on{Fil}}
\nc{\pr}{\on{pr}}
\nc{\Res}{\on{Res}}
\nc{\cusp}{\on{cusp}}
\nc{\Frob}{\on{Frob}}
\nc{\diag}{\Delta}
\nc{\gr}{\on{gr}}
\nc{\Inj}{\on{Inj}}
\nc{\Bl}{\on{Bl}}
\nc{\dem}{\noindent {\bf Proof. }}
\nc{\cqfd}{{\ }\hfill $\square$ \vskip 1mm}
\nc{\s}[1]{\langle #1 \rangle}
\nc{\Cht}{\on{Cht}}
\nc{\isom}{\overset {\thicksim}{\to}}
\nc{\sm}{\smallsetminus}

\begin{document}

\title{Smoothness of cohomology sheaves of stacks of shtukas}
\author{Cong Xue}
\address{Cong Xue: Institut de Mathématiques de Jussieu-Paris Rive Gauche, Université de Paris, CNRS, France.}
\email{cong.xue@imj-prg.fr}

\maketitle

\begin{abstract}
We prove, for all reductive groups,
that the cohomology sheaves with compact support of stacks of shtukas are ind-smooth over $(X \sm N)^I$ and that their geometric generic fibers are equiped with an action of %$\on{Weil}(\eta, \ov{\eta})^I$ factors through 
$\on{Weil}(X \sm N, \ov{\eta})^I$.
Our method does not use any compactification of stacks of shtukas.
\end{abstract}

\section*{Introduction}
Let $X$ be a smooth projective geometrically connected curve over a finite field $\Fq$. We denote by $F$ the function field of $X$.

Let $G$ be a connected reductive group over $F$.
Let $\ell$ be a prime number not dividing $q$. Let $E$ be a finite extension of $\mb Q_{\ell}$ containing a square root of $q$. 
%We denote by $\mc O_E$ the ring of integers of $E$.

\quad

In the introduction, we assume that $G$ is split to simplify the notation. Let $\wh G$ be the Langlands dual group of $G$ over $E$. 
Let $I$ be a finite set and $W$ be a finite dimensional $E$-linear representation of $\wh G^I$. Let $N \subset X$ be a finite subscheme.
Varshavsky (\cite{var}) and V. Lafforgue (\cite{vincent}) defined the stack classifying $G$-shtukas with level structure on $N$: $$\Cht_{G, N, I, W} \text{ over } (X \sm N)^I$$ and its degree $j \in \Z$ cohomology sheaf with compact support with $E$-coefficients: $$\mc H_{G, N, I, W}^{j} \text{ over } (X \sm N)^I $$ which is an inductive limit of constructible $E$-sheaves. 
The cohomology sheaf $\mc H_{G, N, I, W}^{j}$ is equiped with an action of the Hecke algebra, %$C_c(K_{N} \backslash G(\mb A) / K_{N}, E)$ 
and equiped with an action of the partial Frobenius morphisms.

Let $\eta$ be the generic point of $X$ and $\ov{\eta}$ be a geometric point over $\eta$.
Let $\eta_I$ be the generic point of $X^I$ and $\ov{\eta_I}$ be a geometric point over $\eta_I$. Then $\restr{  \mc H_{G, N, I, W}^j  }{  \ov{\eta_I}  }$ is equiped with an action of $\on{Weil}(\eta_I, \ov{\eta_I})$ and an action of the partial Frobenius morphisms.

%When $G$ is split, 
In \cite{coho-filt}, by a variant of a lemma of Drinfeld, 
we equiped $\restr{  \mc H_{G, N, I, W}^j  }{  \ov{\eta_I}  }$ with an action of $\on{Weil}(\eta, \ov{\eta})^I$. The proof uses the fact that $\restr{  \mc H_{G, N, I, W}^j  }{  \ov{\eta_I}  }$ is of finite type as a module over a Hecke algebra (proved in $loc.cit.$). The case of non-split groups of the last assertion has not yet been written.

\quad

In this paper, %we do not need the assumption that $G$ is split.
in Section 1 we give a different proof of the fact that $\restr{  \mc H_{G, N, I, W}^j  }{  \ov{\eta_I}  }$ is equiped with an action of $\on{Weil}(\eta, \ov{\eta})^I$ (Proposition \ref{prop-FWeil-factors-though-Weil-I}). This proof is based on the Eichler-Shimura relations in \cite{vincent}. It does not use the fact that $\restr{  \mc H_{G, N, I, W}^j  }{  \ov{\eta_I}  }$ is of finite type as module over a Hecke algebra. The new proof has the advantage that it is easy to generalize to the case of non-split groups.
 
In Section 1, we also prove that $\restr{\mc H_{G, N, I, W}^j}{(\ov{\eta})^{I}} ,$ the restriction of $\mc H_{G, N, I, W}^j$ to the scheme $(\ov{\eta})^{I} :=\ov{\eta} \times_{\on{Spec} \Fqbar} \cdots \times_{\on{Spec} \Fqbar} \ov{\eta} $ ($I$ times), is ind-smooth (i.e. an inductive limit of smooth $E$-sheaves). Then Proposition \ref{prop-FWeil-factors-though-Weil-I} implies that $\restr{\mc H_{G, N, I, W}^j}{(\ov{\eta})^{I}}$ is a constant sheaf.

In Section 2, we generalize the results in Section 1. For any partition $I = I_1 \sqcup I_2$ and for any geometric point $\ov u$ over a closed point of $(X \sm N)^{I_2}$, we prove in Proposition \ref{prop-H-ov-eta-I-1-ov-s-I-2-is-constant} that the restriction of $\mc H_{G, N, I, W}^j$ to the scheme $(\ov{\eta})^{I_1} \times_{\on{Spec} \Fqbar} \ov u$ is a constant sheaf, where $(\ov{\eta})^{I_1} = \ov{\eta} \times_{\on{Spec} \Fqbar} \cdots \times_{\on{Spec} \Fqbar} \ov{\eta} $ ($I_1$ times).
%$\restr{\mc H_{G, N, I, W}^j}{(\ov{\eta})^{I_1} \times_{\on{Spec} \Fqbar} \ov u}$, 
When $I_2$ is the empty set, we recover Section 1.
 
%In particular, the restriction of $\mc H_{G, N, I, W}^j$ to the scheme $(\ov{\eta})^{I}$ is constant. 

\quad

%With the help of these facts and the "Zorro" lemma, 
In Sections 3-4, we prove that $\mc H_{G, N, I, W}^{j}$ is ind-smooth over $(X \sm N)^I$ in the following sense (which is equivalent to be an inductive limit of smooth $E$-sheaves by Lemma \ref{lem-ind-smooth-specialization}):
\begin{thm} \label{thm-H-I-W-smooth-intro}  (Theorem \ref{thm-H-I-W-indsmooth})
For any geometric point $\ov x$ of $(X \sm N)^I$ and any specialization map 
$\mf{sp}_{\ov x}: \ov{\eta_I} \rightarrow \ov{x}$,
the induced morphism
\begin{equation}  
\mf{sp}_{\ov x}^*: \restr{ \mc H_{G, N, I, W}^j  }{ \ov{x}  } \rightarrow \restr{ \mc H_{G, N, I, W}^j  }{ \ov{\eta_I}  }
\end{equation}
is an isomorphism.
\end{thm}
The idea is that using creation and annihilation operators in \cite{vincent} and 
%fact that $\restr{\mc H_{G, N, I \sqcup I \sqcup I, W \boxtimes W^* \boxtimes W}^j}{ (\ov{\eta})^{I \sqcup I} \times_{\on{Spec} \Fqbar} \ov x }$
%the restriction of $\mc H_{G, N, I, W}^j$ to the scheme $(\ov{\eta})^{I_1} \times_{\on{Spec} \Fqbar} \ov u$ is constant, 
%the results in Section 2, 
Proposition \ref{prop-H-ov-eta-I-1-ov-s-I-2-is-constant},
we construct a morphism $\Upsilon: \restr{ \mc H_{G, N, I, W}^j  }{ \ov{\eta_I}  } \rightarrow \restr{ \mc H_{G, N, I, W}^j  }{ \ov{x}  }$. Then using the "Zorro" lemma, we prove that the composition $\Upsilon \circ \mf{sp}_{\ov x}^*$ and $\mf{sp}_{\ov x}^* \circ \Upsilon$ are isomorphisms.

\quad

As an immediate consequence, in Section 5, we prove
\begin{prop} \label{prop-Weil-eta-factos-through-Weil-X-intro} (Proposition \ref{prop-Weil-eta-factors-through-Weil-X})
The action of $\on{Weil}(\eta, \ov{\eta})^I$ on $\restr{ \mc H_{G, N, I, W}^j }{ \ov{\eta_I} }$ factors through $\on{Weil}(X \sm N, \ov{\eta})^I$.
\end{prop}

\quad

In Section 6, we treat the case where $G$ is not necessarily split. The results in Sections 1-5 still hold.

In Section 7, with the help of Theorem \ref{thm-H-I-W-smooth-intro}, we extend for split $G$ the constant term morphism constructed in \cite{coho-cusp} from $\eta^I$ to $(X \sm N)^I$. Then we define a smooth cuspidal cohomology subsheaf $\mc H_{G, N, I, W}^{j, \, \on{cusp}}$ of $\mc H_{G, N, I, W}^j$ over $(X \sm N)^I$. 

\quad

For the cohomology sheaves with integral coefficients, we still have the main results as in Sections 1-7. For this, we will need the context in \cite{xiao-zhu}.
We will write these in a future version.

%Note that the proof of Theorem \ref{thm-H-I-W-smooth} does not use the compactification of stacks of shtukas.

\quad

Proposition \ref{prop-Weil-eta-factos-through-Weil-X-intro} was already proved by Xinwen Zhu in the MSRI 2019 hot topic workshop. But Theorem \ref{thm-H-I-W-smooth-intro} is new. 

\quad 

\subsection*{Acknowledgments}
I would like to thank Vincent Lafforgue, Gérard Laumon and Jack Thorne for stimulating discussions.

%\subsection*{Notations}

\quad

\section{The cohomology sheaves are constant over $(\ov{\eta})^I$}

In this section, we first recall some general results about ind-smooth sheaves in \ref{subsection-ind-smooth-sheaf} and partial Frobenius morphisms in \ref{recall-partial-Frob}. Then we recall the cohomology sheaves of stacks of shtukas, prove Proposition \ref{prop-FWeil-factors-though-Weil-I} and Proposition \ref{prop-H-x-H-eta-I-isom}. As a result, we deduce Proposition \ref{prop-H-I-W-constant-over-ov-eta-I}.

\subsection{Reminders on ind-smooth $E$-sheaves}   \label{subsection-ind-smooth-sheaf}

\sssec{}
We use the definition in \cite{sga5} VI 1 and \cite{Weil2} 1.1 for constructible $E$-sheaves (for étale topology) and smooth constructible $E$-sheaves. In this paper, a smooth $E$-sheaf always means a smooth constructible $E$-sheaf.

We define an ind-constructible $E$-sheaf to be an inductive limit of constructible $E$-sheaves (i.e. an object in the abstract category of inductive limits of constructible $E$-sheaves). Its fiber at a geometric point is defined to be the $E$-vector space which is the inductive limit of fibers (i.e. in the category of $E$-vector spaces).

\sssec{}
We use \cite{sga4} VIII 7 for the definition of specialization maps.

\sssec{}     \label{subsection-smooth-equivalent-to-sp-bij}
We denote by $\mc O_E$ the ring of integers of $E$ and $\lambda_E$ a uniformizer.
Let $\Lambda = \mc O_E / \lambda_E^s \mc O_E$ for $s \in \N$.

Let $Y$ be a noetherian scheme over $\Fq$.
%Let $\Lambda$ be a noetherian ring.
Let $\mc F$ be a constructible $\Lambda$-sheaf over $Y$.
By \cite{sga4} IX Proposition 2.11, the sheaf $\mc F$ is locally constant if and only if for any geometric points $\ov x$, $\ov y$ of $Y$ and any specialization map $\mf{sp}: \ov y \rightarrow \ov x$, the induced morphism $$\mf{sp}^*: \restr{  \mc F  }{ \ov x }  \rightarrow \restr{  \mc F  }{ \ov y  } $$ is an isomorphism.

Using the definition of constructible $\mc O_E$-sheaves (resp. constructible $E$-sheaves), we deduce that a constructible $\mc O_E$-sheaf (resp. constructible $E$-sheaf) $\mc F$ over $Y$ is smooth if and only if for any geometric points $\ov x$, $\ov y$ of $Y$ and any specialization map $\mf{sp}: \ov y \rightarrow \ov x$, the induced morphism $$\mf{sp}^*: \restr{  \mc F  }{ \ov x }  \rightarrow \restr{  \mc F  }{ \ov y  } $$ is an isomorphism. 

%Note that an $E$-sheaf is isomorphic to a torsion free $\lambda$-adic sheaf in the category of $E$-sheaf. 
%Using the definition of constructible $E$-sheaves, we deduce that a constructible $E$-sheaf $\mc F$ over $Y$ is smooth if and only if for any geometric points $\ov x$, $\ov y$ of $Y$ and any specialization map $\mf{sp}: \ov y \rightarrow \ov x$, the induced morphism $$\mf{sp}^*: \restr{  \mc F  }{ \ov x }  \rightarrow \restr{  \mc F  }{ \ov y  } $$ is an isomorphism. 

\sssec{}    \label{subsection-def-ind-smooth}
Let $Y$ be a normal connected noetherian scheme over $\Fq$.
Let $\mc H = \varinjlim_{\lambda \in \Omega} \mc F_{\lambda}$ be an inductive limit of constructible $E$-sheaves over $Y$, where $\Omega$ is a numerable filtered set. We say that the ind-constructible $E$-sheaf $\mc H$ is {\it{ind-smooth}} if we can write $\mc H$ as an inductive limit of smooth $E$-sheaves over $Y$, i.e. there exists a numerable filtered set $\Omega'$ and smooth $E$-sheaves $\mc G_{\lambda}$ for $\lambda \in \Omega'$ such that $\mc H = \varinjlim_{\lambda \in \Omega'} \mc G_{\lambda}$.

\begin{lem}   \label{lem-ind-smooth-specialization}
For $Y$ as in \ref{subsection-def-ind-smooth}, an ind-constructible $E$-sheaf $\mc H$ over $Y$ is ind-smooth if and only if for any geometric points $\ov x$, $\ov y$ of $Y$ and any specialization map $\mf{sp}: \ov y \rightarrow \ov x$, the induced morphism
$$\mf{sp}^*: \restr{  \mc H  }{ \ov x }  \rightarrow \restr{  \mc H  }{ \ov y  } $$ is an isomorphism. 
\end{lem}

\sssec{}   \label{subsection-def-wt-F}
To prove Lemma \ref{lem-ind-smooth-specialization}, we need some preparations. Let $\mc H = \varinjlim_{\lambda \in \Omega} \mc F_{\lambda}$ as above. For any $\lambda \leq \mu$ in $\Omega$, the kernel $\on{Ker}(\mc F_{\lambda} \rightarrow \mc F_{\mu})$ is a constructible sub-$E$-sheaf of $\mc F_{\lambda}$. For $\lambda \leq \mu_1 \leq \mu_2$, we have $$\on{Ker}(\mc F_{\lambda} \rightarrow \mc F_{\mu_1}) \subset \on{Ker}(\mc F_{\lambda} \rightarrow \mc F_{\mu_2}) \subset \on{Ker}(\mc F_{\lambda} \rightarrow \mc H) \subset \mc F_{\lambda}.$$
Since $\mc F_{\lambda}$ is constructible and $Y$ is noetherian, we deduce that there exists $\lambda_0$, such that for all $\mu \geq \lambda_0$, we have $\on{Ker}(\mc F_{\lambda} \rightarrow \mc F_{\lambda_0}) = \on{Ker}(\mc F_{\lambda} \rightarrow \mc F_{\mu}). $ (The argument is similar to the proof of Lemma 58.73.2 of \cite{stacks-project}.)
%For every $\mu$, let $T_{\mu} \subset Y$ be the set of points $y \in Y$ such that $\restr{\on{Ker}(\mc F_{\lambda} \rightarrow \mc F_{\mu})}{\ov y}  \rightarrow \restr{\on{Ker}(\mc F_{\lambda} \rightarrow \varinjlim \mc F_{\mu})}{\ov y}$ is not surjective. Since $\on{Ker}(\mc F_{\lambda} \rightarrow \mc F_{\mu})$ and  $\mc F_{\lambda}$ are constructible, the set $T_{\mu}$ is a constructible subset of $Y$. Since $\restr{\mc F_{\lambda}}{\ov y}$ is a finite dimensional $E$-vector space, we see that for all $y \in Y$ we have $y \notin T_{\mu}$ for $\mu$ large enough. Since $Y$ is noetherian, we deduce that $T_{\mu} = \emptyset$ for $\mu$ large enough.)
So $\on{Im}(\mc F_{\lambda} \rightarrow \mc F_{\lambda_0}) \isom \on{Im}(\mc F_{\lambda} \rightarrow \mc F_{\mu}).$
We denote by $\wt{\mc F_{\lambda}}:= \on{Im}(\mc F_{\lambda} \rightarrow \mc F_{\lambda_0})$. 
%Since $$\mc F_{\lambda} \rightarrow \wt{\mc F_{\lambda}} \rightarrow \mc F_{\lambda_0}$$
We have 
\begin{equation}   \label{equation-lim-wt-F-equal-lim-F}
\varinjlim_{\lambda \in \Omega} \mc F_{\lambda} = \varinjlim_{\lambda \in \Omega} \wt{ \mc F_{\lambda} }.
\end{equation}

\noindent {\bf Proof of Lemma \ref{lem-ind-smooth-specialization}:}
One direction is obvious.
Let's prove the converse direction. 
%Let $\mc H = \varinjlim \mc F_{\lambda}$ be an inductive limit of constructible $E$-sheaves over $Y$. 
Suppose that all of the specialization homomorphisms $\mf{sp}^*: \restr{  \mc H  }{ \ov x }  \rightarrow \restr{  \mc H  }{ \ov y  } $ are isomorphisms.

%Note that for any $\lambda$, $\restr{\mc F_{\lambda}}{\eta_Y}$ is an $E$-vector space. So there exists $\lambda_0$ large enough such that $\on{Ker}(\restr{\mc F_{\lambda}}{\eta_Y} \rightarrow \restr{\mc F_{\lambda_0}}{\eta_Y} ) = \on{Ker}(\restr{\mc F_{\lambda}}{\eta_Y} \rightarrow \varinjlim \restr{\mc F_{\lambda'} }{\eta_Y} ) $. Replacing $\mc F_{\lambda}$ by its image in $\mc F_{\lambda_0}$, we can suppose that for all $\lambda_1 \leq \lambda_2$, the morphism $\restr{\mc F_{\lambda_1}}{\eta_Y} \rightarrow \restr{\mc F_{\lambda_2}}{\eta_Y} $ is injective.
Since every $\wt{\mc F_{\lambda}}$ is a constructible $E$-sheaf over $Y$, there exists an open dense subscheme $U_{\lambda}$ of $Y$ such that $\wt{\mc F_{\lambda}}$ is smooth over $U_{\lambda}$. Let $j_{\lambda}: U_{\lambda} \hookrightarrow Y$ be the embedding. Let $$\mc G_{\lambda}:= (j_{\lambda})_* ( \restr{ \wt{ \mc F_{\lambda} } }{ U_{\lambda} }  ).$$ To prove that $\mc H$ is ind-smooth, it is enough to prove that

(1) every $\mc G_{\lambda}$ is a smooth $E$-sheaf over $Y$

(2) $\varinjlim_{\lambda \in \Omega} \wt{ \mc F_{\lambda} } = \varinjlim_{\lambda \in \Omega} \mc G_{\lambda}$

\quad

First we prove (1): the adjunction morphism $\Id \rightarrow (j_{\lambda})_*(j_{\lambda})^*$ induces a morphism 
\begin{equation}   \label{equation-wt-F-to-G}
\wt{ \mc F_{\lambda} } \rightarrow (j_{\lambda})_* (j_{\lambda})^* \wt{ \mc F_{\lambda} } = \mc G_{\lambda}
\end{equation}
Taking limit, we deduce a morphism $\varinjlim \wt{ \mc F_{\lambda} } \rightarrow \varinjlim \mc G_{\lambda}$. 
Let $\eta_Y$ be the generic point of $Y$ and $\ov{\eta_Y}$ a geometric point over $\eta_Y$. 
For any geometric point $\ov y$ of $Y$ and any specialization map $\mf{sp}_{\ov y}: \ov{\eta_Y} \rightarrow \ov y$,
we have a commutative diagram
\begin{equation}    \label{equation-lim-F-lim-G}
\xymatrix{
\restr{\varinjlim \wt{ \mc F_{\lambda} } }{\ov y}  \ar[d]  \ar[r]^{\mf{sp}_{\ov y}^*}_{\simeq}
& \restr{\varinjlim \wt{ \mc F_{\lambda} } }{\ov{\eta_Y}} \ar[d]^{\simeq}  \\
\restr{\varinjlim \mc G_{\lambda} }{\ov y}      \ar[r]^{\mf{sp}_{\ov y}^*}
& \restr{\varinjlim \mc G_{\lambda}}{\ov{\eta_Y}} 
}
\end{equation}
By the hypothesis and (\ref{equation-lim-wt-F-equal-lim-F}), the upper line of (\ref{equation-lim-F-lim-G}) is an isomorphism. By the definition of $\mc G_{\lambda}$, the right vertical line of (\ref{equation-lim-F-lim-G}) is an isomorphism. 
Thus the lower line of (\ref{equation-lim-F-lim-G}) is surjective.

Note that for every $\lambda$, the morphism $\mf{sp}_{\ov y}^*: \restr{\mc G_{\lambda} }{\ov y}    \rightarrow  
 \restr{ \mc G_{\lambda}}{\ov{\eta_Y}} $ is injective. In fact, by definition we have $$\restr{\mc G_{\lambda} }{\ov y}   = \Gamma(Y_{(\ov y)}, \mc G_{\lambda}) = \Gamma(Y_{(\ov y)}, (j_{\lambda})_* ( \restr{ \wt{\mc F_{\lambda}} }{ U_{\lambda} } ) ) =
 \Gamma(Y_{(\ov y)} \times_Y U_{\lambda}, \restr{ \wt{ \mc F_{\lambda} } }{U_{\lambda}})$$
 where $Y_{(\ov y)}$ is the strict henselization of $Y$ at $\ov y$.
By \cite{sga1} I Proposition 10.1, since $Y$ is normal connected, the fiber product $Y_{(\ov y)} \times_Y U_{\lambda}$ is connected. Since $\restr{ \wt{ \mc F_{\lambda} } }{U_{\lambda}}$ is smooth, the restriction
$$\Gamma(Y_{(\ov y)} \times_Y U_{\lambda}, \restr{ \wt{ \mc F_{\lambda} } }{U_{\lambda}}) \rightarrow \restr{ \wt{ \mc F_{\lambda}} }{\ov{\eta_Y}}$$ is injective.
 
Now we want to prove that for any $\ov y$ and any $\mf{sp}_{\ov y}$, the induced morphism $\mf{sp}_{\ov y}^*: \restr{\mc G_{\lambda} }{\ov y}   \rightarrow \restr{\mc G_{\lambda}}{\ov{\eta_Y}} $ is surjective.
Let $a \in \restr{ \mc G_{\lambda} }{\ov{\eta_Y}} $.
Since the lower line of (\ref{equation-lim-F-lim-G}) is surjective, there exists $\mu \geq \lambda$ and $b \in \restr{\mc G_{\mu}}{\ov{y}}$ such that the image of $b$ by
$\restr{\mc G_{\mu}}{\ov{y}} \rightarrow \restr{\varinjlim \mc G_{\lambda'} }{\ov y}   \rightarrow  \restr{\varinjlim \mc G_{\lambda'}}{\ov{\eta_Y}} $ coincides with the image of $a$ in $\restr{\varinjlim \mc G_{\lambda'}}{\ov{\eta_Y}} $. Note that $$\restr{\mc G_{\mu}}{\ov{y}} = \Gamma(Y_{(\ov{y})} , \mc G_{\mu}) = \Gamma(Y_{(\ov{y})} \times_Y U_{\mu} , \wt{\mc F_{\mu}}) = \Gamma(Y_{(\ov{y})} \times_Y \eta_Y , \wt{\mc F_{\mu}})$$
where the last equality is because that $\restr{ \wt{\mc F_{\mu}} }{U_{\mu}}$ is smooth.
We have a commutative diagram
$$
\xymatrix{
& \Gamma(Y_{(\ov{y})} \times_Y \eta_Y , \wt{\mc F_{\mu}}) \ar[ld]_{=}  \ar@{_{(}->}[rd]^{\text{rest}}   \\
\restr{\mc G_{\mu}}{\ov{y}}  \ar[rr]^{\mf{sp}_{\ov y}^*}
& & \restr{\mc G_{\mu}}{\ov{\eta_Y}} = \restr{ \wt{\mc F_{\mu}} }{\ov{\eta_Y}} 
}
$$
%$$
%\xymatrix{
%& \Gamma(Y_{(\ov{y})} \times_Y \eta_Y , \mc F_{\mu}) \ar[ld]_{=}  \ar@{_{(}->}[rd]^{\text{res}}  
%& \Gamma(Y_{(\ov{y})} \times_Y \eta_Y , \mc F_{\lambda}) \ar@{_{(}->}[l] \ar@{_{(}->}[rd]^{\text{res}}  \\
%\restr{\mc G_{\mu}}{\ov{y}}  \ar[rr]^{\mf{sp}_{\ov y}^*}
%& & \restr{\mc G_{\mu}}{\ov{\eta_Y}}
%& \restr{\mc G_{\lambda}}{\ov{\eta_Y}}  \ar@{_{(}->}[l]
%}
%$$
As before, since $Y$ is normal connected, the fiber product $Y_{(\ov{y})} \times_Y \eta_Y$ is connected. Over $Y_{(\ov{y})} \times_Y \eta_Y$, $\wt{\mc F_{\lambda}}$ is a smooth subsheaf of the smooth sheaf $\wt{\mc F_{\mu}}$.
For $b \in \Gamma(Y_{(\ov{y})} \times_Y \eta_Y , \wt{\mc F_{\mu}}) $, if $\text{rest}(b)$ is in $\restr{ \wt{\mc F_{\lambda}} }{\ov{\eta_Y}}$, then $b \in \Gamma(Y_{(\ov{y})} \times_Y \eta_Y , \wt{\mc F_{\lambda}}) = \restr{\mc G_{\lambda} }{\ov y}$. We deduce that $\mf{sp}_{\ov y}^*(b) = a$. Thus
 $\mf{sp}_{\ov y}^*: \restr{\mc G_{\lambda} }{\ov y}   \rightarrow \restr{\mc G_{\lambda}}{\ov{\eta_Y}} $ is surjective.

By \ref{subsection-smooth-equivalent-to-sp-bij}, we deduce that $\mc G_{\lambda}$ is a smooth $E$-sheaf over $Y$.

\quad
 
Now we prove (2): by definition, for all $\lambda_1 \leq \lambda_2$, the morphism 
%$\restr{\wt{\mc F_{\lambda_1}}}{\eta_Y} \rightarrow \restr{\wt{\mc F_{\lambda_2}}}{\eta_Y} $ 
$\wt{\mc F_{\lambda_1}} \rightarrow \wt{\mc F_{\lambda_2}}$ 
is injective. So the morphism $\restr{\mc G_{\lambda_1}}{\eta_Y} \rightarrow \restr{\mc G_{\lambda_2}}{\eta_Y} $ is injective. Since $\mc G_{\lambda_1}$ and $\mc G_{\lambda_2}$ are smooth, we deduce that the morphism $\mc G_{\lambda_1} \rightarrow \mc G_{\lambda_2}$ is injective.

We proved that for any $\lambda$, the morphism $\mf{sp}_{\ov y}^*: \restr{\mc G_{\lambda} }{\ov y}    \rightarrow  
 \restr{ \mc G_{\lambda}}{\ov{\eta_Y}} $ is injective.
We deduce that the lower line of (\ref{equation-lim-F-lim-G}) is injective.
%Since $\mc G_{\lambda}$ is smooth, the lower line is an isomorphism. 
So the lower line of (\ref{equation-lim-F-lim-G}), thus the left vertical line, is an isomorphism. 
%Since this is true for any $\ov y$, {\color{blue}How??? }we deduce that the morphism $\varinjlim \mc F_{\lambda} \rightarrow \varinjlim \mc G_{\lambda}$ is an isomorphism. Thus $\mc H = \varinjlim \mc G_{\lambda}$.

For any $\lambda$, we have a commutative diagram
\begin{equation}    
\xymatrix{
\restr{ \wt{ \mc F_{\lambda} } }{\ov y}  \ar[d]  \ar@{^{(}->}[r]
& \restr{\varinjlim \wt{ \mc F_{\lambda} } }{\ov y} \ar[d]^{\simeq}  \\
\restr{ \mc G_{\lambda} }{\ov y}      \ar@{^{(}->}[r]
& \restr{\varinjlim \mc G_{\lambda} }{\ov y} 
}
\end{equation}
%where the lower line is injective because that $\mc G_{\lambda}$ are smooth and that $\restr{ \mc G_{\lambda} }{\ov{\eta_Y}}  \rightarrow \restr{\varinjlim \mc G_{\lambda} }{\ov{\eta_Y}} $ is injective.
We deduce that $\restr{ \wt{ \mc F_{\lambda} } }{\ov y} \rightarrow \restr{ \mc G_{\lambda} }{\ov y}  $ is injective. Since this is true for any $\ov y$, we deduce that (\ref{equation-wt-F-to-G}) is injective.

Now fix $\lambda$. For any $\mu$, consider the subset of $Y$
$$C_{\mu}:=\{ y \in Y \text{ such that }  \restr{ \mc G_{\lambda} }{\ov y} \nsubseteq \on{Im}(  \restr{ \wt{ \mc F_{\mu} } }{\ov y} \rightarrow \restr{ \mc G_{\mu} }{\ov y} )   \}$$
It is constructible. For any $\mu_1 \leq \mu_2$, we have $C_{\mu_1} \supset C_{\mu_2}$. We have $\cap_{\mu} C_{\mu} = \emptyset$. We deduce that there exists $\mu'(\lambda)$, such that for any $\mu \geq \mu'(\lambda)$, we have $C_{\mu} = \emptyset$. In particular, $ \mc G_{\lambda}  \subset \on{Im}(   \wt{ \mc F_{\mu'(\lambda)} }  \rightarrow  \mc G_{\mu'(\lambda)} ) $. Thus for any $\lambda$, we have $$\wt{\mc F_{\lambda}} \subset \mc G_{\lambda} \subset \wt{\mc F_{\mu'(\lambda)}}$$
This implies 
\begin{equation}    \label{equation-lim-wt-F-equal-lim-G}
\varinjlim_{\lambda \in \Omega} \wt{ \mc F_{\lambda} } = \varinjlim_{\lambda \in \Omega} \mc G_{\lambda} .
\end{equation}

\quad

%Combining (\ref{equation-lim-wt-F-equal-lim-F}) and (\ref{equation-lim-wt-F-equal-lim-G})

%{\color{blue}New idea: let $\wt{\mc F}_{\lambda}$ be the image of $\restr{\mc F_{\lambda}}{\eta_Y}$ in $\restr{\varinjlim \mc F_{\lambda}}{\eta_Y}$. By the hypothesis the action of all inertia groups $\ms I_v$ on $\restr{ \wt{\mc F}_{\lambda} }{\ov{\eta_Y}}$ is trivial. So $\mc G_{\lambda}:= j_* ( \wt{\mc F}_{\lambda} )$ is a smooth sheaf. Then the above argument works.  }

\sssec{}    \label{subsection-smooth-sheaf-pi-1-rep}
Let $Y$ be a normal connected noetherian scheme. Let $\eta_Y$ be the generic point of $Y$ and $\ov{\eta_Y}$ a geometric point over $\eta_Y$. Then the functor $\mc F \mapsto \restr{\mc F}{\ov{\eta_Y}}$ from
the category of smooth $E$-sheaves over $Y$ to the category of finite dimensional $E$-vector spaces with continuous $\pi_1(Y, \ov{\eta_Y})$-actions is an equivalence. (cf. \cite{sga5} VI 1.)

Let $\mc H = \varinjlim \mc F_{\lambda}$ be an inductive limit of smooth $E$-sheaves $\mc F_{\lambda}$ over $Y$. Replacing $\mc F_{\lambda}$ by $\wt{ \mc F_{\lambda} }$ as in \ref{subsection-def-wt-F}, we can suppose that for all $\lambda_1 \leq \lambda_2$, the morphism $\restr{\mc F_{\lambda_1}}{\eta_Y} \rightarrow \restr{\mc F_{\lambda_2}}{\eta_Y} $ is injective.

Suppose that the action of $\pi_1(Y, \ov{\eta_Y})$ on $\restr{\mc H}{\ov{\eta_Y}}$ is trivial. Then for every $\lambda$ the action of $\pi_1(Y, \ov{\eta_Y})$ on $\restr{\mc F_{\lambda}}{\ov{\eta_Y}}$ is trivial. Since $\mc F_{\lambda}$ is a smooth $E$-sheaf, by the above equivalence every $\mc F_{\lambda}$ is a constant sheaf. Thus $\mc H$ is a constant sheaf.

\begin{rem}
In this paper, we use the étale topology.
One can also use the pro-étale topology instead of the étale topology. In the context of pro-étale topology, Lemma \ref{lem-ind-smooth-specialization} would be more obvious. %We will write this in a future version.
\end{rem}

\subsection{Reminders on partial Frobenius morphisms}   \label{recall-partial-Frob}

\sssec{}
In this paper, $\times_{\Fq}$ means $\times_{\on{Spec}\Fq}$ and 
$\times_{\Fqbar}$ means $\times_{\on{Spec}\Fqbar}$.

\sssec{}
We denoted by $F$ the function field of $X$. Fix an algebraic closure $\ov F$ of $F$ and an embedding $\Fqbar \subset \ov F$. We denote by $\eta$ the generic point of $X$ and by $\ov{\eta}$ the geometric point over $\eta$.

Let $I$ be a finite set. We denote by $X^I:=X \times_{\Fq} \cdots \times_{\Fq} X$ ($I$ copies).
We denote by $F_I$ the function field of $X^I$. Fix an algebraic closure $\ov {F_I}$ of $F_I$. We denote by $\eta_I$ the generic point of $X^I$ and by $\ov{\eta_I}$ the geometric point over $\eta_I$.

Our notation is slightly different from \cite{vincent}, where the function field of $X^I$ is denoted by $F^I$ and the generic point of $X^I$ is denoted by $\eta^I$. We will never use the notations $F^I$ and $\eta^I$ in this paper.

\begin{notation}
We denote by
$$(\eta)^I:= \eta \times_{\Fq}  \cdots \times_{\Fq} \eta  \quad \text{ and } \quad  (\ov{\eta})^I:= \ov{\eta} \times_{\Fqbar}  \cdots \times_{\Fqbar} \ov{\eta}.$$
\end{notation}

\sssec{}   \label{subsection-ov-eta-I-integral-scheme}
Note that $(\ov{\eta})^I$ is an integral scheme over $\on{Spec} \Fqbar$. In fact, we have
$$(\ov{\eta})^I = \varprojlim Y \otimes_{\Fqbar} \cdots \otimes_{\Fqbar} Y $$
where the projective limit is over affine étale $Y=\on{Spec} A$ over $X_{\Fqbar}$. By hypothesis $X_{\Fqbar}$ is irreducible. 
Every $A$ is an integral domain which is finitely generated as $\Fqbar$-algebra, 
%Since $Y$ is an irreducible finite type scheme over $\Fqbar$, 
thus the product $A \otimes_{\Fqbar} \cdots \otimes_{\Fqbar} A$ is still an integral domain (see for example \cite{stacks-project} Part 2, Lemma 33.3.3). We deduce that $\ov F \otimes_{\Fqbar} \cdots \otimes_{\Fqbar} \ov F$ is an integral domain.
%$$\ov F \otimes_{\Fqbar} \cdots \otimes_{\Fqbar} \ov F = \varinjlim A \otimes_{\Fqbar} \cdots \otimes_{\Fqbar} A $$ where the injective limit is over $A$ integral domain finitely generated $\Fqbar$-algebra.

As a consequence, $\ov{\eta_I}$ is a geometric generic point of $(\ov{\eta})^I$.

\sssec{}   \label{subsection-def-Frob-J}
For any scheme $Y$ over $\Fq$, we denote by $\Frob: Y \rightarrow Y$ the Frobenius morphism over $\Fq$.

For any subset $J \subset I$, we denote by $$\on{Frob}_{J}: X^I \rightarrow X^I$$ the morphism sending $(x_i)_{i \in I}$ to $(x'_i)_{i \in I}$, with $x'_i = \on{Frob}(x_i)$ if $i \in J$ and $x_i' = x_i$ if $i \notin J$.

In particular, for $J = \{i\}$ a singleton, we have the morphism $$\on{Frob}_{\{i\}}: X^I \rightarrow X^I.$$
%be the morphism sending $(x_j)_{j \in I}$ to $(x'_j)_{j \in I}$, with $x'_i = \on{Frob}(x_i)$ and $x_j' = x_j$ if $j \neq i$.

\sssec{}  \label{subsection-def-FWeil-eta-I}
As in \cite{vincent} Remarque 8.18, we define
$$\on{FWeil}(\eta_I, \ov{\eta_I}):=\{ \varepsilon \in \on{Aut}_{\ov{\Fq}}(\ov{F_I}) \; | \; \exists (n_i)_{i \in I} \in \Z^I, \restr{\varepsilon}{(F_I)^{\on{perf}}} = \prod_{i \in I} (\on{Frob}_{\{i\}})^{n_i}  \} .$$
where $(F_I)^{\on{perf}}$ is the perfection of $F_I$. %It is an extension of $\Z^I$ by $\pi_1^{\on{geom}}(\eta_I, \ov{\eta_I})$.

\sssec{} (\cite{vincent} Remarque 8.18)
Let $\Delta: X \rightarrow X^I$ be the diagonal inclusion.
We fix a specialization map in $X^I$
\begin{equation*}
\mf{sp}: \ov{\eta_I} \rightarrow \Delta(\ov{\eta})
\end{equation*}
i.e. a morphism from $\ov{\eta_I}$ to the strict henselization of $X^I$ at $\Delta(\ov{\eta})$. In particular, for every étale neighbourhood $U$ of $\Delta(\ov{\eta})$ in $X^I$, the specialization map $\mf{sp}$ induces a $X^I$-morphism $\ov{\eta_I} \rightarrow U$. The scheme $(\ov{\eta})^I$ equiped with the diagonal morphism $\Delta(\ov{\eta}) \rightarrow (\ov{\eta})^I$ is a projective limit of étale neighbourhoods of $\Delta(\ov{\eta})$ in $X^I$ of the form $V^I$ with $V$ étale neighbourhood of $\ov{\eta}$ in $X$. We deduce a morphism
$$\ov{\eta_I} \rightarrow (\ov{\eta})^I,$$
i.e. an inclusion
%As in \cite{vincent} Remarque 8.18, the specialization map $\mf{sp}$ induces an inclusion 
\begin{equation}   \label{equation-ov-F-times-ov-F-includ-ov-F-I}
\ov{F} \otimes_{\ov{\Fq}} \cdots \otimes_{\ov{\Fq}} \ov F \subset \ov{F_I} .
\end{equation}

Let $\varepsilon \in \on{FWeil}(\eta_I, \ov{\eta_I}) $.
For any $i \in I$, let $\varepsilon_i$ be the restriction of $\varepsilon$ to the algebraic closure in $\ov{F_I}$ of
$\Fq \otimes_{\Fq} \cdots  \otimes_{\Fq} F \otimes_{\Fq}  \cdots \otimes_{\Fq} \Fq$, where $F$ is the $i$-th factor.
We have a surjective morphism
$$
\begin{aligned}
\Psi: \on{FWeil}(\eta_I, \ov{\eta_I}) & \rightarrow \on{Weil}(\eta, \ov{\eta})^I  \\
\varepsilon \quad & \mapsto \big( (\on{Frob}_{\{i\}})^{-n_i} \circ \varepsilon_i \big)_{i \in I}
\end{aligned}
$$
%Its kernel is equal to $\on{Ker}\big( \pi_1^{\on{geom}}(\eta_I, \ov{\eta_I}) \twoheadrightarrow \pi_1^{\on{geom}}(\eta, \ov{\eta})^I  \big)$.
We have exact sequences:
$$
\xymatrix{
0 \ar[r]
& \pi_1^{\on{geom}}(\eta^I, \ov{\eta^I}) \ar[r]  \ar@{->>}[d]
& \on{FWeil}(\eta^I, \ov{\eta^I}) \ar[r]  \ar@{->>}[d]^{\Psi}
& \Z^I \ar[r]   \ar[d]^{=}
& 0 \\
0 \ar[r]
& \pi_1^{\on{geom}}(\eta, \ov{\eta})^I \ar[r]
& \on{Weil}(\eta, \ov{\eta})^I \ar[r]
& \Z^I \ar[r]
& 0
}
$$

\sssec{}   \label{subsection-def-FWeil-action}
Let $\mc G$ be an ind-constructible $E$-sheaf over $(\eta)^I$, equiped with an action of the partial Frobenius morphisms, i.e. equiped with isomorphisms $F_{\{i\}}: \on{Frob}_{\{i\}}^* \mc G \isom \mc G$ commuting to each other and whose composition is the total Frobenius isomorphism $\on{Frob}^* \mc G \isom \mc G$ over $(\eta)^I$.

Then $\restr{ \mc G }{ \ov{\eta_I} }$ is equiped with an action of $\on{FWeil}(\eta_I, \ov{\eta_I})$ in the following way:
for any $\varepsilon \in \on{FWeil}(\eta_I, \ov{\eta_I})$ with $\restr{\varepsilon}{(F_I)^{\on{perf}}} = \prod_{i \in I} (\on{Frob}_{\{i\}})^{n_i}$, it induces a commutative diagram
$$
\xymatrixrowsep{2pc}
\xymatrixcolsep{6pc}
\xymatrix{
\ov{F_I}  \ar[d]
& \ov{F_I} \ar[l]_{\varepsilon}  \ar[d]   \\
(F_I)^{\on{perf}}
& (F_I)^{\on{perf}}  \ar[l]_{\prod_{i \in I} (\on{Frob}_{\{i\}})^{n_i}}
}
$$
In other words, a commutative diagram
\begin{equation}   \label{equation-spec-varepsilon-over-F-I-perf}
\xymatrixrowsep{2pc}
\xymatrixcolsep{6pc}
\xymatrix{
\ov{\eta_I}  \ar[d]    \ar[r]^{\on{Spec} \varepsilon}  
& \ov{\eta_I}  \ar[d]   \\
\on{Spec}  (F_I)^{\on{perf}} \ar[r]^{\prod_{i \in I} (\on{Frob}_{\{i\}})^{n_i}}
& \on{Spec} (F_I)^{\on{perf}} 
}
\end{equation}
Moreover, we have a commutative diagram
\begin{equation}    \label{equation-spec-varepsilon-over-Fqbar}
\xymatrix{
\ov{\eta_I}  \ar[dr]    \ar[rr]^{\on{Spec} \varepsilon} 
& & \ov{\eta_I}  \ar[dl]   \\
& \on{Spec}  \Fqbar
}
\end{equation}
By (\ref{equation-spec-varepsilon-over-F-I-perf}), we deduce a specialization map in $X^I$:
$$\mf{sp}_{\varepsilon}: \prod_{i \in I} \on{Frob}_{\{i\}}^{n_i}(\ov{\eta_I}) \rightarrow \ov{\eta_I}$$
which is in fact an isomorphism of schemes over $\on{Spec} (F_I)^{\on{perf}}$ and over $\on{Spec} \Fqbar$.
%It induces a morphism $\mf{sp}_{\varepsilon}^*: \restr{\mc G}{\ov{\eta_I}} \rightarrow \restr{\mc G}{\prod_{i \in I} (\on{Frob}_{\{i\}})^{n_i}(\ov{\eta_I}) }$. Besides, 

The action of $\varepsilon$ on $\restr{ \mc G }{ \ov{\eta_I} }$ is defined to be the composition:
\begin{equation}
\restr{\mc G}{\ov{\eta_I}} \xrightarrow{\mf{sp}_{\varepsilon}^*} \restr{\mc G}{\prod_{i \in I} \on{Frob}_{\{i\}}^{n_i}(\ov{\eta_I}) } 
 = \restr{\big( \prod_{i \in I} (\on{Frob}_{\{i\}}^{n_i})^* \mc G \big)}{ \ov{\eta_I} }
\xrightarrow{ \prod_{i \in I} F_{\{i\}}^{n_i} } \restr{\mc G}{\ov{\eta_I}}
\end{equation}
We deduce an action of $\on{FWeil}(\eta_I, \ov{\eta_I})$.

\sssec{}
An action of $\on{FWeil}(\eta_I, \ov{\eta_I})$ on a finite dimensional $E$-vector space is said to be continuous if the action of $\pi_1^{\on{geom}}(\eta_I, \ov{\eta_I})$ is continuous. 

More generally, an action of $\on{FWeil}(\eta_I, \ov{\eta_I})$ on an $E$-vector space $M$ is said to be continuous if $M$ is an inductive limit of finite dimensional $E$-vector subspaces which are stable under $\pi_1^{\on{geom}}(\eta_I, \ov{\eta_I})$ and on which the action of $\pi_1^{\on{geom}}(\eta_I, \ov{\eta_I})$ is continuous.

\quad

We will need the following two variants of a lemma of Drinfeld (these variants are proved by Drinfeld and recalled in \cite{coho-filt}).

\begin{lem}  \label{lem-Drinfeld-E-vs}  (\cite{coho-filt} Lemma 3.2.10)
A continuous action of $\on{FWeil}(\eta_I, \ov{\eta_I})$ on a finite dimensional $E$-vector space factors through $\on{Weil}(\eta, \ov{\eta})^I$.
\cqfd
\end{lem}

\begin{lem}  \label{lem-Drinfeld-A-mod}  (\cite{coho-filt} Lemma 3.2.13)
Let $A$ be a finitely generated $E$-algebra. Let $M$ be an $A$-module of finite type. Then a continuous $A$-linear action of $\on{FWeil}(\eta_I, \ov{\eta_I})$ on $M$ factors through $\on{Weil}(\eta, \ov{\eta})^I$. %where an action of $\on{FWeil}(\eta_I, \ov{\eta_I})$ on $M$ is said to be continuous if 
\cqfd
\end{lem}

We will also need

\begin{lem}  \label{lem-Drinfeld-smooth-over-U-I}  (\cite{coho-filt} Lemma 3.3.4)
Let $\mc F$ be a constructible sheaf over $X^I$, equiped with an action of the partial Frobenius morphisms. Then there exists an open dense subscheme $U$ of $X$ such that $\mc F$ is smooth over $U^I$.
\cqfd
\end{lem}

\quad

\sssec{}   \label{subsection-when-smooth-sheaf-is-constant-abstract}
We begin by a too simple case to illustrate the case \ref{subsection-F-A-mod-tf-smooth-then-constant}. 
Let $\mc G$ be a constructible $E$-sheaf over $(\eta)^I$, equiped with an action of the partial Frobenius morphisms. %i.e. isomorphisms $F_{\{i\}}: \on{Frob}_{\{i\}}^* \mc G \isom \mc G$ commuting to each other and whose composition is the total Frobenius isomorphism $\on{Frob}^* \mc G \isom \mc G$ over $(\eta)^I$.

Firstly, by Lemma \ref{lem-Drinfeld-smooth-over-U-I}, $\mc G$ is smooth over $(\eta)^I$. In particular, $\mc G$ is smooth over 
$(\ov{\eta})^I$.

Secondly, as in \ref{subsection-def-FWeil-action}, $\restr{ \mc G }{ \ov{\eta_I} }$ is equiped with a continuous action of $\on{FWeil}(\eta_I, \ov{\eta_I})$. 
By Lemma \ref{lem-Drinfeld-E-vs}, this action factors through $\on{Weil}(\eta, \ov{\eta})^I$. 
%i.e. the action of $\Ker \Psi$ is trivial. %Note that 
%\begin{equation}
%\Ker \Psi = \on{Weil}((\ov{\eta})^I, \ov{\eta_I}).
%\end{equation}
We deduce that the action of $\on{Weil}((\ov{\eta})^I, \ov{\eta_I})$ on $\restr{ \mc G }{ \ov{\eta_I} }$ is trivial. 

%The category of smooth constructible $E$-sheaves over $(\ov{\eta})^I$ is equivalent to the category of finite dimensional $E$-vector spaces with continuous $\pi_1((\ov{\eta})^I, \ov{\eta_I})$-actions (cf. \cite{sga5} VI 1). 
By the equivalence in \ref{subsection-smooth-sheaf-pi-1-rep},
we deduce that $\restr{\mc G}{  (\ov{\eta})^I  }$ is a constant sheaf over $(\ov{\eta})^I$.

\sssec{}  \label{subsection-F-A-mod-tf-smooth-then-constant}
Let $A$ be a finitely generated $E$-algebra.
Let $\mc G$ be an ind-constructible $E$-sheaf over $(\eta)^I$ equiped with an action of the partial Frobenius morphisms and an action of $A$, such that 
\begin{itemize}
\item these two actions commute with each other

\item $\restr{ \mc G }{ \ov{\eta_I} }$ is an $A$-module of finite type

%\item $\mc G$ is an inductive limit of constructible $E$-sheaves

%the action of $\pi_1^{\on{geom}}(\eta_I, \ov{\eta_I})$ on $\restr{ \mc F }{ \ov{\eta_I} }$ is continuous, 
\end{itemize}

As in \ref{subsection-def-FWeil-action}, $\restr{ \mc G }{ \ov{\eta_I} }$ is equiped with a continuous $A$-linear action of $\on{FWeil}(\eta_I, \ov{\eta_I})$. 
%whose restriction to $\pi_1^{\on{geom}}(\eta_I, \ov{\eta_I})$ is , that is to say, $\restr{ \mc F }{ \ov{\eta_I} }$ is a union of finite dimensional $E$-vector subspaces which are stable under $\pi_1^{\on{geom}}(\eta_I, \ov{\eta_I})$ and on which the action of $\pi_1^{\on{geom}}(\eta_I, \ov{\eta_I})$ are continuous.
By Lemma \ref{lem-Drinfeld-A-mod}, this action factors through $\on{Weil}(\eta, \ov{\eta})^I$. Thus the action of $\on{Weil}((\ov{\eta})^I, \ov{\eta_I})$ on $\restr{ \mc G }{ \ov{\eta_I} }$ is trivial. 

However, in general $\mc G$ may not be ind-smooth over $(\ov{\eta})^I$. For example, 
let $I = \{1, 2\}$ and $\mc E_{n}$ be the extension by zero of the constant sheaf $E$ over $\Frob_{\{1\}}^{n}(\Delta)$, where $\Delta$ is the image of the diagonal morphism $X \hookrightarrow X^2$. Then $\mc G = \oplus_{n \in \Z} \mc E_{n}$ satisfies the above condition, but is not ind-smooth.

\sssec{}
In our situation \ref{subsection-action-Weil-I-on-coho}-\ref{subsection-constant-over-ov-eta-I} below, we will apply \ref{subsection-F-A-mod-tf-smooth-then-constant} to the cohomology sheaves of stacks of shtukas, and prove that for other reasons the cohomology sheaves are ind-smooth over $(\ov{\eta})^I$. Then they are constant sheaves over $(\ov{\eta})^I$.

%In the following, we will apply this situation to an $E$-sheaf over $(\eta)^I$ which is an inductive limit of subsheaves $\mc F$ as above.

\subsection{Action of $\on{Weil}(\eta, \ov{\eta})^I$ on cohomology sheaves of stacks of shtukas}   \label{subsection-action-Weil-I-on-coho}  

\sssec{}
Let $W \in \on{Rep}_{E}(\wh G^I)$, where $\on{Rep}_{E}(\wh G^I)$ denotes the category of finite dimensional $E$-linear representations of $\wh G^I$.

As in \cite{vincent} Définition 4.7 and \cite{coho-cusp} Section 2.5, 
let $\Cht_{G, N, I, W} / \Xi$ be the stack of $G$-shtukas, $\mc F_{G, N, I, W}^{\Xi}$ be the canonical perverse sheaf over $\Cht_{G, N, I, W} / \Xi$
with $E$-coefficients and $\mf p_G: \Cht_{G, N, I, W} / \Xi \rightarrow (X \sm N)^I$ be the morphism of paws. We denote by $\wh \Lambda_{G^{\mr{ad}}}^+$ the set of dominant coweights of $G^{\mr{ad}}$.
For any $\mu \in \wh \Lambda_{G^{\mr{ad}}}^+$ and any $j \in \Z$, we define the sheaf of degree $j$ cohomology with compact support for the Harder-Narasimhan truncation indexed by $\mu$
$$\mc H_{G, N, I, W}^{j, \, \leq \mu}:= R^j (\mf p_G)_! \restr{ \mc F_{G, N, I, W}^{\Xi} }{ \Cht_{G, N, I, W}^{\leq \mu} / \Xi} $$
%Since $\Cht_{G, N, I, W}^{\leq \mu} / \Xi$ is a Deligne-Mumford stack of finite type, $\mc H_{G, N, I, W}^{j, \, \leq \mu}$ 
It is a constructible $E$-sheaf over $(X \sm N)^I$.
We define the degree $j$ cohomology sheaf
$$\mc H_{G, N, I, W}^j:= \varinjlim _{\mu \in \wh \Lambda_{G^{\mr{ad}}}^+} \mc H_{G, N, I, W}^{j, \, \leq \mu}$$
in the abstract category of inductive limits of constructible $E$-sheaves over $(X \sm N)^I$. %That is to say, $\mc H_{G, N, I, W}^j$ is the sheafification of the presheaf associated to $\varinjlim _{\mu \in \wh \Lambda_{G^{\mr{ad}}}^+} \mc H_{G, N, I, W}^{j, \, \leq \mu}$. For any geometric point $\ov x$ of $X^I$, we have $$\restr{ \mc H_{G, N, I, W}^j  }{ \ov x } =\varinjlim _{\mu \in \wh \Lambda_{G^{\mr{ad}}}^+} ( \restr{ \mc H_{G, N, I, W}^{j, \, \leq \mu} }{\ov x} ).$$ Thus the definition in this paper is coherent with \cite{vincent}.

\sssec{}
In the following, to shorten the notation, we omit the indices $G$ and $N$, i.e. we write $\mc H_{I, W}^{j, \, \leq \mu}$ for $\mc H_{G, N, I, W}^{j, \, \leq \mu}$ and $\mc H_{I, W}^{j}$ for $\mc H_{G, N, I, W}^{j}$.

\sssec{}   \label{subsection-def-partial-Frob-on-H-I-W}
As in \cite{vincent} Section 4.3, there exists $\kappa \in \wh \Lambda_{G^{\mr{ad}}}^+$ big enough (depending on $W$) such that for any $i \in I$ and any $\mu$, we have morphisms of constructible $E$-sheaves over $(X \sm N)^I$:
$$F_{\{i\}}: \on{Frob}_{\{i\}}^* \mc H_{I, W}^{j, \, \leq \mu} \rightarrow \mc H_{I, W}^{j, \, \leq \mu + \kappa} $$
where $\on{Frob}_{\{i\}}$ is defined in \ref{subsection-def-Frob-J}.
The composition for all $i \in I$ is the total Frobenius morphism followed by an increase of the Harder-Narasimhan truncation: $$\on{Frob}^* \mc H_{I, W}^{j, \, \leq \mu} \rightarrow \mc H_{I, W}^{j, \, \leq \mu + (\sharp I) \kappa}$$

Taking inductive limit, we have canonical morphisms
$$F_{\{i\}}: \on{Frob}_{\{i\}}^* \mc H_{I, W}^{j}  \isom \mc H_{I, W}^{j} $$
whose composition is the total Frobenius morphism
$$\on{Frob}^* \mc H_{I, W}^{j} \isom \mc H_{I, W}^{j}.$$
Hence $\mc H_{I, W}^j$ is equiped with an action of the partial Frobenius morphisms.

As a consequence, as in \ref{subsection-def-FWeil-action}, $\restr{\mc H_{I, W}^j}{\ov{\eta_I}}$ is equiped with a continuous action of $\on{FWeil}(\eta_I, \ov{\eta_I})$.

\begin{prop}   \label{prop-FWeil-factors-though-Weil-I}
The action of $\on{FWeil}(\eta_I, \ov{\eta_I})$ on $\restr{\mc H_{I, W}^j}{\ov{\eta_I}}$ factors through $\on{Weil}(\eta, \ov{\eta})^I$.
\end{prop}

\begin{rem}
Proposition \ref{prop-FWeil-factors-though-Weil-I} is already proved in \cite{coho-filt} by using the main result in $loc.cit.$ that $\restr{\mc H_{I, W}^j}{\ov{\eta_I}}$ is of finite type as a module over a Hecke algebra. 

In the following, we give a new proof of Proposition \ref{prop-FWeil-factors-though-Weil-I}, which follows the arguments of \cite{vincent} Proposition 8.27. The proof uses a weaker result (Lemma \ref{lem-H-is-union-of-sub-modules} below) than \cite{coho-filt}. The advantage of this proof is that it is easily generalised 
%in Section \ref{section-cohomology-constant-I-1-I-2} and 
to not necessarily split groups in Section \ref{section-non-split}.
\end{rem}

To prove Proposition \ref{prop-FWeil-factors-though-Weil-I}, we begin by some preparations.

\sssec{}
For any family $(v_i)_{i \in I}$ of closed points of $X$, we denote by $\times_{i \in I} v_i$ their product over $\on{Spec} \Fq $. This is a finite union of closed points of $X^I$. %By induction on the cardinality of $I$ we prove that any dense open subscheme of $(X \sm N)^I$ contains such a product $\times_{i \in I} v_i$.

\begin{lem}  (\cite{vincent} Proposition 7.1) \label{lem-ES-relation}
Let $W = \boxtimes_{i \in I} W_i$ with $W_i \in \on{Rep}_E(\wh G)$. Let $(v_i)_{i \in I}$ be a family of closed points of $X \sm N$.
Then there exists $\kappa$, such that for any $\mu$ and any $i \in I$, we have
$$\sum_{\alpha =0}^{\on{dim}W_i} (-1)^{\alpha} S_{\wedge^{\on{dim}W_i - \alpha}W_i, v_i}(F_{\{i\}}^{\on{deg}(v_i)})^{\alpha} =0 \; \text{ in } \; \on{Hom}(\restr{\mc H_{I, W}^{j, \, \leq \mu} }{ \times_{i \in I} v_i}, \restr{\mc H_{I, W}^{j, \, \leq \mu + \kappa} }{ \times_{i \in I} v_i})$$
where $S_{\wedge^{\on{dim}W_i - \alpha}W_i, v_i}: \mc H_{I, W}^{j, \, \leq \mu} \rightarrow \mc H_{I, W}^{j, \, \leq \mu+ \kappa}$ is defined in \cite{vincent} Section 6, and we restrict it to $\restr{\mc H_{I, W}^{j, \, \leq \mu} }{ \times_{i \in I} v_i}$. 
\cqfd
\end{lem}

\sssec{}   \label{subsection-def-Hecke-local}
For any place $v$ of $X \sm N$, we denote by $\mc O_v$ the complete local ring at $v$ and $F_v$ its field of fractions.
Let $\ms H_{G, v}:=C_c(G(\mc O_v) \backslash G(F_v) / G(\mc O_v), E)$ be the Hecke algebra of $G$ at the place $v$.

\sssec{}
(\cite{vincent} Section 4.4)
Let $$T(h_{\wedge^{\on{dim}W_i - \alpha}W_i, v_i}) : \restr{\mc H_{I, W}^{j, \, \leq \mu} }{ (X \sm (N \cup v_i))^I} \rightarrow \restr{ \mc H_{I, W}^{j, \, \leq \mu+\kappa'}}{(X \sm (N \cup v_i))^I }$$ 
be the Hecke operator in $\ms H_{G, v_i}$ defined by Hecke correspondence.

\begin{lem} \label{lem-S-equal-T} (\cite{vincent} Proposition 6.2)
The operator $S_{\wedge^{\on{dim}W_i - \alpha}W_i, v_i}$, which is a morphism of sheaves over $(X \sm N)^I$, extends the action of the Hecke operator $T(h_{\wedge^{\on{dim}W_i - \alpha}W_i, v_i}) \in \ms H_{G, v_i}$, which is a morphism of sheaves over $(X \sm (N \cup v_i))^I$.
%acts by an element of $\ms H_{G, v_i}$.

\cqfd
\end{lem}

The combination of Lemma \ref{lem-ES-relation} and Lemma \ref{lem-S-equal-T} is called the Eichler-Shimura relations.

\quad

We will use Lemma \ref{lem-ES-relation} and Lemma \ref{lem-S-equal-T} to prove the following Lemma \ref{lem-H-is-union-of-sub-modules}.

\begin{lem}   \label{lem-H-is-union-of-sub-modules}
$\restr{\mc H_{I, W}^j}{\ov{\eta_I}}$ is an increasing union of $E$-vector subspaces $\mf M$ which are stable by $\on{FWeil}(\eta_I, \ov{\eta_I})$, and for which there exists a family $(v_i)_{i \in I}$ of closed points in $X \sm N$ (depending on $\mf M$) such that $\mf M$ is stable under the action of $\otimes_{i \in I} \ms H_{G, v_i}$ and is of finite type as module over $\otimes_{i \in I} \ms H_{G, v_i}$.
\end{lem}
\dem
Since the category $\on{Rep}_E(\wh G^I)$ is semisimple, it is enough to prove the lemma for $W$ irreducible, which is of the form $W = \boxtimes_{i \in I} W_i$ with $W_i \in \on{Rep}_E(\wh G)$ (after increasing $E$).

For any $\mu \in \wh \Lambda_{G^{\mr{ad}}}^{+}$, we choose a dense open subscheme $\Omega$ of $(X \sm N)^I$ such that $\restr{\mc H_{I, W}^{j, \, \leq \mu}}{\Omega}$ is smooth. We choose a closed point $v$ of $\Omega$.
Let $v_i$ be the image of $v$ under $(X \sm N)^I \xrightarrow{\on{pr}_i} X \sm N$, where $\on{pr}_i$ is the projection to the $i$-th factor. Then $\times_{i \in I} v_i$ is a finite union of closed points containing $v$. 
Let $\mf{M}_{\mu}$ be the image of
\begin{equation}    \label{equation-def-mf-M-mu}
\sum_{(n_i)_{i \in I} \in \N^I} (\otimes_{i \in I} \ms H_{G, v_i} ) \cdot \restr{\big( \prod_{i \in I} F_{\{i\}}^{n_i} ( (\prod_{i \in I} \on{Frob}_{\{i\}}^{n_i})^* \mc H_{I, W}^{j, \, \leq \mu}   )  \big)}{\ov{\eta_I}}  
\end{equation}
in $\restr{\mc H_{I, W}^j}{\ov{\eta_I}}$.
We have $\on{Im}( \restr{\mc H_{I, W}^{j, \, \leq \mu}}{\ov{\eta_I}} \rightarrow \restr{\mc H_{I, W}^{j}}{\ov{\eta_I}} ) \subset \mf M_{\mu}$ and 
\begin{equation}  \label{equation-H-I-W-egale-lim-M-mu}
\restr{\mc H_{I, W}^j}{\ov{\eta_I}} = \bigcup_{\mu} \mf M_{\mu}.
\end{equation}
By definition, $\mf{M}_{\mu}$ is stable under the action of the partial Frobenius morphisms and the action of $\on{Weil}(\eta_I, \ov{\eta_I})$, so is stable by $\on{FWeil}(\eta_I, \ov{\eta_I})$. We only need to prove that $\mf{M}_{\mu}$ is of finite type as $\otimes_{i \in I} \ms H_{G, v_i}$-module.

We fix a geometric point $\ov v$ over $v$ and a specialization map $\mf{sp}_v: \ov{\eta_I} \rightarrow \ov v$.
For any $n_i$, since $$F_{\{i\}}^{\on{deg}(v_i)n_i}: (\on{Frob}_{\{i\}}^{\on{deg}(v_i)n_i})^* \mc H_{I, W}^{j, \, \leq \mu} \rightarrow \mc H_{I, W}^j $$ is a morphism of sheaves, the specialization map $\mf{sp}_v$ induces a commutative diagram
\begin{equation}   \label{equation-Frob-v-to-Frob-eta-I-commute}
\xymatrix{
\restr{ (\on{Frob}_{\{i\}}^{\on{deg}(v_i)n_i})^* \mc H_{I, W}^{j, \, \leq \mu} }{\ov v}  \ar[r]^{\mf{sp}_v^*}_{\simeq}  \ar[d]_{F_{\{i\}}^{\on{deg}(v_i)n_i}}
& \restr{ (\on{Frob}_{\{i\}}^{\on{deg}(v_i)n_i})^* \mc H_{I, W}^{j, \, \leq \mu} }{\ov{\eta_I}}  \ar[d]^{F_{\{i\}}^{\on{deg}(v_i)n_i}}  \\
\restr{  \mc H_{I, W}^j }{\ov v}  \ar[r]^{\mf{sp}_v^*}
& \restr{  \mc H_{I, W}^j }{\ov{\eta_I}}
}
\end{equation}
The upper line of (\ref{equation-Frob-v-to-Frob-eta-I-commute}) is an isomorphism because $(\on{Frob}_{\{i\}}^{\on{deg}(v_i)n_i})^* \mc H_{I, W}^{j, \, \leq \mu}$ is smooth over $(\on{Frob}_{\{i\}}^{\on{deg}(v_i)n_i})^{-1} \Omega$. Note that $$\on{Frob}_{\{i\}}^{\on{deg}(v_i)n_i} (v) = v \in \Omega$$ thus $v \in (\on{Frob}_{\{i\}}^{\on{deg}(v_i)n_i})^{-1} \Omega$.

%Thus $$\restr{ (\on{Frob}_{\{i\}}^{\on{deg}(v_i)n_i})^* \mc H_{I, W}^{j, \, \leq \mu} }{\ov v} = \restr{  \mc H_{I, W}^{j, \, \leq \mu} }{\on{Frob}_{\{i\}}^{\on{deg}(v_i)n_i}(\ov v)} = \restr{  \mc H_{I, W}^{j, \, \leq \mu} }{\ov v}.$$

By Lemma \ref{lem-ES-relation}, for each $i \in I$, we have 
$$\sum_{\alpha =0}^{\on{dim}W_i} (-1)^{\alpha} S_{\wedge^{\on{dim}W_i - \alpha}W_i, v_i}(F_{\{i\}}^{\on{deg}(v_i)})^{\alpha} =0 \quad \text{ in }\on{Hom}(\restr{\mc H_{I, W}^{j, \, \leq \mu} }{ \ov v }, \restr{\mc H_{I, W}^{j} }{ \ov v })$$
We deduce that
$$
\begin{aligned}
& F_{\{i\}}^{\on{deg}(v_i)\on{dim}W_i} \big( \restr{(\on{Frob}_{\{i\}}^{\on{deg}(v_i)\on{dim}W_i})^* \mc H_{I, W}^{j, \, \leq \mu} }{ \ov v} \big) \\
\subset & \sum_{\alpha =0}^{\on{dim}W_i-1}  S_{\wedge^{\on{dim}W_i - \alpha}W_i, v_i} F_{\{i\}}^{\on{deg}(v_i)\alpha} \big( \restr{(\on{Frob}_{\{i\}}^{\on{deg}(v_i)\alpha})^* \mc H_{I, W}^{j, \, \leq \mu} }{\ov v} \big)
\end{aligned}
$$
Since $S_{\wedge^{\on{dim}W_i - \alpha}W_i, v_i}$ and $F_{\{i\}}$ are morphisms of sheaves, they commute with $\mf{sp}_v^*$. We have
$$
\begin{aligned}
& F_{\{i\}}^{\on{deg}(v_i)\on{dim}W_i} \big( \mf{sp}_v^* \restr{(\on{Frob}_{\{i\}}^{\on{deg}(v_i)\on{dim}W_i})^* \mc H_{I, W}^{j, \, \leq \mu} }{\ov v} \big)  \\
\subset & \sum_{\alpha =0}^{\on{dim}W_i-1}  S_{\wedge^{\on{dim}W_i - \alpha}W_i, v_i} F_{\{i\}}^{\on{deg}(v_i)\alpha} \big( \mf{sp}_v^* \restr{(\on{Frob}_{\{i\}}^{\on{deg}(v_i)\alpha})^* \mc H_{I, W}^{j, \, \leq \mu} }{\ov v} \big)
\end{aligned}
$$
Since the upper line of (\ref{equation-Frob-v-to-Frob-eta-I-commute}) is an isomorphism, we deduce that
\begin{equation}   \label{equation-Frob-dim-W-includ-sum-Frob-dim-small-eta-I}
\begin{aligned}
& F_{\{i\}}^{\on{deg}(v_i)\on{dim}W_i} \big(  \restr{(\on{Frob}_{\{i\}}^{\on{deg}(v_i)\on{dim}W_i})^* \mc H_{I, W}^{j, \, \leq \mu} }{\ov{\eta_I}} \big)  \\
\subset & \sum_{\alpha =0}^{\on{dim}W_i-1}  S_{\wedge^{\on{dim}W_i - \alpha}W_i, v_i} F_{\{i\}}^{\on{deg}(v_i)\alpha} \big( \restr{(\on{Frob}_{\{i\}}^{\on{deg}(v_i)\alpha})^* \mc H_{I, W}^{j, \, \leq \mu} }{\ov{\eta_I}} \big)  
\end{aligned} 
\end{equation}
By \cite{vincent}, the action of the partial Frobenius morphisms commute with the action of Hecke algebras. Moreover, by Lemma \ref{lem-S-equal-T}, $S_{\wedge^{\on{dim}W_i - \alpha}W_i, v_i}$ acts over $\ov{\eta_I}$ by an element of $\ms H_{G, v_i}$. We deduce that
\begin{equation}   \label{equation-RHS-equal-Hecke-action}
\begin{aligned}
\text{RHS of } (\ref{equation-Frob-dim-W-includ-sum-Frob-dim-small-eta-I}) = & \sum_{\alpha =0}^{\on{dim}W_i-1}   F_{\{i\}}^{\on{deg}(v_i)\alpha} \big( \restr{(\on{Frob}_{\{i\}}^{\on{deg}(v_i)\alpha})^* (S_{\wedge^{\on{dim}W_i - \alpha}W_i, v_i} \mc H_{I, W}^{j, \, \leq \mu} ) }{\ov{\eta_I}} \big) \\
\subset & \sum_{\alpha =0}^{\on{dim}W_i-1}   F_{\{i\}}^{\on{deg}(v_i)\alpha} \big( \restr{(\on{Frob}_{\{i\}}^{\on{deg}(v_i)\alpha})^* (\ms{H}_{G, v_i} \cdot \mc H_{I, W}^{j, \, \leq \mu} ) }{\ov{\eta_I}} \big) 
\end{aligned}
\end{equation}

We deduce from (\ref{equation-Frob-dim-W-includ-sum-Frob-dim-small-eta-I}) and (\ref{equation-RHS-equal-Hecke-action}) that $\mf{M}_{\mu}$ is equal to the image of
\begin{equation}   \label{equation-mf-M-equal-Hecke-Frob-finite}
\sum_{(n_i)_{i \in I} \in \prod_{i \in I}\{0, 1, \cdots, \on{deg}(v_i)(\on{dim}W_i -1) \}} (\otimes_{i \in I} \ms H_{G, v_i} ) \cdot \restr{\big( \prod_{i \in I} F_{\{i\}}^{n_i} ( (\prod_{i \in I} \on{Frob}_{\{i\}}^{n_i})^* \mc H_{I, W}^{j, \, \leq \mu}   )  \big)}{\ov{\eta_I}}  
\end{equation}
in $\restr{\mc H_{I, W}^j}{\ov{\eta_I}}$.
Thus $\mf{M}_{\mu}$ is of finite type as $\otimes_{i \in I} \ms H_{G, v_i}$-module.
\cqfd

\quad

\noindent {\bf Proof of Proposition \ref{prop-FWeil-factors-though-Weil-I}:}
For every $\mu$, applying Lemma \ref{lem-Drinfeld-A-mod} to $A = \otimes_{i \in I} \ms H_{G, v_i}$ and $M = \mf M_{\mu}$ (which is possible because of Lemma \ref{lem-H-is-union-of-sub-modules}), we deduce that
the action of $\on{FWeil}(\eta_I, \ov{\eta_I})$ on $\mf{M}_{\mu}$ factors through $\on{Weil}(\eta, \ov{\eta})^I$.

Since $\restr{\mc H_{I, W}^j}{\ov{\eta_I}}$ is an inductive limit of $\mf{M}_{\mu}$, we deduce Proposition \ref{prop-FWeil-factors-though-Weil-I}.
\cqfd

\subsection{Smoothness of cohomology sheaves over $(\ov{\eta})^I$}  

\sssec{}
In \cite{vincent} Proposition 8.32, V. Lafforgue proved that for a specialization map $\ov{\eta_I} \rightarrow \Delta(\ov{\eta})$, the induced morphism $$\restr{  \mc H_{I, W}^j   }{ \Delta(\ov{\eta}) }  \rightarrow \restr{  \mc H_{I, W}^j   }{ \ov{\eta_I}  } $$
is injective.
In fact, the same argument gives a more general result: the injectivity of morphism (\ref{equation-sp-x-sheaf}) below.

Moreover, a similar argument as \cite{vincent} Proposition 8.31 gives a more general result: the surjectivity of morphism (\ref{equation-sp-x-sheaf}) below.

\sssec{}
%Let $x$ be a point of $(\ov{\eta})^I$. We fix a geometric point $\ov x$ over $x$.
Let $\ov x$ be a geometric point of $(\ov{\eta})^I$.
We fix a specialization map
\begin{equation*}
\mf{sp}_{\ov x}: \ov{\eta_I} \rightarrow \ov x
\end{equation*}
It induces a morphism
\begin{equation}    \label{equation-sp-x-sheaf}
\mf{sp}_{\ov x}^*: \restr{  \mc H_{I, W}^j   }{ \ov x }  \rightarrow \restr{  \mc H_{I, W}^j   }{ \ov{\eta_I}  } 
\end{equation}

\begin{prop}   \label{prop-H-x-H-eta-I-isom}
The morphism (\ref{equation-sp-x-sheaf}) is an isomorphism.
\end{prop}

To prove Proposition \ref{prop-H-x-H-eta-I-isom}, we need the following lemma.

\begin{lem}   \label{lem-Eike-Lau}  (\cite{eike-lau} Lemma 9.2.1)
Let $x$ be a point of $(\eta)^I$. The set $\{ (\prod_{i \in I} \Frob_{ \{i \} }^{m_i}  )(x) , (m_i)_{i \in I} \in \N^I \}$ is Zariski dense in $X^I$.
\end{lem}
\dem
In fact, the Zariski closure of this set is a closed subscheme $Z$ of $X^I$, invariant by the partial Frobenius morphisms. If $Z$ is not equal to $X^I$, by Lemma 9.2.1 of \cite{eike-lau} (recalled in the proof of Lemme 8.12 of \cite{vincent}), $Z$ is included in a finite union of vertical divisors (i.e. the inverse image of a closed point by one of the projections $X^I \rightarrow X$). However, the image of $x$ in $X^I$ is not included in any vertical divisor. This is a contradiction. We deduce that $Z = X^I$.
\cqfd

\quad

\noindent {\bf Proof of Proposition \ref{prop-H-x-H-eta-I-isom}:} 
%It is enough to prove for $W = \boxtimes_{j \in I} W_j$.
Since the category $\on{Rep}_E(\wh G^I)$ is semisimple, it is enough to prove the proposition for $W$ irreducible, which is of the form $W = \boxtimes_{i \in I} W_i$ with $W_i \in \on{Rep}_E(\wh G)$ (after increasing $E$).

\quad

{\bf Injectivity}: the proof is the same as Proposition 8.32 of \cite{vincent}, except that we replace everywhere $\Delta(\ov{\eta})$ by $\ov x$ and replace everywhere $\Delta(\ov v)$ by $\ov y$ (defined below). For the convenience of the reader, we briefly recall the proof. Let $a \in \on{Ker} (\mf{sp}_{\ov x}^*)$. We want to prove that $a=0$. 
%The proof consists of three steps: 

%(1) Use the Eichler-Shimura relations to get (\ref{equation-sum-S-V-a-n-i-zero}):

There exists $\mu_0$ large enough and $\wt a \in \restr{  \mc H_{I, W}^{j, \, \leq \mu_0}   }{ \ov x }$, such that $a$ is the image of $\wt a$ in $\restr{  \mc H_{I, W}^j   }{ \ov x }$. We denote by $x$ the image of $\ov x$ in $(X \sm N)^I$ and $\ov{\{x\}}$ the Zariski closure of $x$. Let $\Omega_0$ be a dense open subscheme of $\ov{\{x\}}$ such that $\restr{  \mc H_{I, W}^{j, \, \leq \mu_0}   }{ \Omega_0 }$ is smooth. 
Let $y$ be a closed point in $\Omega_0$. Let $\ov y$ be a geometric point over $y$ and $\mf{sp}_y: \ov{x} \rightarrow \ov{y}$ a specialization map over $\Omega_0$. 
We have a commutative diagram
\begin{equation}
\xymatrixrowsep{2pc}
\xymatrixcolsep{5pc}
\xymatrix{
\restr{  \mc H_{I, W}^{j, \, \leq \mu_0}   }{ \ov y }   \ar[r]^{ \mf{sp}_y^*  }_{\simeq}   \ar[d]
& \restr{  \mc H_{I, W}^{j, \, \leq \mu_0}   }{ \ov{x}  }  \ar[d] \\
\restr{  \mc H_{I, W}^{j}   }{ \ov y }    \ar[r]^{ \mf{sp}_y^*  }
& \restr{  \mc H_{I, W}^{j}   }{ \ov{x}  } 
}
\end{equation}
The upper horizontal morphism is an isomorphism because $\restr{  \mc H_{I, W}^{j, \, \leq \mu_0}   }{ \Omega_0 }$ is smooth.
Thus there exists $\wt b \in \restr{  \mc H_{I, W}^{j, \, \leq \mu_0}   }{ \ov y } $ such that $\wt a = \mf{sp}_y^* (\wt b)$. 
Let $b$ be the image of $\wt b$ in $\restr{  \mc H_{I, W}^{j}   }{ \ov y }$. We have $a = \mf{sp}_y^* (b)$. 

Let $y_i$ be the image of $y$ by $(X \sm N)^I \xrightarrow{\on{pr}_i} X \sm N$. Then $\times_{i \in I} y_i$ is a finite union of closed points containing $y$. 
Let $d_i = \on{deg}(y_i)$.
For any $(n_i)_{i \in I} \in \N^I$, we have $\prod_{i \in I} \on{Frob}_{\{i\}}^{d_i n_i}(\ov y) = \ov y$. (Note that in general $\prod_{i \in I} \on{Frob}_{\{i\}}^{d_i n_i}(\ov x) \neq \ov x$.)
We have the partial Frobenius morphism 
$$\prod_{i \in I} F_{\{i\}}^{d_i n_i}: \restr{  \mc H_{I, W}^{j}   }{ \ov y } =
\restr{ (\prod_{i \in I} \on{Frob}_{\{i\}}^{d_i n_i})^* \mc H_{I, W}^{j}   }{ \ov y } \rightarrow \restr{  \mc H_{I, W}^{j}   }{ \ov y } $$ 
%We have morphisms
%\begin{equation}
%\xymatrix{
%\restr{  \mc H_{I, W}^{j}   }{ \ov y }   \ar[r]\ar[r]^{ \mf{sp}_y^*  }   \ar[d]_{ \prod_{i \in I} F_{\{i\}}^{d_i n_i}  }
%& \restr{  \mc H_{I, W}^{j}   }{ \ov{x}  }  \\
%\restr{  \mc H_{I, W}^{j}   }{ \ov y }   \ar[r]\ar[r]^{ \mf{sp}_y^*  }
%& \restr{  \mc H_{I, W}^{j}   }{ \ov x }  
%}
%\end{equation}
Let 
$$b_{  (n_i)_{i \in I}   } = \prod_{i \in I} F_{\{i\}}^{d_i n_i} (b) \in \restr{  \mc H_{I, W}^{j}   }{ \ov y } \quad \text{ and } \quad a_{  (n_i)_{i \in I}   } = \mf{sp}_y^* (  b_{  (n_i)_{i \in I}   }  ) \in \restr{  \mc H_{I, W}^{j}   }{ \ov x }  .$$
In particular, $b_{(0)_{i \in I}} = b$ and $a_{(0)_{i \in I}} = a$.

Let $d = \on{deg} (y) = ppcm (\{d_i\}_{i \in I})$.
Note that $\prod_{i \in I} \on{Frob}_{\{i\}}$ is the total Frobenius morphism, thus the morphism
$$\prod_{i \in I} F_{\{i\}}^{dn}: \restr{  \mc H_{I, W}^{j}   }{ \ov x }  \rightarrow
\restr{  \mc H_{I, W}^{j}   }{ \ov x }  $$
is bijective. We have 
\begin{equation}   \label{equation-a-plus-n-equal-Frob-a}
a_{  (n_i + nd / d_i)_{i \in I}   } = \prod_{i \in I} F_{\{i\}}^{dn} ( a_{  (n_i)_{i \in I}   } ) .
\end{equation}

\quad

\cite{vincent} Lemme 8.33 is still true if we replace everywhere $\Delta(\ov{\eta})$ by $\ov x$ and replace everywhere $\Delta(\ov v)$ by $\ov y$. Thus we have:
%
%
%By Lemma \ref{lem-ES-relation}, for every $k \in I$ and every $(n_i)_{i \in I} \in \N^I$, we have
%\begin{equation}   
%\sum_{\alpha =0}^{\on{dim}W_k} (-1)^{\alpha} S_{\bigwedge^{\on{dim} W_k - \alpha}W_k, y_k} (b_{(n_i + \alpha \delta_{i, k})_{i \in I}}) =0 \; \text{ in } \restr{ \mc H_{I, W}^j }{ \ov y } .
%\end{equation}
%We deduce that

(1) for all $k \in I$ and for all $(n_i)_{i \in I} \in \N^I$,
\begin{equation}   \label{equation-sum-S-V-a-n-i-zero}
\sum_{\alpha =0}^{\on{dim}W_k} (-1)^{\alpha} S_{\bigwedge^{\on{dim} W_k - \alpha}W_k, y_k} (a_{(n_i + \alpha \delta_{i, k})_{i \in I}}) =0 \; \text{ in } \restr{ \mc H_{I, W}^j }{ \ov x } .
\end{equation}

(2) Let $\mu_1 \geq \mu_0$ such that $\mf{sp}_{\ov x}^* (\wt a) \in \restr{  \mc H_{I, W}^{j, \, \leq \mu_0}   }{ \ov{\eta_I}  } $ has zero image in $\restr{  \mc H_{I, W}^{j, \, \leq \mu_1}   }{ \ov{\eta_I}  } $.
Let $\Omega_1$ be a dense open subscheme of $(X \sm N)^I$ such that $\restr{  \mc H_{I, W}^{j, \, \leq \mu_1}   }{ \Omega_1 }$ is smooth. 
%As in \cite{vincent} Lemme 8.33 (b) (replace everywhere $\Delta(\ov{\eta})$ by $\ov x$ and $\Delta(\ov v)$ by $\ov y$), we prove that 
Then for every $(m_i)_{i \in I} \in \N^I$ such that $\prod_{i \in I} \on{Frob}_{\{i\}}^{d_i m_i} (x) \in \Omega_1$, we have $a_{(m_i)_{i \in I}} = 0$ in $\restr{ \mc H_{I, W}^j }{ \ov x }$. 

\quad

Note that the open subscheme
\begin{equation}    \label{equation-intersect-Frob-Omega-1}
\bigcap_{ (\alpha_i)_{i \in I} \in \prod_{i \in I} \{0, \cdots, \on{dim}W_i - 1 \}   } (\prod_{i \in I} \on{Frob}_{\{i\}}^{d_i \alpha_i})^{-1}(\Omega_1)
\end{equation}
is also dense in $X^I$.
By Lemma \ref{lem-Eike-Lau}, there exists $(N_i)_{i \in I} \in \N^I$ such that $\prod_{i \in I} \on{Frob}_{\{i\}}^{d_iN_i} (x)$ is in (\ref{equation-intersect-Frob-Omega-1}). We deduce
$$\prod_{i \in I} \on{Frob}_{\{i\}}^{d_i(N_i + \alpha_i)} (x) \in \Omega_1 \text{ for all } (\alpha_i)_{i \in I} \in \prod_{i \in I} \{0, \cdots, \on{dim}W_i - 1 \} .$$
By (2), we deduce that
\begin{equation}   \label{equation-a-alpha-equal-zero-finite-alpha}
a_{(N_i+\alpha_i)_{i \in I}} =0 \text{ for all } (\alpha_i)_{i \in I} \in \prod_{i \in I} \{0, \cdots, \on{dim}W_i - 1 \} .
\end{equation}
By (1), for every $k \in I$ and $(n_i)_{i \in I} \in \N^I$,
\begin{equation}   \label{equation-a-dim-W-equal-sum-leq-dim-W}
a_{(n_i + \on{dim}W_k \delta_{i, k})_{i \in I}} =
\sum_{\alpha =0}^{\on{dim}W_k -1} (-1)^{\alpha+\dim W_k} S_{\bigwedge^{\on{dim} W_k - \alpha}W_k, y_k} (a_{(n_i + \alpha \delta_{i, k})_{i \in I}})
\end{equation}
%\begin{equation}  \label{equation-a-N-i-alpha-i-dim-W-i-zero}
%a_{(N_i+\alpha_i)_{i \in I}} =0 \text{ for all } (\alpha_i)_{i \in I} \in \prod_{i \in I} \{0, \cdots, \on{dim}W_i \} .
%\end{equation}
%Let $N'_i = N_i +1$, (\ref{equation-a-N-i-alpha-i-dim-W-i-zero}) implies
%\begin{equation*}
%a_{(N_i'+\alpha_i)_{i \in I}} =0 \text{ for all } (\alpha_i)_{i \in I} \in \prod_{i \in I} \{0, \cdots, \on{dim}W_i - 1 \} .
%\end{equation*}
%We deduce from Lemma \ref{lem-8.33-vincent} (a) that
%\begin{equation*}  
%a_{(N_i'+\alpha_i)_{i \in I}} =0 \text{ for all } (\alpha_i)_{i \in I} \in \prod_{i \in I} \{0, \cdots, \on{dim}W_i \} .
%\end{equation*}
%Let $N''_i = N'_i +1$, continue the process... Finally 
Using (\ref{equation-a-alpha-equal-zero-finite-alpha}) and (\ref{equation-a-dim-W-equal-sum-leq-dim-W}), by induction we deduce that
\begin{equation*}  
a_{(n_i)_{i \in I}} =0 \text{ for all } (n_i)_{i \in I} \in \N^I \text{ such that } n_i \geq N_i , \, \forall i \in I
\end{equation*}
Thus for $n \geq N_i$ for all $i \in I$, we have $a_{(nd/ d_i)_{i \in I}} =0$.
Then (\ref{equation-a-plus-n-equal-Frob-a}) implies $a_{(0)_{i \in I}}=0.$ This proves the injectivity of $\mf{sp}_{\ov x}^*$.

\quad

{\bf Surjectivity}: the proof is similar as \cite{vincent} Proposition 8.31. 
By (\ref{equation-H-I-W-egale-lim-M-mu}), we have $\restr{\mc H_{I, W}^j}{\ov{\eta_I}} = \bigcup_{\mu} \mf M_{\mu}$, where $\mf M_{\mu}$ is defined in (\ref{equation-def-mf-M-mu}).
To prove that $\mf{sp}_{\ov x}^*$ is surjective, 
it is enough to prove %that for every sub-$E$-module $\mf M \subset \restr{\mc H_{I, W}^j}{\ov{\eta_I}}$ which is stable by $\on{FWeil}(\eta_I, \ov{\eta_I})$ and is of finite type as $\otimes_{i \in I} \ms H_{G, v_i}$-module for some family $(v_i)_{i \in I}$ of closed points in $X \sm N$, we have $\mf M \subset \on{Im}(\mf{sp}_x^*)$.
that for every $\mu$, we have $\mf M_{\mu} \subset \on{Im}(\mf{sp}_{\ov x}^*)$.

%By (\ref{equation-mf-M-equal-Hecke-Frob-finite}), there exists $\mu_2$ large enough such that 
%\begin{equation}   \label{equation-mf-M-equal-Hecke-Frob-finite-mu-2}
%\mf M_{\mu} \subset (\otimes_{i \in I} \ms H_{G, v_i}) \cdot \restr{\mc H_{I, W}^{j, \, \leq \mu_2}}{\ov{\eta_I}} .
%\end{equation}
%Let $\omega_1, \cdots, \omega_m$ be a family of generators of $\mf M_{\mu}$ as $\otimes_{i \in I} \ms H_{G, v_i}$-module. 
There exist $\mu_2$ large enough and $\omega_1, \cdots, \omega_m \in \restr{\mc H_{I, W}^{j, \, \leq \mu_2}}{\ov{\eta_I}} $
%$\wt{\omega_1}, \cdots, \wt{\omega_m} \in \restr{\mc H_{I, W}^{j, \, \leq \mu_2}}{\ov{\eta_I}} $, 
such that $\omega_1, \cdots, \omega_m$ is a family of generators of $\mf M_{\mu}$ as $\otimes_{i \in I} \ms H_{G, v_i}$-module. 
%for every $i \in \{1, \cdots, m\}$, $\omega_i$ is the image of $\wt{\omega_i}$ in $\restr{\mc H_{I, W}^{j}}{\ov{\eta_I}}$.
Let $\Omega_2$ be a dense open subscheme of $(X \sm N)^I$ such that $\restr{\mc H_{I, W}^{j, \, \leq \mu_2}}{\Omega_2}$ is smooth.  
%Let $\omega_1, \cdots, \omega_m \in \restr{\mc H_{I, W}^{j, \, \leq \mu_2}}{\ov{\eta_I}} $ be a family of generators of $\mf M_{\mu}$. Let $\mf G$ be the smooth subsheaf of $\restr{\mc H_{I, W}^{j, \, \leq \mu_2}}{\Omega_2}$ generated by $\omega_1, \cdots, \omega_m$.
By Lemma \ref{lem-Eike-Lau}, the set $\{ (\prod_{i \in I} \Frob_{ \{i \} }^{m_i}  )(x) , (m_i)_{i \in I} \in \N^I \}$ is Zariski dense in $X^I$. We deduce that there exists $(n_i)_{i \in I} \in \N^I$, such that
$(\prod_{i \in I} \Frob_{ \{i \} }^{n_i}  )(x) \in \Omega_2$.

Let the specialization map
$$(\prod_{i \in I} \Frob_{ \{i \} }^{n_i}  )\mf{sp}_{\ov x}: (\prod_{i \in I} \Frob_{ \{i \} }^{n_i}  )( \ov{\eta_I} )  \rightarrow  (\prod_{i \in I} \Frob_{ \{i \} }^{n_i}  )(\ov x) $$
be the image of $\mf{sp}_{\ov x}$ by $\prod_{i \in I} \Frob_{ \{i \} }^{n_i}$.
It induces a morphism
$$((\prod_{i \in I} \Frob_{ \{i \} }^{n_i}  )\mf{sp}_{\ov x})^*: \restr{\mc H_{I, W}^{j, \, \leq \mu_2} }{(\prod_{i \in I} \Frob_{ \{i \} }^{n_i}  )(\ov x)}   \rightarrow  \restr{\mc H_{I, W}^{j, \, \leq \mu_2} }{ (\prod_{i \in I} \Frob_{ \{i \} }^{n_i}  )( \ov{\eta_I} )  } $$
%We have a commutative diagram
%\begin{equation}
%\xymatrixrowsep{3pc}
%\xymatrixcolsep{5pc}
%\xymatrix{
%\restr{\mc H_{I, W}^{j, \, \leq \mu_2} }{(\prod_{i \in I} \Frob_{ \{i \} }^{n_i}  )(\ov x)}   \ar[r]^{ (\mf{sp}_{x, (n_i)_{i \in I}})^*  }_{\simeq}   \ar[d]_{\prod_{i \in I} F_{ \{i \} }^{n_i}}
%& \restr{\mc H_{I, W}^{j, \, \leq \mu_2} }{ (\prod_{i \in I} \Frob_{ \{i \} }^{n_i}  )( \ov{\eta_I} )  } \ar[d]^{\prod_{i \in I} F_{ \{i \} }^{n_i}}   \\
%\restr{\mc H_{I, W}^j }{\ov x}    \ar[r]^{ \mf{sp}_{x}^*  }
%& \restr{\mc H_{I, W}^j }{  \ov{\eta_I}   }  
%}
%\end{equation}
%The upper horizontal morphism 
This is an isomorphism because $\restr{  \mc H_{I, W}^{j, \, \leq \mu_2}   }{ \Omega_2 }$ is smooth. In particular, $$ \restr{\mc H_{I, W}^{j, \, \leq \mu_2} }{ (\prod_{i \in I} \Frob_{ \{i \} }^{n_i}  )( \ov{\eta_I} )  } \subset \on{Im} ((\prod_{i \in I} \Frob_{ \{i \} }^{n_i}  )\mf{sp}_{\ov x})^*.$$
%As a consequence, $\omega_1, \cdots, \omega_m \in \on{Im}(\mf{sp}_{x}^*)$. 
Since the action of the Hecke algebra is given by morphisms of sheaves, it commutes with $((\prod_{i \in I} \Frob_{ \{i \} }^{n_i}  )\mf{sp}_{\ov x})^*$. We deduce that 
\begin{equation}   \label{equation-Hecke-H-leq-mu-2-in-imega}
(\otimes_{i \in I} \ms H_{G, v_i}) \cdot \restr{\mc H_{I, W}^{j, \, \leq \mu_2} }{ (\prod_{i \in I} \Frob_{ \{i \} }^{n_i}  )( \ov{\eta_I} )  }   \subset \on{Im}  (  (\prod_{i \in I} \Frob_{ \{i \} }^{n_i}  )\mf{sp}_{\ov x})^*.
\end{equation}

%$$\otimes_{i \in I} \ms H_{G, v_i} \cdot \restr{\mc H_{I, W}^{j, \, \leq \mu_2} }{(\prod_{i \in I} \Frob_{ \{i \} }^{n_i}  )(\ov x)}   \xrightarrow{  (\prod_{i \in I} \Frob_{ \{i \} }^{n_i}  )(\mf{sp}_{x})^*  }  \otimes_{i \in I} \ms H_{G, v_i} \cdot \restr{\mc H_{I, W}^{j, \, \leq \mu_2} }{ (\prod_{i \in I} \Frob_{ \{i \} }^{n_i}  )( \ov{\eta_I} )  } $$ is an isomorphism.

%We have the partial Frobenius morphism $$\prod_{i \in I} F_{ \{i \} }^{- n_i}: (\prod_{i \in I} \Frob_{ \{i \} }^{ - n_i})^* \mc H_{I, W}^{j, \, \leq \mu_2} \rightarrow \mc H_{I, W}^{j, \, \leq \mu_2+\kappa}$$
%It induces
%$$\prod_{i \in I} F_{ \{i \} }^{ - n_i}: \restr{\mc H_{I, W}^{j, \, \leq \mu_2} }{ \ov{\eta_I} } \rightarrow \restr{ \mc H_{I, W}^{j, \, \leq \mu_2+\kappa} }{(\prod_{i \in I} \Frob_{ \{i \} }^{n_i}  )(\ov{\eta_I})}$$

%As a consequence, $\omega_1, \cdots, \omega_m \in \on{Im} (\prod_{i \in I} \Frob_{ \{i \} }^{n_i}  )(\mf{sp}_{x})^*$. 

\quad

As in the proof of \cite{vincent} Proposition 8.31. we have a commutative diagram
\begin{equation}  \label{equation-Frob-n-i-and-sp-x}
\xymatrixrowsep{2pc}
\xymatrixcolsep{8pc}
\xymatrix{
\restr{\mc H_{I, W}^{j} }{(\prod_{i \in I} \Frob_{ \{i \} }^{n_i}  )(\ov x)}  \ar[r]^{ ((\prod_{i \in I} \Frob_{ \{i \} }^{n_i}  )\mf{sp}_{\ov x})^* }  \ar[d]^{\prod_{i \in I} F_{ \{i \} }^{n_i} }_{\simeq}
& \restr{\mc H_{I, W}^{j} }{(\prod_{i \in I} \Frob_{ \{i \} }^{n_i}  )(\ov{\eta_I})} \ar[d]^{\prod_{i \in I} F_{ \{i \} }^{n_i} }_{\simeq} 
& \mf M_{\mu}'  \ar@{^{(}->}[l] \ar[d]^{\simeq}   \\
\restr{\mc H_{I, W}^{j} }{ \ov x }  \ar[r]^{ \mf{sp}_{\ov x}^* }
& \restr{\mc H_{I, W}^{j} }{ \ov{\eta_I} } 
& \mf M_{\mu}  \ar@{^{(}->}[l] 
}
\end{equation}
%Let $\mf M' \subset \restr{\mc H_{I, W}^{j} }{(\prod_{i \in I} \Frob_{ \{i \} }^{n_i}  )(\ov{\eta_I})}$ such that $(\prod_{i \in I} F_{ \{i \} }^{n_i} ) \mf M' = \mf M$. 
where $\mf M_{\mu}' := (\prod_{i \in I} F_{ \{i \} }^{n_i} )^{-1} \mf M_{\mu}$ is the inverse image of $\mf M_{\mu}$. %by the isomorphism on the right vertical arrow of (\ref{equation-Frob-n-i-and-sp-x}).
%Since the action of the partial Frobenius morphisms is a morphism of sheaves $\prod_{i \in I} F_{ \{i \} }^{n_i}: (\prod_{i \in I} \Frob_{ \{i \} }^{n_i})^* \mc H_{I, W}^{j} \isom \mc H_{I, W}^{j}$, and it commutes with the action of the Hecke algebra, we deduce from (\ref{equation-mf-M-equal-Hecke-Frob-finite-mu-2}) that
Note that $\mf M_{\mu}$ is stable under 
%the action of the partial Frobenius morphism and $\on{Weil}(\eta_I, \ov{\eta_I})$, thus it is stable under 
the action of $\on{FWeil}(\eta_I, \ov{\eta_I})$ (recall that the action is given in \ref{subsection-def-FWeil-action}). Thus any choice of isomorphism $\beta$ between $(\prod_{i \in I} \Frob_{ \{i \} }^{n_i}  )(\ov{\eta_I})$ and $\ov{\eta_I}$ identifies $\mf M_{\mu}'$ and $\mf M_{\mu}$:
$$
\xymatrix{
\restr{ \mc H_{I, W}^{j}   }{\ov{\eta_I}}   \ar[d]^{\beta^*}
& \mf M_{\mu}  \ar@{^{(}->}[l]  \ar[d]^{\simeq} 
& \omega_1, \cdots, \omega_m \\
\restr{\mc H_{I, W}^{j} }{(\prod_{i \in I} \Frob_{ \{i \} }^{n_i}  )(\ov{\eta_I})} \ar[d]^{\prod_{i \in I} F_{ \{i \} }^{n_i} }_{\simeq} 
& \mf M_{\mu}'  \ar@{^{(}->}[l] \ar[d]^{\simeq}  
& \wt{\omega_1}, \cdots, \wt{\omega_m}  \\
\restr{\mc H_{I, W}^{j} }{ \ov{\eta_I} } 
& \mf M_{\mu}  \ar@{^{(}->}[l] 
}
$$
Let $\wt{\omega_1}, \cdots, \wt{\omega_m} \in \restr{\mc H_{I, W}^{j, \, \leq \mu_2}}{(\prod_{i \in I} \Frob_{ \{i \} }^{n_i}  )(\ov{\eta_I})} $ denote the images of $\omega_1, \cdots, \omega_m$ by $\beta^*$. Then
$\wt{\omega_1}, \cdots, \wt{\omega_m} $ is a family of generators of $\mf M_{\mu}'$ as $\otimes_{i \in I} \ms H_{G, v_i}$-module. 
We deduce that
\begin{equation}
\mf M_{\mu}' \subset (\otimes_{i \in I} \ms H_{G, v_i}) \cdot \restr{\mc H_{I, W}^{j, \, \leq \mu_2}}{(\prod_{i \in I} \Frob_{ \{i \} }^{n_i}  )(\ov{\eta_I})}.
\end{equation}
By (\ref{equation-Hecke-H-leq-mu-2-in-imega}), $\mf M_{\mu}' \subset \on{Im} ((\prod_{i \in I} \Frob_{ \{i \} }^{n_i}  )\mf{sp}_{\ov x})^*$. We deduce from (\ref{equation-Frob-n-i-and-sp-x}) that $\mf M_{\mu} \subset \on{Im}(\mf{sp}_{\ov x}^*)$.
% $\otimes_{i \in I} \ms H_{G, v_i} \cdot \langle \omega_1, \cdots, \omega_m \rangle \subset \on{Im}(\prod_{i \in I} \Frob_{ \{i \} }^{n_i}  )(\mf{sp}_{x})^*$, i.e. $\mf M \subset \on{Im} (\prod_{i \in I} \Frob_{ \{i \} }^{n_i}  )(\mf{sp}_{x})^*$. 

\cqfd

\subsection{The cohomology sheaves are constant over $(\ov{\eta})^I$}    \label{subsection-constant-over-ov-eta-I}

\begin{prop}  \label{prop-H-I-W-constant-over-ov-eta-I}
$\restr{\mc H_{I, W}^j}{(\ov{\eta})^I}$ is a constant sheaf over $(\ov{\eta})^I$.
\end{prop}
\dem 
By Proposition \ref{prop-H-x-H-eta-I-isom}, the sheaf $\restr{\mc H_{I, W}^j}{(\ov{\eta})^I}$ is ind-smooth over $(\ov{\eta})^I$.
By Proposition \ref{prop-FWeil-factors-though-Weil-I}, the action of $\on{Weil}((\ov{\eta})^I, \ov{\eta_I})$ on $\restr{ \mc H_{I, W}^j }{ \ov{\eta_I} }$ is trivial. 

By Lemma \ref{lem-ind-smooth-specialization} and \ref{subsection-smooth-sheaf-pi-1-rep}, we deduce Proposition \ref{prop-H-I-W-constant-over-ov-eta-I}.
\cqfd

\quad

\section{The cohomology sheaves are constant over $(\ov{\eta})^{I_1} \times_{\Fqbar} \ov u$}   \label{section-cohomology-constant-I-1-I-2}

\sssec{}   \label{subsection-def-I-1-I-2-and-ov-u}
Let $I$ be a disjoint union $I_1 \sqcup I_2$. %Let $(u_i)_{i \in I_2}$ be a family of closed points of $X \sm N$. Let $u$ be a closed point of $\times_{i \in I_2} u_i$. Fix a geometric point $\ov{u}$ over $u$. 
Let $u$ be a closed point of $(X \sm N)^{I_2}$. 
%so that for every $i \in I_2$, the image of $u$ by the projection to the $i$-th factor $\on{pr}_i: (X \sm N)^{I_2} \rightarrow X \sm N$ is a closed point. 
Let $\ov u = \on{Spec} \Fqbar$ be a geometric point over $u$.

%The main results in this section are Proposition \ref{prop-FWeil-factors-though-Weil-I-1}, Proposition \ref{prop-H-x-H-eta-I-1-times-u-isom} and Proposition \ref{prop-H-ov-eta-I-1-ov-s-I-2-is-constant}.

In this section, we will generalize the results in Section 1 from $(\ov{\eta})^I$ to $(\ov{\eta})^{I_1} \times_{\Fqbar} \ov u$. When $I_1 = I$ and $I_2 = \emptyset$, the results in Section 2 recover the results in Section 1.

%prove in Proposition \ref{prop-H-ov-eta-I-1-ov-s-I-2-is-constant} that $\restr{\mc H_{I, W}^j}{(\ov{\eta})^{I_1} \times_{\Fqbar} \ov u}$ is a constant sheaf over $(\ov{\eta})^{I_1} \times_{\Fqbar} \ov u$.

%\begin{rem}
%When $I_1 = I$ and $I_2 = \emptyset$, the results in Section 2 recover the results in Section 1.
%\end{rem}

\subsection{Reminders on partial Frobenius morphisms}

\sssec{}
Let $\eta_{I_1}$ be the generic point of $X^{I_1}$. Fix a geometric point $\ov{ \eta_{I_1} }$ over $\eta_{I_1}$.
By \ref{subsection-ov-eta-I-integral-scheme}, the scheme $(\ov{\eta})^{I_1} \times_{\Fqbar} \ov u$ is an integral scheme. Note that $\ov{ \eta_{I_1} } \times_{\Fqbar} \ov u$ is a geometric generic point of $(\ov{\eta})^{I_1} \times_{\Fqbar} \ov u$. 

%Recall that in \ref{subsection-def-Frob-J} we defined morphisms $\Frob_{\{i\}}$ for any $i \in I_1$ and $\Frob_{I_2}$.

\sssec{}   \label{subsection-def-FWeil-action-I-1}
Let $\mc F$ be an ind-constructible $E$-sheaf over $(\eta)^{I_1} \times_{\Fq} u$, equiped with an action of the partial Frobenius morphisms, i.e. for every $i \in I_1$, an isomorphism $F_{\{i\}}: \on{Frob}_{\{i\}}^* \mc F \isom \mc F$
and an isomorphism $F_{I_2}: \on{Frob}_{I_2}^* \mc F \isom \mc F$,
commuting to each other and whose composition is the total Frobenius isomorphism $\on{Frob}^* \mc F \isom \mc F$ over $(\eta)^{I_1} \times_{\Fq} u$.

Then $\restr{ \mc F }{ \ov{\eta_{I_1} } \times_{\Fqbar} \ov u }$ is equiped with an action of $\on{FWeil}(\eta_{I_1}, \ov{\eta_{I_1}})$ in the following way:
let $\varepsilon \in \on{FWeil}(\eta_{I_1}, \ov{\eta_{I_1}})$ with $\restr{\varepsilon}{(F_{I_1})^{\on{perf}}} = \prod_{i \in {I_1}} (\on{Frob}_{\{i\}})^{n_i}$. By \ref{subsection-def-FWeil-action}, we have a morphism $\mf{sp}_{\varepsilon}: \prod_{i \in I_1} (\on{Frob}_{\{i\}})^{n_i}(\ov{\eta_{I_1}})  \rightarrow \ov{\eta_{I_1}} $ over $\on{Spec} \Fqbar$.
%Note that $\restr{\varepsilon}{\Fqbar}$ is identity.
Thus, tensoring with $\Id_{\ov u}$, we deduce a specialization map in $X^I$ (which is in fact an isomorphism of schemes):
$$\mf{sp}_{\varepsilon} \otimes \Id_{\ov u}: \prod_{i \in I_1} (\on{Frob}_{\{i\}})^{n_i}(\ov{\eta_{I_1}}) \times_{\Fqbar} \ov u \rightarrow \ov{\eta_{I_1}} \times_{\Fqbar} \ov u $$
The action of $\varepsilon$ is defined to be the composition:
$$\restr{\mc F}{\ov{\eta_{I_1}} \times_{\Fqbar} \ov u } \xrightarrow{(\mf{sp}_{\varepsilon} \otimes \Id_{\ov u})^*} \restr{\mc F}{ \prod_{i \in I_1} (\on{Frob}_{\{i\}})^{n_i}(\ov{\eta_{I_1}}) \times_{\Fqbar} \ov u  } \xrightarrow{ \prod_{i \in I_1} (F_{\{i\}})^{n_i} } \restr{\mc F}{ \ov{\eta_{I_1}} \times_{\Fqbar} \ov u  }$$
We deduce an action of $\on{FWeil}(\eta_{I_1}, \ov{\eta_{I_1}})$.

\sssec{}   \label{subsection-when-smooth-sheaf-is-constant-abstract-eta-times-u}
As in \ref{subsection-when-smooth-sheaf-is-constant-abstract}, we begin by considering a too simple case to illustrate the case \ref{subsection-F-A-mod-tf-smooth-then-constant-I-1-I-2}. %which we will need in Sections 2.2-2.3.
Let $\mc F$ be a constructible $E$-sheaf over $(\eta)^{I_1} \times_{\Fq} u$, equiped with an action of the partial Frobenius morphisms.

By Lemma \ref{lem-Drinfeld-smooth-over-U-I} applied to $I_1$, we deduce that $\mc F$ is smooth over $(\ov{\eta})^{I_1} \times_{\Fqbar} \ov u$.

As in \ref{subsection-def-FWeil-action-I-1}, $\restr{ \mc F }{ \ov{ \eta_{I_1} } \times_{\Fqbar} \ov u }$ is equiped with a continuous action of $\on{FWeil}(\eta_{I_1}, \ov{\eta_{I_1}})$. By Lemma \ref{lem-Drinfeld-E-vs}, this action factors through $\on{Weil}(\eta, \ov{\eta})^{I_1}$. We deduce that the action of $\on{Weil}((\ov{\eta})^{I_1}, \ov{\eta_{I_1}})$ on $\restr{ \mc F }{ \ov{\eta_{I_1}} \times_{\Fqbar} \ov u }$ is trivial. 
Note that $$\on{Weil}((\ov{\eta})^{I_1} \times_{\Fqbar} \ov u, \ov{\eta_{I_1}} \times_{\Fqbar} \ov u ) \simeq \on{Weil}((\ov{\eta})^{I_1}, \ov{\eta_{I_1}}) .$$ Thus the action of $\on{Weil}((\ov{\eta})^{I_1} \times_{\Fqbar} \ov u, \ov{\eta_{I_1}} \times_{\Fqbar} \ov u )$ on $\restr{ \mc F }{ \ov{\eta_{I_1}} \times_{\Fqbar} \ov u }$ is also trivial.

By \ref{subsection-smooth-sheaf-pi-1-rep}, $\mc F$ is constant over $(\ov{\eta})^{I_1} \times_{\Fqbar} \ov u$.

\sssec{}  \label{subsection-F-A-mod-tf-smooth-then-constant-I-1-I-2}
As in \ref{subsection-F-A-mod-tf-smooth-then-constant},
let $A$ be a finitely generated $E$-algebra.
Let $\mc F$ be an ind-constructible $E$-sheaf over $(\eta)^{I_1} \times_{\Fq} u$ equiped with an action of the partial Frobenius morphisms and an action of $A$, such that 
\begin{itemize}
\item these two actions commute with each other

\item $\restr{ \mc F }{ \ov{\eta_{I_1}} \times_{\Fqbar} \ov u }$ is an $A$-module of finite type

%\item $\mc F$ is an inductive limit of constructible $E$-sheaves

%the action of $\pi_1^{\on{geom}}(\eta_I, \ov{\eta_I})$ on $\restr{ \mc F }{ \ov{\eta_I} }$ is continuous, 
\end{itemize}
Then as in \ref{subsection-def-FWeil-action-I-1}, $\restr{ \mc F }{ \ov{\eta_{I_1}} \times_{\Fqbar} \ov u }$ is equiped with a continuous $A$-linear action of $\on{FWeil}(\eta_{I_1}, \ov{\eta_{I_1}})$. 
%whose restriction to $\pi_1^{\on{geom}}(\eta_{I_1}, \ov{\eta_{I_1}})$ is continuous, that is to say, $\restr{ \mc F }{ \ov{\eta}_{I_1} \times_{\Fqbar} \ov u }$ is a union of finite dimensional $E$-vector subspaces which are stable under $\pi_1^{\on{geom}}(\eta_{I_1}, \ov{\eta_{I_1}})$ and on which the action of $\pi_1^{\on{geom}}(\eta_{I_1}, \ov{\eta_{I_1}})$ is continuous.
By Lemma \ref{lem-Drinfeld-A-mod}, this action factors through $\on{Weil}(\eta, \ov{\eta})^{I_1}$. Thus the action of $\on{Weil}((\ov{\eta})^{I_1} \times_{\Fqbar} \ov u, \ov{\eta_{I_1}} \times_{\Fqbar} \ov u )$ on $\restr{ \mc F }{ \ov{\eta_{I_1}} \times_{\Fqbar} \ov u }$ is trivial. 

%If in addition $\mc F$ is "ind-smooth" over $(\eta)^{I_1} \times_{\Fq} u$, 
%%that is to say, for any geometric point $\ov{x}$ of $(\ov{\eta})^{I_1}$ and any specialization map $\mf{sp}_x: \ov{\eta}_{I_1} \times_{\Fqbar} \ov u \rightarrow \ov x \times_{\Fqbar} \ov u$, the induced morphism $$(\mf{sp}_x)^*: \restr{  \mc F  }{ \ov x \times_{\Fqbar} \ov u }  \rightarrow \restr{  \mc F  }{ \ov{\eta}_{I_1} \times_{\Fqbar} \ov u  } $$ is an isomorphism,
%then $\mc F$ is constant over $(\ov{\eta})^{I_1} \times_{\Fqbar} \ov u$.

However, %I do not know a generalization of Lemma \ref{lem-Drinfeld-smooth-over-U-I} to ind-constructible sheaves. So we cannot deduce that 
in general $\mc F$ is not ind-smooth over $(\ov{\eta})^{I_1} \times_{\Fqbar} \ov u$.

\sssec{}
In our situation below, we will apply \ref{subsection-F-A-mod-tf-smooth-then-constant-I-1-I-2} to the cohomology sheaves of stacks of shtukas, and prove that for other reasons the cohomology sheaves are ind-smooth over $(\ov{\eta})^{I_1} \times_{\Fqbar} \ov u$. Then they are constant over $(\ov{\eta})^{I_1} \times_{\Fqbar} \ov u$.

%In the following, we will apply this situation to an $E$-sheaf over $(\eta)^I$ which is an inductive limit of subsheaves $\mc F$ as above.

\subsection{Action of $\on{Weil}(\eta, \ov{\eta})^{I_1}$ on cohomology of stacks of shtukas}   

%\begin{prop}   \label{prop-FWeil-factors-though-Weil-I-1}
%The action of $\on{FWeil}(\eta_{I_1}, \ov{\eta_{I_1}})$ on $\restr{\mc H_{I, W}^j}{  (\ov{\eta})^{I_1} \times_{\Fqbar} (\ov s)^{I_2}  }$ factors through $\on{Weil}(\eta, \ov{\eta})^{I_1}$.
%\end{prop}

\sssec{}
Consider the sheaf $\restr{\mc H_{I, W}^j}{ (\eta)^{I_1} \times_{\Fq} u }$ over $(\eta)^{I_1} \times_{\Fq} u$. By \ref{subsection-def-partial-Frob-on-H-I-W}, it is equiped with an action of the partial Frobenius morphisms. As in \ref{subsection-def-FWeil-action-I-1}, the fiber $\restr{  \mc H_{I, W}^j   }{ \ov{\eta_{I_1}} \times_{\Fqbar} \ov u  } $ is equiped with an action of $\on{FWeil}(\eta_{I_1}, \ov{\eta_{I_1}})$.

\begin{prop}   \label{prop-FWeil-factors-though-Weil-I-1}
The action of $\on{FWeil}(\eta_{I_1}, \ov{\eta_{I_1}})$ on $\restr{  \mc H_{I, W}^j   }{ \ov{\eta_{I_1}} \times_{\Fqbar} \ov u  } $ factors through $\on{Weil}(\eta, \ov{\eta})^{I_1}$.
\end{prop}
\dem
The proof is similar to the proof of Proposition \ref{prop-FWeil-factors-though-Weil-I}. 
We use Lemma \ref{lem-H-is-union-of-sub-modules-I-1-I-2} below. Then to each $\mf M$, we apply \ref{subsection-F-A-mod-tf-smooth-then-constant-I-1-I-2} to $A = \otimes_{i \in I_1} \ms H_{G, v_i}$ and $\mc F = \mf M$.
\cqfd

\begin{lem}   \label{lem-H-is-union-of-sub-modules-I-1-I-2}
$\restr{\mc H_{I, W}^j}{\ov{\eta_{I_1}} \times_{\Fqbar} \ov u}$ is an inductive limit of sub-$E$-modules $\mf M$ which are stable under $\on{FWeil}(\eta_{I_1}, \ov{\eta_{I_1}})$, and for which there exists a family $(v_i)_{i \in I_1}$ of closed points in $X \sm N$ (depending on $\mf M$) such that $\mf M$ is stable under the action of $\otimes_{i \in I_1} \ms H_{G, v_i}$ and is of finite type as module over $\otimes_{i \in I_1} \ms H_{G, v_i}$.
\end{lem}
\dem
The proof is similar to the proof of Lemma \ref{lem-H-is-union-of-sub-modules}. 

For any $\mu$, we choose a dense open subscheme $\Omega $ of $(X \sm N)^{I_1} $ such that $\restr{\mc H_{I, W}^{j, \, \leq \mu}}{\Omega \times_{\Fq} u }$ is smooth. We choose a closed point $v$ of $\Omega$.
Let $v_i$ be the image of $v$ under $(X \sm N)^{I_1} \xrightarrow{\on{pr}_i} X \sm N$. Let $\mf{M}_{\mu}$ be the image of
$$\sum_{(n_i)_{i \in I_1} \in \N^{I_1}} (\otimes_{i \in I_1} \ms H_{G, v_i} ) \cdot \restr{\big( \prod_{i \in I_1} F_{\{i\}}^{n_i} ( (\prod_{i \in I_1} \on{Frob}_{\{i\}}^{n_i})^* \mc H_{I, W}^{j, \, \leq \mu}   )  \big)}{\ov{\eta_{I_1}} \times_{\Fqbar} \ov u}  $$
in $\restr{\mc H_{I, W}^j}{  \ov{\eta_{I_1}} \times_{\Fqbar} \ov u  }$.

Then the proof of Lemma \ref{lem-H-is-union-of-sub-modules} works if we replace everywhere $\ov v$ by $\ov v \times_{\Fqbar} \ov u $ and replace $\ov{\eta_I}$ by $\ov{\eta_{I_1}} \times_{\Fqbar} \ov u$.
\cqfd

\subsection{Smoothness of cohomology sheaves over $(\ov{\eta})^{I_1} \times_{\Fqbar} \ov u$}

%\begin{prop}
%$\restr{\mc H_{I, W}^j}{(\ov{\eta})^{I_1} \times_{\Fqbar} (\ov s)^{I_2}}$ is ind-smooth over $(\ov{\eta})^{I_1} \times_{\Fqbar} (\ov s)^{I_2}$.
%\end{prop}

\sssec{}
Let $\ov x$ be a geometric point of $(\ov{\eta})^{I_1}$. Then $\ov{\eta_{I_1}} \times_{\Fqbar} \ov u$ and $\ov x  \times_{\Fqbar} \ov u$ are geometric points of $(X \sm N)^I$.
Let
\begin{equation*}
\mf{sp}_{\ov x}: \ov{\eta_{I_1}} \times_{\Fqbar} \ov u  \rightarrow \ov x  \times_{\Fqbar} \ov u
\end{equation*}
be a specialization map in $(X \sm N)^I$.
It induces a morphism
\begin{equation}    \label{equation-sp-x-u-sheaf}
\mf{sp}_{\ov x}^*: \restr{  \mc H_{I, W}^j   }{ \ov x  \times_{\Fqbar} \ov u }  \rightarrow \restr{  \mc H_{I, W}^j   }{ \ov{\eta_{I_1}} \times_{\Fqbar} \ov u  } 
\end{equation}

\begin{prop}   \label{prop-H-x-H-eta-I-1-times-u-isom}
The morphism (\ref{equation-sp-x-u-sheaf}) is an isomorphism.
\end{prop}
\dem
The proof is similar to the proof of Proposition \ref{prop-H-x-H-eta-I-isom}. We replace everywhere $\ov{x}$ by $\ov x  \times_{\Fqbar} \ov u$ %replace $\ov y$ by $\ov y  \times_{\Fqbar} \ov u$ 
and replace everywhere $\ov{\eta_I}$ by $\ov{\eta_{I_1}} \times_{\Fqbar} \ov u$. 
%The slight difference with the proof of Proposition \ref{prop-H-x-H-eta-I-isom} is that here
%\begin{itemize}
%\item $\Omega_0$ is a dense open subscheme of $\ov{\{x\}}$ such that $\restr{\mc H_{I, W}^{j, \, \leq \mu_0}}{\Omega_0 \times_{\Fq} u}$ is smooth
%
%\item $\Omega_1$ is a dense open subscheme of $(X \sm N)^{I_1}$ such that $\restr{\mc H_{I, W}^{j, \, \leq \mu_1}}{\Omega_1 \times_{\Fq}  u}$ is smooth
%
%\item $\Omega_2$ is a dense open subscheme of $(X \sm N)^{I_1}$ such that $\restr{\mc H_{I, W}^{j, \, \leq \mu_2}}{\Omega_2 \times_{\Fq}  u}$ is smooth.
%\end{itemize}
\cqfd

\subsection{The cohomology sheaves are constant over $(\ov{\eta})^{I_1} \times_{\Fqbar} \ov u$}

\begin{prop}   \label{prop-H-ov-eta-I-1-ov-s-I-2-is-constant}
$\restr{\mc H_{I, W}^j}{(\ov{\eta})^{I_1} \times_{\Fqbar} \ov u}$ is a constant sheaf over $(\ov{\eta})^{I_1} \times_{\Fqbar} \ov u$.
\end{prop}
\dem
By Proposition \ref{prop-H-x-H-eta-I-1-times-u-isom}, the sheaf $\restr{\mc H_{I, W}^j}{(\ov{\eta})^{I_1} \times_{\Fqbar} \ov u}$ is ind-smooth over $(\ov{\eta})^{I_1} \times_{\Fqbar} \ov u$.
By Proposition \ref{prop-FWeil-factors-though-Weil-I-1}, the action of $\on{Weil}((\ov{\eta})^{I_1} \times_{\Fqbar} \ov u, \ov{\eta_{I_1}} \times_{\Fqbar} \ov u )$ on $\restr{ \mc H_{I, W}^j }{ \ov{\eta_{I_1}} \times_{\Fqbar} \ov u }$ is trivial.
By Lemma \ref{lem-ind-smooth-specialization} and \ref{subsection-smooth-sheaf-pi-1-rep}, we deduce Proposition \ref{prop-H-ov-eta-I-1-ov-s-I-2-is-constant}.
\cqfd

\sssec{}
Let $s$ be a closed point of $X \sm N$ and $\ov s = \on{Spec} \Fqbar$ a geometric point over $s$. Let $$(\ov{s})^{I_2}:= \ov{s} \times_{\Fqbar}  \cdots \times_{\Fqbar} \ov{s} .$$
Then $\ov u = (\ov s)^{I_2}$ is a geometric point of $(X \sm N)^{I_2}$ as in \ref{subsection-def-I-1-I-2-and-ov-u}. By Proposition  \ref{prop-H-ov-eta-I-1-ov-s-I-2-is-constant}, $\restr{\mc H_{I, W}^j}{(\ov{\eta})^{I_1} \times_{\Fqbar} (\ov s)^{I_2}}$ is a constant sheaf over $(\ov{\eta})^{I_1} \times_{\Fqbar} (\ov s)^{I_2}$.

\quad

\section{A condition of smoothness}

The goal of this section is to prove Proposition \ref{prop-condition-of-smooth-multi-traits}, which will be used in the next section. To illustrate, we begin with a simple case, that is Proposition \ref{prop-condition-of-smooth-trait}.

\subsection{pseudo-product sheaves}

\sssec{}   \label{subsection-def-S-trait}
As in [SGA7 XIII], %let $k$ be an algebraically closed field. 
let $S$ be a trait (i.e. spectrum of a DVR) %over $k$ 
which is henselian. Let $s = \on{Spec} k$ be the closed point and $\delta = \on{Spec} K$ the generic point. Fix an algebraic closure $\ov K$ of $K$. We denote by $\ov{\delta} = \on{Spec} \ov K$. 
%To simplify the notation, we suppose that $k$ is algebraically closed, i.e. $\ov s = s$. 
It will be enough for us to consider only the case where $k$ is separately closed, i.e. we assume $\ov s = s$. 

In this section, we write $\times$ for $\times_{\on{Spec} k}$.

\sssec{}   
Let $I$ be a finite set. Let $(\ov{x_i})_{i \in I}$ be a family of geometric points of $S$ such that $\ov{x_i} \in \{\ov s, \ov{\delta} \}$. 
We denote by $\times_{i \in I} \ov{x_i}$ the fiber product over $\on{Spec} k$.
As in \ref{subsection-ov-eta-I-integral-scheme}, $\times_{i \in I} \ov{x_i}$ is an integral scheme over $\on{Spec} k$. 

\sssec{}
Let $\Lambda = \mc O_E / \lambda_E^s \mc O_E$ or $\mc O_E$ or $E$.

\begin{defi}   \label{def-pseudo-product-sheaf}
Let $\mc G$ be an ind-constructible $\Lambda$-sheaf over $S^I$. We say that $\mc G$ is a pseudo-product if for any family $(\ov{x_i})_{i \in I}$ of geometric points of $S$ such that $\ov{x_i} \in \{\ov s, \ov{\delta} \}$, the restriction $\restr{\mc G}{ \times_{i \in I} \ov{x_i}}$ is a constant sheaf over $\times_{i \in I} \ov{x_i}$.
\end{defi}

\begin{notation}
For pseudo-product sheaf $\mc G$, we denote
$$\restr{\mc G}{ \times_{i \in I} \ov{x_i}} := \Gamma(\times_{i \in I} \ov{x_i}, \mc G).$$
\end{notation}

\begin{example}
If $\mc G = \boxtimes_{i \in I} \mc F_i$ where every $\mc F_i$ is an ind-constructible $\Lambda$-sheaf over $S$, then $\mc G$ is a pseudo-product.
\end{example}

\quad

%\begin{thm} \label{thm-milne-sheaf-open-closed} ([Milne, II, Theorem 3.10])
%\end{thm}

\begin{lem}  \label{lem-equivalence-over-S}
([SGA7 XIII 1.2.2])
We have a functor $\Theta_S$ from
$$\ms A: = \{ \text{ind-constructible } \Lambda\text{-sheaf } \mc G \text{ over } S\} $$
to
$$\ms B:= \{(\restr{\mc G}{\ov s}, \restr{\mc G}{\ov{\eta}}, \phi), \restr{\mc G}{\ov s} \text{ is a trivial } \on{Gal}(\ov{\delta} / \delta) \text{-module}, \restr{\mc G}{\ov{\eta}} \text{ is a } \on{Gal}(\ov{\delta} / \delta) \text{-module},$$ 
$$\phi: \restr{\mc G}{\ov s} \rightarrow \restr{\mc G}{\ov{\eta}} \text{ is a } \on{Gal}(\ov{\delta} / \delta) \text{-equivariant morphism}\} .$$
It is an equivalence of categories.
\end{lem}

\begin{rem}   \label{rem-locally-constant-equal-constant-S}
In general, we have an exact sequence
$$1 \rightarrow \ms I \rightarrow \on{Gal}(\ov{\eta} / \eta) \rightarrow \on{Gal}(\ov{s} / s) \rightarrow 1$$
where $\ms I$ is the inertia group. 
The morphism $\phi$ factors through
$$\phi: \restr{\mc G}{\ov s} \rightarrow \restr{\mc G}{\ov{\eta}}^{\ms I}  \hookrightarrow \restr{\mc G}{\ov{\eta}} $$
Thus if $\phi$ is an isomorphism, then the action of $\ms I$ on $\restr{\mc G}{\ov{\eta}}$ is trivial.

Here, under our hypothesis that $\ov s = s$, we have $\ms I = \on{Gal}(\ov{\eta} / \eta)$. If $\phi$ is an isomorphism, then the action of $\on{Gal}(\ov{\eta} / \eta)$ on $\restr{\mc G}{\ov{\eta}}$ is trivial. We deduce that an object $(\restr{\mc G}{\ov s}, \restr{\mc G}{\ov{\eta}}, \phi: \restr{\mc G}{\ov s} \isom \restr{\mc G}{\ov{\eta}})$ in the category $\ms B$ corresponds to a constant sheaf over $S$.
Thus a smooth sheaf over $S$ is the same thing as a constant sheaf over $S$. 

Another way to see this fact is that $S$ does not have non trivial etale covering. 
\end{rem}

\quad

\sssec{}
Let $Y$ be a scheme over $s$. We have a commutative diagram
$$
\xymatrix{
Y \times s  \ar[r]^i  \ar[d]
& Y \times S    \ar[d]
& Y \times \delta  \ar[l]_j  \ar[d] \\
s  \ar[r]^i
& S
& \delta  \ar[l]_j
}
$$

\sssec{}
Applying [SGA7 XIII 1.2.4] to $Y = S$ and using Lemma \ref{lem-equivalence-over-S}, we deduce

\begin{lem}   \label{lem-equivalence-over-S-times-S}
We have a functor $\Theta_{S \times S}$ from
$$\ms A: = \{ \text{pseudo-product ind-constructible } \Lambda\text{-sheaf } \mc G \text{ over } S \times S \} $$
to
$$
\begin{aligned}
\ms B:= & \{ (  \restr{\mc G}{\ov s \times \ov s}  , \restr{\mc G}{\ov{\delta} \times \ov s} , \restr{\mc G}{\ov s \times \ov{\delta}} , \restr{\mc G}{\ov{\delta} \times \ov{\delta}}, \phi_{00, 10}, \phi_{00, 01}, \phi_{10, 11}, \phi_{01, 11} ) \\
& \text{ where } \restr{\mc G}{\ov s \times \ov s}  , \restr{\mc G}{\ov{\delta} \times \ov s} , \restr{\mc G}{\ov s \times \ov{\delta}} , \restr{\mc G}{\ov{\delta} \times \ov{\delta}} \text{ are } \on{Gal}(\ov{\delta} / \delta)^2 \text{-modules }, \\
& \text{ the action of } \on{Gal}(\ov{\delta} / \delta)^2 \text{ is } 
\text{trivial on }  \restr{\mc G}{\ov s \times \ov s} , \\
& \text{ factors through the first factor on } \restr{\mc G}{\ov{\delta} \times \ov s} , \\
& \text{ factors through the second factor on } \restr{\mc G}{\ov s \times \ov{\delta}} ; \\
& \phi_{00, 10}, \phi_{00, 01}, \phi_{10, 11}, \phi_{01, 11} \text{ are }\on{Gal}(\ov{\delta} / \delta)^2 \text{-equivariant morphisms } \\
& \text{ such that the following diagram is commutative } \}
\end{aligned}
$$
\begin{equation}   \label{equation-square-sheaf-over-S-S}
\xymatrixrowsep{2pc}
\xymatrixcolsep{4pc}
\xymatrix{
\restr{\mc G}{\ov s \times \ov{\delta} } \ar[r]^{\phi_{01, 11}}
& \restr{\mc G}{\ov{\delta} \times \ov{\delta}}   \\
\restr{\mc G}{\ov s \times \ov s} \ar[r]^{\phi_{00, 10}} \ar[u]^{\phi_{00, 01}} 
& \restr{\mc G}{\ov{\delta} \times \ov s} \ar[u]_{\phi_{10, 11}}
}
\end{equation}
This functor is an equivalence of categories.
\end{lem}

\begin{rem}     \label{rem-locally-constant-equal-constant-S-times-S}
The morphism $\phi_{00, 10}$ (resp. $\phi_{00, 01}$) factors through $$\restr{\mc G}{\ov s \times \ov s} \rightarrow \restr{\mc G}{\ov{\delta} \times \ov s}^{\ms I} \hookrightarrow \restr{\mc G}{\ov{\delta} \times \ov s}$$
$$( \text{resp. } \quad \restr{\mc G}{\ov s \times \ov s} \rightarrow \restr{\mc G}{\ov s \times \ov{\delta}  }^{\ms I} \hookrightarrow \restr{\mc F}{\ov s\times \ov{\delta}  } )$$
where $\ms I$ is the inertia group.
Thus $\phi_{00, 11}: = \phi_{10, 11} \circ \phi_{00, 10} = \phi_{01, 11} \circ \phi_{00, 01}$ factors through
$$\restr{\mc G}{\ov s \times \ov s} \rightarrow \restr{\mc G}{\ov{\delta} \times \ov{\delta}}^{\ms I \times \ms I} \hookrightarrow \restr{\mc G}{\ov{\delta} \times \ov{\delta}} $$
When they are all isomorphism, the action of $\ms I$ on $\restr{\mc G}{\ov{\delta} \times \ov s}$ (resp. $\restr{\mc F}{\ov s \times \ov{\delta} }$) is trivial and the action of $\ms I^2$ on $\restr{\mc F}{\ov{\delta} \times \ov{\delta}} $ is trivial. 
In our case $\ov s = s$, we have $\ms I = \on{Gal}(\ov{\delta} / \delta)$. Via the equivalence of categories in Lemma \ref{lem-equivalence-over-S-times-S}, we deduce that a smooth pseudo-product sheaf over $S \times S$ is the same as a constant sheaf.
\end{rem}

\begin{example}   \label{example-F-1-otimes-F-2}
In Lemma \ref{lem-equivalence-over-S-times-S}, when $\mc G = \mc F_1 \boxtimes \mc F_2$, (\ref{equation-square-sheaf-over-S-S}) coincides with
$$ 
\xymatrixrowsep{2pc}
\xymatrixcolsep{4pc}
\xymatrix{
\restr{\mc F_1}{\ov s} \otimes \restr{\mc F_2}{\ov{\delta}} \ar[r]^{ \phi_{0, 1} \otimes \Id }
& \restr{\mc F_1}{\ov{\delta}} \otimes \restr{\mc F_2}{\ov{\delta}}   \\
\restr{\mc F_1}{\ov s} \otimes \restr{\mc F_2}{\ov s}  \ar[r]^{\phi_{0, 1} \otimes \Id} \ar[u]^{\Id \otimes \phi_{0, 1}} 
& \restr{\mc F_1}{\ov{\delta}} \otimes \restr{\mc F_2}{\ov s}  \ar[u]_{\Id \otimes \phi_{0, 1}}
}
$$
\end{example}

\quad

\sssec{}
Applying [SGA7 XIII 1.2.4] to $Y = S \times S$ and using Lemma \ref{lem-equivalence-over-S-times-S}, we deduce 

\begin{lem}   \label{lem-equivalence-over-S-times-S-times-S}
We have a functor $\Theta_{S \times S \times S}$ from
$$\ms A: = \{ \text{pseudo-product ind-constructible } \Lambda\text{-sheaf } \mc G \text{ over } S \times S \times S \} $$
to
$$
\begin{aligned}
\ms B:= & \{ ( \on{Gal}(\ov{\delta} / \delta)^3 \text{-modules }: \\
& \restr{\mc G}{\ov s \times \ov s \times \ov s}  , \restr{\mc G}{\ov{\delta} \times \ov s \times \ov s} , \restr{\mc G}{\ov s \times \ov{\delta} \times \ov s} , \restr{\mc G}{\ov s \times
 \ov s \times \ov{\delta}}, \restr{\mc G}{\ov{\delta} \times \ov{\delta} \times \ov s}, \restr{\mc G}{\ov{\delta} \times \ov s \times \ov{\delta} }, \restr{\mc G}{\ov s \times \ov{\delta} \times \ov{\delta} }, \restr{\mc G}{\ov{\delta} \times \ov{\delta} \times \ov{\delta}} , \\
& \on{Gal}(\ov{\delta} / \delta)^3 \text{-equivariant morphisms }: \\
 & \phi_{000, 100}, \phi_{000, 010}, \phi_{000, 001}, \cdots, \phi_{110, 111}, \phi_{101, 111}, \phi_{011, 111} ) \\
& \text{ where the action of } \on{Gal}(\ov{\delta} / \delta)^3 \text{ is } 
\text{trivial on }  \restr{\mc G}{\ov s \times \ov s} , \\
& \text{ factors through the first (resp. second, third) factor on } \restr{\mc G}{\ov{\delta} \times \ov s \times \ov s} \text{ (resp.} \restr{\mc G}{\ov s \times \ov{\delta} \times \ov s}, \restr{\mc G}{\ov s \times \ov s \times \ov{\delta} } ), \\
& \text{ factors through the } \{1, 2\} (resp. \{1, 3\}, \{2, 3\})\text{-factors on } \restr{\mc G}{\ov{\delta} \times \ov{\delta} \times \ov s} (resp. \restr{\mc G}{\ov{\delta} \times \ov s \times \ov{\delta} }, \restr{\mc G}{\ov s \times \ov{\delta} \times \ov{\delta} }); \\
& \text{morphisms are such that the following diagram is commutative } \}
\end{aligned}
$$
\begin{equation}   \label{equation-cube-sheaf-over-S-S-S}
\xymatrixrowsep{2pc}
\xymatrixcolsep{3pc}
\xymatrix{
& \restr{\mc G}{\ov{\delta} \times \ov s \times \ov s}  \ar[rr]^{\phi_{100, 110}} \ar[ld]_{\phi_{100, 101}}
& & \restr{\mc G}{ \ov{\delta} \times \ov{\delta} \times \ov s} \ar[ld]^{\phi_{110, 111}} \\
\restr{\mc G}{ \ov{\delta} \times \ov s \times \ov{\delta} }   \ar[rr]^{ \quad \quad \quad \quad \phi_{101, 111}}
& & \restr{\mc G}{ \ov{\delta} \times \ov{\delta} \times \ov{\delta} } \\
& \restr{\mc G}{\ov s \times \ov s \times \ov s}  \ar[rr]^{ \phi_{000, 010}   \quad \quad \quad \quad   }  \ar[uu]^{\phi_{000, 100}}  \ar[ld]_{\phi_{000, 001}}
& & \restr{\mc G}{\ov s \times \ov{\delta} \times \ov s} \ar[uu]_{\phi_{010, 110}} \ar[ld]^{\phi_{010, 011}} \\
\restr{\mc G}{\ov s \times \ov s \times \ov{\delta} }  \ar[uu]^{\phi_{001, 101}}  \ar[rr]^{\phi_{001, 011}}
& & \restr{\mc G}{\ov s \times \ov{\delta} \times \ov{\delta}}  \ar[uu]_{\phi_{011, 111}}
}
\end{equation}
This functor is an equivalence of categories.
\end{lem}

\begin{rem}   
Similarly to Remark \ref{rem-locally-constant-equal-constant-S-times-S}, a smooth pseudo-product sheaf over $S \times S \times S$ is the same as a constant sheaf.
\end{rem}

\sssec{}   \label{subsection-sheaf-S-S-S-restrict-to-S-S}
Let $\mc G$ be a pseudo-product sheaf over $S \times S \times S$. 
We have a partial diagonal morphism
\begin{equation}   \label{equation-partial-diag-Delta-1-2}
S \times S \xrightarrow{(\Delta^{\{1, 2\}}, \, \Id)} S \times S \times S.
\end{equation}
%We denote by $\Delta^{\{1, 2\}}: S \rightarrow S \times S$ the diagonal morphism. 
Then $\restr{\mc G}{\Delta^{\{1, 2\}}(S) \times S}$ is a pseudo-product sheaf over $S \times S$. Applying [SGA7 XIII 1.2.7 (a)] to $Y = S$, $Y' = S \times S$ and the diagonal morphism $\Delta^{\{1, 2\}}: S \rightarrow S \times S$, we deduce that $\Theta_{S \times S}(\restr{\mc G}{\Delta^{\{1, 2\}}(S) \times S})$ coincides with the commutative sub-diagram
$$
\xymatrix{
\restr{\mc G}{\ov{\delta} \times \ov{\delta} \times \ov s}   \ar[r]
& \restr{\mc G}{\ov{\delta} \times \ov{\delta} \times \ov{\delta} } \\
\restr{\mc G}{\ov s \times \ov s \times \ov s}  \ar[r]  \ar[u]
& \restr{\mc G}{\ov s \times \ov s \times \ov{\delta} }  \ar[u]
}$$
of (\ref{equation-cube-sheaf-over-S-S-S}). 

Similarly, we have a partial diagonal morphism
\begin{equation}   \label{equation-partial-diag-Delta-2-3}
S \times S \xrightarrow{(\Id, \, \Delta^{\{2, 3\}})} S \times S \times S.
\end{equation}
$\restr{\mc G}{ S \times \Delta^{\{2, 3\}}(S) }$ is a pseudo-product sheaf over $S \times S$ %Applying [SGA7 XIII 1.2.7 (a)] to $Y = S$, $Y' = S \times S$ and the diagonal morphism $\Delta^{\{2, 3\}}: S \rightarrow S \times S$, we deduce that 
and $\Theta_{S \times S}(\restr{\mc G}{ S \times \Delta^{\{2, 3\}}(S) })$ coincides with the commutative sub-diagram
$$
\xymatrix{
\restr{\mc G}{\ov s \times \ov{\delta} \times \ov{\delta} }   \ar[r]
& \restr{\mc G}{\ov{\delta} \times \ov{\delta} \times \ov{\delta} } \\
\restr{\mc G}{\ov s \times \ov s \times \ov s}  \ar[r]  \ar[u]
& \restr{\mc G}{\ov{\delta} \times \ov s \times \ov s}  \ar[u]
}$$
of (\ref{equation-cube-sheaf-over-S-S-S}).

\subsection{A condition of smoothness}

\begin{prop}   \label{prop-condition-of-smooth-trait}
Let $S$ be a strictly henselian trait as in \ref{subsection-def-S-trait}. Let $\mc F$ be an ind-constructible $E$-sheaf over $S$. 
Suppose that there exists a pseudo-product ind-constructible $E$-sheaf $\mc G$ over $S^{\{1, 2, 3\}}$, and morphisms over $S \times S$
$$\mc C^{\sharp}: \mc F \boxtimes E \rightarrow \restr{\mc G}{ S \times \Delta^{\{2, 3\}}(S) } \quad \text{ and }  \quad \mc C^{\flat}: \restr{\mc G}{\Delta^{\{1, 2\}}(S) \times S} \rightarrow E \boxtimes \mc F$$
where $E$ is the constant sheaf over $S$,
such that the composition of the restriction over $S$
$$\mc F \simeq \mc F \otimes E \xrightarrow{\mc C^{\sharp}} \restr{\mc G}{\Delta^{\{1, 2, 3\}}(S)} \xrightarrow{\mc C^{\flat}} E \otimes \mc F \simeq \mc F $$
is an isomorphism. 
Then $\mc F$ is constant over $S$, i.e. the canonical morphism 
\begin{equation}   \label{equation-F-s-to-F-delta}
\restr{\mc F}{\ov s} \rightarrow \restr{\mc F}{\ov{\delta}}
\end{equation}
is an isomorphism.
\end{prop}
\dem
First we construct a morphism $\restr{\mc F}{\ov{\delta}} \rightarrow \restr{\mc F}{\ov s}$.
Applying the functor $\Theta_{S \times S \times S}$ in Lemma \ref{lem-equivalence-over-S-times-S-times-S} to $\mc G$, we obtain a commutative diagram
\begin{equation}   \label{equation-comm-triangle-sss}
\xymatrix{
\restr{\mc G}{  \ov s \times \ov s \times \ov s  } \ar[r]^{\phi_{000, 100}}  \ar[rd]_{\phi_{000, 110}}
& \restr{\mc G}{  \ov{\delta} \times \ov s \times \ov s  } \ar[d]^{\phi_{100, 110}} \\
& \restr{\mc G}{\ov{\delta} \times \ov{\delta} \times \ov s} 
}
\end{equation}
In particular, we obtain a morphism $\phi_{100, 110}: \restr{\mc G}{\ov{\delta} \times \Delta^{\{2, 3\}}(\ov s)} \rightarrow \restr{\mc G}{\Delta^{\{1, 2\}}(\ov \delta) \times \ov s} $.
We construct the composition of morphisms %$\Upsilon$:
\begin{equation}   \label{equation-inverse-I-singleton}
\restr{\mc F}{\ov{\delta}} \simeq \restr{\mc F}{\ov{\delta}} \otimes \restr{E}{\ov s}   \xrightarrow{\mc C^{\sharp}} \restr{\mc G}{\ov{\delta} \times \Delta^{\{2, 3\}}(\ov s)} \xrightarrow{\phi_{100, 110}} \restr{\mc G}{\Delta^{\{1, 2\}}(\ov \delta) \times \ov s}  \xrightarrow{\mc C^{\flat}} \restr{E}{\ov{\delta}} \otimes \restr{\mc F}{\ov s} \simeq \restr{\mc F}{\ov s}
\end{equation}
We want to prove that (\ref{equation-inverse-I-singleton}) is the inverse of (\ref{equation-F-s-to-F-delta}), up to isomorphism.

\quad

Now we prove that (\ref{equation-inverse-I-singleton}) $\circ$ (\ref{equation-F-s-to-F-delta}) is an isomorphism. This will imply that (\ref{equation-F-s-to-F-delta}) is injective.

Applying the functor $\Theta_{S \times S}$ in Lemma \ref{lem-equivalence-over-S-times-S} to $\mc C^{\sharp}: \mc F \boxtimes E \rightarrow \restr{\mc G}{ S \times \Delta^{\{2, 3\}}(S) }$, we obtain a commutative diagram
\begin{equation}   \label{equation-comm-C-sharp-s-s-to-delta-s}
\xymatrix{
\restr{\mc F}{\ov s} \otimes \restr{E}{\ov s}  \ar[r] \ar[d]^{\mc C^{\sharp}}
& \restr{\mc F}{\ov{\delta}} \otimes \restr{E}{\ov s}  \ar[d]^{\mc C^{\sharp}}  \\
\restr{\mc G}{\ov s \times \Delta^{\{2, 3\}}(\ov s)} \ar[r] 
& \restr{\mc G}{\ov{\delta} \times \Delta^{\{2, 3\}}(\ov s)} 
}
\end{equation}
where the horizontal morphisms are the canonical morphisms $\phi_{00, 10}$ in Lemma \ref{lem-equivalence-over-S-times-S}.

Applying the functor $\Theta_{S \times S}$ in Lemma \ref{lem-equivalence-over-S-times-S} to $\mc C^{\flat}: \restr{\mc G}{\Delta^{\{1, 2\}}(S) \times S} \rightarrow E \boxtimes \mc F$, we obtain a commutative diagram
\begin{equation}   \label{equation-comm-C-flat-s-s-to-delta-s}
\xymatrix{
\restr{\mc G}{\Delta^{\{1, 2\}}(\ov s) \times \ov s} \ar[r] \ar[d]^{\mc C^{\flat}} 
& \restr{\mc G}{\Delta^{\{1, 2\}}(\ov{\delta}) \times \ov s} \ar[d]^{\mc C^{\flat}} \\
\restr{E}{\ov s} \otimes \restr{\mc F}{\ov s} \ar[r]
& \restr{E}{ \ov{\delta} } \otimes \restr{\mc F}{\ov s}
}
\end{equation}
where the horizontal morphisms are the canonical morphisms $\phi_{00, 10}$ in Lemma \ref{lem-equivalence-over-S-times-S}.

By \ref{subsection-sheaf-S-S-S-restrict-to-S-S}, the lower line of (\ref{equation-comm-C-sharp-s-s-to-delta-s}) and the upper line of (\ref{equation-comm-triangle-sss}) coincide. 
The left lower line of (\ref{equation-comm-triangle-sss}) and the upper line of (\ref{equation-comm-C-flat-s-s-to-delta-s}) coincide. 
Combining (\ref{equation-comm-C-sharp-s-s-to-delta-s}), (\ref{equation-comm-triangle-sss}) and (\ref{equation-comm-C-flat-s-s-to-delta-s}), we deduce that the following diagram is commutative:
\begin{equation}   \label{equation-grand-diagram-injectivity}
\xymatrix{
\restr{\mc F}{\ov s} \otimes \restr{E}{\ov s}  \ar[r] \ar[d]_{\mc C^{\sharp}}
& \restr{\mc F}{\ov{\delta}} \otimes \restr{E}{\ov s}  \ar[d]^{\mc C^{\sharp}}  \\
\restr{\mc G}{\Delta^{\{1, 2, 3\}}(\ov s)} \ar[r] \ar[dd]_{\mc C^{\flat}} \ar[rd]
& \restr{\mc G}{\ov{\delta} \times \Delta^{\{2, 3\}}(\ov s)} \ar[d] \\
& \restr{\mc G}{\Delta^{\{1, 2\}}(\ov \delta) \times \ov s} \ar[d]^{\mc C^{\flat}} \\
\restr{E}{\ov s} \otimes \restr{\mc F}{\ov s} \ar[r]
& \restr{E}{\ov{\delta}} \otimes \restr{\mc F}{\ov s}
}
\end{equation}

The right vertical line is (\ref{equation-inverse-I-singleton}).
Taking into account $\restr{\mc F}{\ov s} \otimes \restr{E}{\ov s} \simeq \restr{\mc F}{\ov s}$ and $\restr{\mc F}{\ov{\delta}} \otimes \restr{E}{\ov s} \simeq \restr{\mc F}{\ov{\delta}} $, by Example \ref{example-F-1-otimes-F-2}, we deduce that the upper line of (\ref{equation-grand-diagram-injectivity}) is nothing but the canonical morphism (\ref{equation-F-s-to-F-delta}). 

By the hypothesis, the composition of the left vertical morphisms is an isomorphism. Since $E$ is a constant sheaf over $S$, the lower line is identity. These imply that (\ref{equation-F-s-to-F-delta}) is injective.

\quad

Now we prove that (\ref{equation-F-s-to-F-delta}) $\circ$ (\ref{equation-inverse-I-singleton}) is an isomorphism. This will imply that (\ref{equation-F-s-to-F-delta}) is surjective.

Applying the functor $\Theta_{S \times S}$ in Lemma \ref{lem-equivalence-over-S-times-S} to $\mc C^{\sharp}: \mc F \boxtimes E \rightarrow \restr{\mc G}{ S \times \Delta^{\{2, 3\}}(S) } $, we obtain a commutative diagram:
\begin{equation}   \label{equation-comm-C-sharp-delta-s-to-delta-delta}
\xymatrix{
\restr{\mc F}{\ov{\delta}} \otimes \restr{E}{\ov s} \ar[r] \ar[d]^{ \mc C^{\sharp} }
& \restr{\mc F}{\ov{\delta}} \otimes \restr{E}{\ov{\delta}}   \ar[d]^{ \mc C^{\sharp}  } \\
\restr{\mc G}{\ov{\delta} \times \Delta^{\{2, 3\}}(\ov s)}  \ar[r]
& \restr{\mc G}{\ov{\delta} \times \Delta^{\{2, 3\}}(\ov{\delta})}
}
\end{equation}
where the horizontal morphisms are the canonical morphisms $\phi_{10, 11}$ in Lemma \ref{lem-equivalence-over-S-times-S}.

Applying the functor $\Theta_{S \times S \times S}$ in Lemma \ref{lem-equivalence-over-S-times-S-times-S} to $\mc G$, we obtain a commutative diagram
\begin{equation}   \label{equation-comm-triangle-delta-delta-delta}
\xymatrixrowsep{2pc}
\xymatrixcolsep{4pc}
\xymatrix{
\restr{\mc G}{\ov{\delta} \times \ov s \times \ov s}  \ar[r]^{\phi_{100, 111}}  \ar[d]_{\phi_{100, 110}}
& \restr{\mc G}{\ov{\delta} \times \ov{\delta} \times \ov{\delta}} \\
\restr{\mc G}{\ov{\delta} \times \ov{\delta} \times \ov s} \ar[ru]_{\phi_{110, 111}}
}
\end{equation}

Applying the functor $\Theta_{S \times S}$ in Lemma \ref{lem-equivalence-over-S-times-S} to $\mc C^{\flat}: \restr{\mc G}{\Delta^{\{1, 2\}}(S) \times S} \rightarrow E \boxtimes \mc F$, we obtain a commutative diagram
\begin{equation}   \label{equation-comm-C-flat-delta-s-to-delta-delta}
\xymatrix{
\restr{\mc G}{\Delta^{\{1, 2\}}(\ov{\delta}) \times \ov s}  \ar[r] \ar[d]^{\mc C^{\flat}}
& \restr{\mc G}{ \Delta^{\{1, 2\}}(\ov{\delta}) \times \ov{\delta} } \ar[d]^{\mc C^{\flat}} \\
\restr{E}{\ov{\delta}} \otimes \restr{\mc F}{\ov s} \ar[r]
& \restr{E}{\ov{\delta}} \otimes \restr{\mc F}{\ov{\delta}} 
}
\end{equation}
where the horizontal morphisms are the canonical morphisms $\phi_{10, 11}$ in Lemma \ref{lem-equivalence-over-S-times-S}.

By \ref{subsection-sheaf-S-S-S-restrict-to-S-S}, the lower line of (\ref{equation-comm-C-sharp-delta-s-to-delta-delta}) and the upper line of (\ref{equation-comm-triangle-delta-delta-delta}) coincide. The right lower line of (\ref{equation-comm-triangle-delta-delta-delta}) and the upper line of (\ref{equation-comm-C-flat-delta-s-to-delta-delta}) coincide. Note that $\restr{\mc G}{\ov{\delta} \times \ov{\delta} \times \ov{\delta}}$ is constant, thus
$ \Gamma(\ov{\delta} \times \ov{\delta} \times \ov{\delta}, \mc G) = \restr{\mc G}{\Delta^{\{1, 2, 3\}}(\ov{\delta})}$. 
Combining (\ref{equation-comm-C-sharp-delta-s-to-delta-delta}), (\ref{equation-comm-triangle-delta-delta-delta}) and (\ref{equation-comm-C-flat-delta-s-to-delta-delta}), we obtain a commutative diagram:
\begin{equation}   \label{equation-grand-diagram-surjectivity}
\xymatrixcolsep{4pc}
\xymatrix{
\restr{\mc F}{\ov{\delta}} \otimes \restr{E}{\ov s} \ar[r] \ar[d]^{ \mc C^{\sharp} }
& \restr{\mc F}{\ov{\delta}} \otimes \restr{E}{\ov{\delta}}   \ar[d]^{ \mc C^{\sharp}  } \\
\restr{\mc G}{\ov{\delta} \times \Delta^{\{2, 3\}}(\ov s)}  \ar[r]  \ar[d]
& \restr{\mc G}{\Delta^{\{1, 2, 3\}}(\ov{\delta})} \ar[dd]^{\mc C^{\flat}} \\
\restr{\mc G}{\Delta^{\{1, 2\}}(\ov{\delta}) \times \ov s}  \ar[ru] \ar[d]^{\mc C^{\flat}}  \\
\restr{E}{\ov{\delta}} \otimes \restr{\mc F}{\ov s}   \ar[r]
& \restr{E}{\ov{\delta}} \otimes \restr{\mc F}{\ov{\delta}}
}
\end{equation} 

The left vertical line is (\ref{equation-inverse-I-singleton}).
By Example \ref{example-F-1-otimes-F-2}, we deduce that the lower line of (\ref{equation-grand-diagram-surjectivity}) is nothing but the canonical morphism (\ref{equation-F-s-to-F-delta}). 

By hypothesis, the composition of right vertical morphisms is an isomorphism. Since $E$ is a constant sheaf over $S$, the upper line is identity. 
These imply that (\ref{equation-F-s-to-F-delta}) is surjective.
\cqfd

\subsection{A generalization}

\sssec{}   \label{subsection-F-multi-v-to-multi-eta}
Let $I = \{1, 2, \cdots, n\}$ be a finite set. Let $(S_i)_{i \in I}$ be a family of strictly henselian traits with
$$\ov{s_i} = s_i \rightarrow S_i \leftarrow \delta_i \leftarrow \ov{\delta_i}$$

As a generalization of Lemmas \ref{lem-equivalence-over-S}, \ref{lem-equivalence-over-S-times-S} and \ref{lem-equivalence-over-S-times-S-times-S}, we have

\begin{lem}   \label{lem-equivalence-over-times-S-i}
We have a functor $\Theta_{\times_{i \in I} S_i}$ from
$$\ms A: = \{ \text{pseudo-product ind-constructible } \Lambda\text{-sheaf } \mc F \text{ over } \times_{i \in I} S_i \} $$
to
$$
\begin{aligned}
\ms B:= & \{ ( \times_{i \in I} \on{Gal}(\ov{\delta_i} / \delta_i) \text{-modules } \restr{\mc F}{\times_{i \in I} \ov{x_i} }  \text{ for all families } (\ov{x_i})_{i \in I} \text{ with } \ov{x_i} \in \{\ov{s_i}, \ov{\delta_i} \} , \\
& \times_{i \in I} \on{Gal}(\ov{\delta_i} / \delta_i) \text{-equivariant morphisms } \phi_{(\ov{x_i}), (\ov{y_i})}:\restr{\mc F}{\times_{i \in I} \ov{x_i} }  \rightarrow \restr{\mc F}{\times_{i \in I} \ov{y_i} }  \\
& \text{ for all families }  (\ov{x_i})_{i \in I}, (\ov{y_i})_{i \in I}  \text{ with } \ov{x_i}, \ov{y_i} \in \{\ov{s_i}, \ov{\delta_i} \} \text{ and } \ov{x_i} \text{ a specialization of } \ov{y_i} ) \\
& \text{ where the action of }  \times_{i \in I} \on{Gal}(\ov{\delta_i} / \delta_i) \text{ on } \restr{\mc F}{\times_{i \in I} \ov{x_i} }  \text{ factors through the factors } i \\
& \text{ where } \ov{x_i} = \ov{\delta_i} \, , \text{ the morphisms }
\phi_{(\ov{x_i}), (\ov{y_i})} \text{ are such that the corresponding } \\
& \text{ diagram is commutative }  \}
\end{aligned}
$$
This functor is an equivalence of categories.
\end{lem}
\dem
We argument by induction on $I$. When $I$ is a singleton, this is Lemma \ref{lem-equivalence-over-S}. In general, suppose that we have the equivalence $\Theta_{S_1 \times \cdots \times S_m}$ for some $m \geq 1$. Apply [SGA7 XIII 1.2.4] to $S = S_{m+1}$ and $Y = S_1 \times \cdots \times S_m$. Using the induction hypothesis and the hypothesis that $\mc G$ is pseudo-product, we deduce the equivalence $\Theta_{S_1 \times \cdots  \times S_{m+1}}$.
\cqfd

\begin{rem}   
Similarly to Remark \ref{rem-locally-constant-equal-constant-S-times-S}, a smooth pseudo-product sheaf over $\times_{i \in I} S_i$ corresponds to an object in $\ms B$ where all $\phi_{(\ov{x_i}), (\ov{y_i})}$ are isomorphism, thus the $\times_{i \in I} \on{Gal}(\ov{\delta_i} / \delta_i)$-modules are trivial. As a consequence, a smooth pseudo-product sheaf over $\times_{i \in I} S_i$ is the same as a constant sheaf.
\end{rem}

%Let $\mc F$ be a pseudo-product ind-constructible $\Lambda$-sheaf over $\times_{i \in I} S_i$. 
%Let $(x_i)_{i \in I}$ and $(y_i)_{i \in I}$ be two families of geometric points with $x_i, y_i \in \{\ov{s_i}, \ov{\delta_i} \}$ such that $x_i$ is a specialization of $y_i$.
%As a genelization of \ref{sebsection-sheaf-over-S-times-S-times-S}, we construct by induction an equivalence of categories from
%$$\{ \text{pseudo-product sheaf } \mc F \text{ over } \times_{i \in I} S_i \}$$
%to
%$$\{ \big( \text{for every } (x_i)_{i \in I}, \restr{\mc F}{\times_{i \in I} x_i} = \Gamma(  \times_{i \in I} x_i, \mc F  ), $$ 
%$$\text{morphisms } \phi_{(x_i), (y_i)} : \restr{\mc F}{\times_{i \in I} x_i} \rightarrow \restr{\mc F}{\times_{i \in I}y_i} \text{ such that the induced diagram is commutative }\big)  \}$$
%
%In particular, we have the canonical morphism 
%\begin{equation}   \label{equation-F-multi-s-to-F-multi-delta}
%\restr{\mc F}{\times_{i \in I} x_i} \rightarrow \restr{\mc F}{\times_{i \in I}y_i}
%\end{equation}
%%where $\restr{\mc F}{\times_{i \in I} x_i} = \Gamma(  \times_{i \in I} x_i, \mc F  )$ and $\restr{\mc F}{\times_{i \in I} y_i} = \Gamma(  \times_{i \in I} y_i, \mc F  )$.

\begin{prop}   \label{prop-condition-of-smooth-multi-traits}
Let $\mc F$ be a pseudo-product ind-constructible $E$-sheaf over $\times_{i \in I} S_i$.
Suppose that there exists a pseudo-product ind-constructible $E$-sheaf $\mc G$ over $(\times_{i \in I} S_i)^{\{1, 2, 3\}}$, and morphisms over $(\times_{i \in I} S_i) \times (\times_{i \in I} S_i)$
$$\mc C^{\sharp}: \mc F \boxtimes E \rightarrow \restr{\mc G}{ (\times_{i \in I} S_i) \times \Delta^{\{2, 3\}}(\times_{i \in I} S_i) } \quad \text{ and }  \quad \mc C^{\flat}: \restr{\mc G}{\Delta^{\{1, 2\}}(\times_{i \in I} S_i) \times (\times_{i \in I} S_i)} \rightarrow E \boxtimes \mc F$$
where $E$ is the constant sheaf over $(\times_{i \in I} S_i)$,
such that the composition of the restriction over $\times_{i \in I} S_i$
$$\mc F = \mc F \otimes E \xrightarrow{\mc C^{\sharp}} \restr{\mc G}{\Delta^{\{1, 2, 3\}}(\times_{i \in I} S_i) } \xrightarrow{\mc C^{\flat}} E \otimes \mc F = \mc F $$
is an isomorphism. 

Then for any families of $(\ov{x_i})_{i \in I}$ and $(\ov{y_i})_{i \in I}$ with $\ov{x_i}, \ov{y_i} \in \{\ov{s_i}, \ov{\delta_i} \}$ and $\ov{x_i}$ a specialization of $\ov{y_i}$ as in Lemma \ref{lem-equivalence-over-times-S-i}, the canonical morphism 
\begin{equation}   \label{equation-phi-x-i-y-i}
\phi_{(\ov{x_i}), (\ov{y_i})}:\restr{\mc F}{\times_{i \in I} \ov{x_i} }  \rightarrow \restr{\mc F}{\times_{i \in I} \ov{y_i} }  
\end{equation}
is an isomorphism. In particular, $\mc F$ is a constant sheaf.
\end{prop}
\dem
The proof is similar to the proof of Proposition \ref{prop-condition-of-smooth-trait}.
First we construct a morphism of the inverse direction
\begin{equation}   \label{equation-inverse-I-general}
\begin{aligned}
 \restr{\mc F}{\times_{i \in I} \ov{y_i} } \simeq 
 \restr{\mc F}{\times_{i \in I} \ov{y_i} } \otimes \restr{E}{\times_{i \in I} \ov{x_i} }   
  \xrightarrow{\mc C^{\sharp}} \restr{\mc G}{\times_{i \in I} \ov{y_i}  \times \Delta^{\{2, 3\}}(\times_{i \in I} \ov{x_i} )} \rightarrow  \\
\restr{\mc G}{\Delta^{\{1, 2\}}(\times_{i \in I} \ov{y_i} ) \times \times_{i \in I} \ov{x_i} }  \xrightarrow{\mc C^{\flat}} \restr{E}{\times_{i \in I} \ov{y_i} } \otimes \restr{\mc F}{\times_{i \in I} \ov{x_i} } \simeq \restr{\mc F}{\times_{i \in I} \ov{x_i} } 
 \end{aligned}
\end{equation}

Injectivity of (\ref{equation-phi-x-i-y-i}): using Lemma \ref{lem-equivalence-over-times-S-i}, we construct the following commutative diagram
\begin{equation}  
\xymatrix{
\restr{\mc F}{\times_{i \in I} \ov{x_i} } \otimes \restr{E}{\times_{i \in I} \ov{x_i} }  \ar[r] \ar[d]_{\mc C^{\sharp}}
& \restr{\mc F}{\times_{i \in I} \ov{y_i} } \otimes \restr{E}{\times_{i \in I} \ov{x_i}}  \ar[d]^{\mc C^{\sharp}}  \\
\restr{\mc G}{\Delta^{\{1, 2, 3\}}(\times_{i \in I} \ov{x_i})} \ar[r] \ar[dd]_{\mc C^{\flat}} \ar[rd]
& \restr{\mc G}{\times_{i \in I} \ov{y_i} \times \Delta^{\{2, 3\}}(\times_{i \in I} \ov{x_i})} \ar[d] \\
& \restr{\mc G}{\Delta^{\{1, 2\}}(\times_{i \in I} \ov{y_i} ) \times \times_{i \in I} \ov{x_i}} \ar[d]^{\mc C^{\flat}} \\
\restr{E}{\times_{i \in I} \ov{x_i}} \otimes \restr{\mc F}{\times_{i \in I} \ov{x_i}} \ar[r]
& \restr{E}{\times_{i \in I} \ov{y_i} } \otimes \restr{\mc F}{\times_{i \in I} \ov{x_i}}
}
\end{equation}
The upper line identifies with (\ref{equation-phi-x-i-y-i}). The composition of the left vertical morphisms is an isomorphism. The lower line is identity. Thus (\ref{equation-phi-x-i-y-i}) is injective.

Surjectivity of (\ref{equation-phi-x-i-y-i}): using Lemma \ref{lem-equivalence-over-times-S-i}, we construct the following commutative diagram
\begin{equation}   
\xymatrixcolsep{4pc}
\xymatrix{
\restr{\mc F}{\times_{i \in I} \ov{y_i}} \otimes \restr{E}{ \times_{i \in I} \ov{x_i}} \ar[r]  \ar[d]^{ \mc C^{\sharp} }
&  \restr{\mc F}{\times_{i \in I} \ov{y_i}} \otimes \restr{E}{\times_{i \in I} \ov{y_i}} \ar[d]^{ \mc C^{\sharp}  } \\
\restr{\mc G}{\times_{i \in I} \ov{y_i} \times \Delta^{\{2, 3\}}(\times_{i \in I} \ov{x_i})}  \ar[r] \ar[d]
& \restr{\mc G}{\Delta^{\{1, 2, 3\}}(\times_{i \in I} \ov{y_i})}  \ar[dd]^{\mc C^{\flat}} \\
\restr{\mc G}{\Delta^{\{1, 2\}}(\times_{i \in I} \ov{y_i}) \times \times_{i \in I} \ov{x_i}} \ar[d]^{\mc C^{\flat}}  \ar[ru] \ar[d]^{\mc C^{\flat}}  \\
\restr{E}{\times_{i \in I} \ov{y_i}} \otimes \restr{\mc F}{\times_{i \in I} \ov{x_i}}   \ar[r]
& \restr{E}{\times_{i \in I} \ov{y_i}} \otimes \restr{\mc F}{\times_{i \in I} \ov{y_i}}
}
\end{equation} 
The lower line identifies with (\ref{equation-phi-x-i-y-i}). The composition of the right vertical morphisms is an isomorphism. The upper line is identity. Thus (\ref{equation-phi-x-i-y-i}) is surjective.

\cqfd

%\begin{cor}   \label{cor-condition-of-smooth-multi-traits}
%The pseudo-product ind-constructible $\Lambda$-sheaf $\mc F$ in Proposition \ref{prop-condition-of-smooth-multi-traits} is a constant sheaf over $\times_{i \in I} S_i$. In particular, it is smooth.
%\end{cor}
%\dem By Proposition \ref{prop-condition-of-smooth-multi-traits}, all the morphisms $\phi_{(x_i), (y_i)}$ are isomorphisms.
%By the equivalence of categories in Lemma \ref{lem-equivalence-over-times-S-i}, the sheaf $\mc F$ is a constant sheaf.
%\cqfd

\quad

\section{Smoothness of the cohomology of stacks of shtukas}

The goal of this section is to prove Theorem \ref{thm-H-I-W-indsmooth}. We illustrate the proof by a simple case in Section  \ref{subsection-smppthness-I-singleton}.

\subsection{Example: when $I$ is a singleton}   \label{subsection-smppthness-I-singleton}

\sssec{}
Let $I = \{1\}$ and $W \in \on{Rep}_{E}(\wh G)$. %We have the complexe of cohomology sheaves $\mc H_{\{1\}, W}$ over $X \sm N$. 
For any $j \in \Z$, we have the degree $j$ cohomology sheaf $\mc H_{\{1\}, W}^j$ over $X \sm N$.

\sssec{}   \label{subsection-creation-annihilation-op}
We recall the definition of the creation operator and the annihilation operator in \cite{vincent} Section 5.

Let $W^*$ be the dual representation of $W$. Let $\{e_k\}$ be a basis of $W$ and $\{e_k^*\}$ be the dual basis. 
We denote by ${\bf 1}$ the trivial representation of $\wh G$.
Let $\delta: {\bf 1} \rightarrow W^* \otimes W$ be the morphism sending $1$ to $\sum_{k} e_k^* \otimes e_k$.
The creation operator 
$\mc C_{\delta}^{\sharp, \{2, 3\}}$ is defined to be the composition of morphisms of sheaves over $(X \sm N)  \times (X \sm N)$:
\begin{equation}    \label{equation-creation-op}
\begin{aligned}
\mc H_{\{1\}, W}^j \boxtimes E_{(X \sm N)} \isom \mc H_{\{1, 0\}, W \boxtimes \bf 1}^j \xrightarrow{\mc H(\Id_W \boxtimes \delta)} \mc H_{\{1, 0\}, W \boxtimes (W^* \otimes W) }^j  \\
\underset{\sim}{\xrightarrow{ \chi_{\{2, 3\}}^{-1}} }  \restr{ \mc H_{\{1, 2, 3\}, W \boxtimes W^* \boxtimes W}^j }{ (X \sm N) \times \Delta^{\{2, 3\}}(X \sm N) }
\end{aligned}
\end{equation}
where $\mc H(\Id_W \boxtimes \delta)$ follows from the functoriality and $\chi_{\{2, 3\}}$ is the fusion isomorphism (\cite{vincent} Proposition 4.12) associated to the map
$$\{1, 2, 3\} \twoheadrightarrow \{1, 0\}, \; 1 \mapsto 1, 2 \mapsto 0, 3 \mapsto 0$$

Let $\on{ev}: W \otimes W^* \rightarrow \bf 1$ be the evaluation map. The annihilation operator $\mc C_{\on{ev}}^{\flat, \{1, 2\}}$ is defined to be the composition of morphisms of sheaves over $(X \sm N) \times (X \sm N)$:
\begin{equation}    \label{equation-annihilation-op}
\begin{aligned}
\restr{ \mc H_{\{1, 2, 3\}, W \boxtimes W^* \boxtimes W }^j }{ \Delta^{\{1, 2\}}(X \sm N) \times (X \sm N)  } \underset{\sim}{\xrightarrow{\chi_{\{1, 2\}}}}
\mc H_{\{0, 3\}, (W \otimes W^*) \boxtimes W}^j \\
\xrightarrow{\mc H(\on{ev} \boxtimes \Id_W )} \mc H_{\{0, 3\}, {\bf 1} \boxtimes W}^j \isom  E_{(X \sm N)} \boxtimes \mc H_{\{3\}, W}^j
\end{aligned}
\end{equation}
where $\chi_{\{1, 2\}}$ is the fusion isomorphism associated to the map
$$\{1, 2, 3\} \twoheadrightarrow \{0, 3\}, \; 1 \mapsto 0, 2 \mapsto 0, 3 \mapsto 3$$

\begin{lem}   \label{lem-Zorro}  ("Zorro" lemma, \cite{vincent} (6.18))
The composition of morphisms of sheaves over $X \sm N$:
\begin{equation}
\mc H_{\{1\}, W}^j = \mc H_{\{1\}, W}^j \otimes E \xrightarrow{\mc C_{\delta}^{\sharp, \{2, 3\}}} 
 \restr{ \mc H_{\{1, 2, 3\}, W \boxtimes W^* \boxtimes W}^j }{  \Delta^{\{1, 2, 3\}}(X \sm N) }
 \xrightarrow{\mc C_{\on{ev}}^{\flat, \{1, 2\}}} E \otimes \mc H_{\{3\}, W}^j = \mc H_{\{3\}, W}^j
\end{equation}
is the identity.
\cqfd
\end{lem}

\quad

\begin{prop}   \label{prop-H-1-W-indsmooth}
$\mc H_{\{1\}, W}^j$ is ind-smooth over $X \sm N$.
\end{prop}
\dem
It is enough to prove that for any geometric point $\ov v$ of $X \sm N$ and any specialization map $\mf{sp}: \ov{\eta} \rightarrow \ov v$, the induced morphism $$\mf{sp}^*: \restr{ \mc H_{\{1\}, W}^j  }{\ov v} \rightarrow \restr{ \mc H_{\{1\}, W}^j  }{\ov{\eta}}$$ is an isomorphism.

Since $X \sm N$ is smooth of dimension $1$, the strict henselization $(X \sm N)_{\ov v}$ of $X \sm N$ at $\ov v$ is a trait. 
%Taking into account \ref{subsection-H-1-2-3-is-pseudo-product} and Lemma \ref{lem-Zorro} 
Apply Proposition \ref{prop-condition-of-smooth-trait} to
\begin{itemize}
\item $S = (X \sm N)_{\ov v}$, $\ov s = s = \ov v$, $\ov{\delta} = \ov{\eta}$, $\delta$ the image of $\ov{\eta}$ in $S$,

\item $\mc F = \mc H_{\{1\}, W}^j$,

\item $\mc G = \mc H_{\{1, 2, 3\}, W \boxtimes W^* \boxtimes W}^j$. Note that by Proposition \ref{prop-H-ov-eta-I-1-ov-s-I-2-is-constant}, $\mc H_{\{1, 2, 3\}, W \boxtimes W^* \boxtimes W}^j$ is a pseudo-product sheaf over $S^{\{1, 2, 3\}}$.

\item $\mc C^{\sharp}$ the creation morphism $\mc C_{\delta}^{\sharp, \{2, 3\}}$

\item $\mc C^{\flat}$ the annihilation morphism $\mc C_{\on{ev}}^{\flat, \{1, 2\}}$.
\end{itemize}
By Lemma \ref{lem-Zorro}, the hypothesis of Proposition \ref{prop-condition-of-smooth-trait} is satisfied. We deduce that $\mf{sp}^*$ is an isomorphism.
\cqfd

\quad

\subsection{General case}

Let $I$ be a finite set and $W \in \on{Rep}_{E}(\wh G^I)$.

\sssec{} 
Let $I_0=I_1 =I_2 = I_3= I$.
We denote by
$$
\begin{aligned}
\Delta^{I_1 \sqcup I_2 \sqcup I_3}: (X \sm N)^I & \rightarrow (X \sm N)^{I_1} \times (X \sm N)^{I_2} \times (X \sm N)^{I_3}, \\
(x_i)_{i \in I} & \mapsto \big( (x_i)_{i \in I_1} , (x_i)_{i \in I_2}, (x_i)_{i \in I_3}  \big)
\end{aligned}
$$
$$\Delta^{I_1 \sqcup I_2}: (X \sm N)^I \rightarrow (X \sm N)^{I_1} \times (X \sm N)^{I_2}, \quad (x_i)_{i \in I} \mapsto \big( (x_i)_{i \in I_1} , (x_i)_{i \in I_2}  \big)$$
$$\Delta^{I_2 \sqcup I_3}: (X \sm N)^I \rightarrow (X \sm N)^{I_2} \times (X \sm N)^{I_3}, \quad (x_i)_{i \in I} \mapsto \big( (x_i)_{i \in I_2} , (x_i)_{i \in I_3}  \big)$$

We denote by ${\bf 1}$ the trivial representation of $\wh G^I$. 
%Let $\delta: {\bf 1} \rightarrow W^* \otimes W$ be the morphism sending $1$ to $\sum_{k} e_k^* \otimes e_k$.
As in \cite{vincent} Section 5, we define the creation operator 
$\mc C_{\delta}^{\sharp, I_2 \sqcup I_3}$ (creating paw $I_2 \sqcup I_3$) to be the composition of morphisms of complexes over $(X \sm N)^I  \times (X \sm N)^I$:
\begin{equation}    \label{equation-creation-op-general}
\begin{aligned}
\mc H_{I_1, W} \boxtimes E_{(X \sm N)^I} \isom \mc H_{I_1 \sqcup I_0, W \boxtimes \bf 1} \xrightarrow{\mc H(\Id_W \boxtimes \delta)} \mc H_{I_1 \sqcup I_0, W \boxtimes (W^* \otimes W) }  \\
\underset{\sim}{\xrightarrow{ \chi_{I_2 \sqcup I_3}^{-1}} }  \restr{ \mc H_{I_1 \sqcup I_2 \sqcup I_3, W \boxtimes W^* \boxtimes W} }{ (X \sm N)^{I_1} \times \Delta^{I_2 \sqcup I_3}\big((X \sm N)^I\big) }
\end{aligned}
\end{equation}
where $\chi_{I_2 \sqcup I_3}$ is the fusion isomorphism (\cite{vincent} Proposition 4.12) associated to the map
$$I_1 \sqcup I_2 \sqcup I_3 \twoheadrightarrow I_1 \sqcup I_0$$
sending $I_1$ to $I_1$ by identity, $I_2$ to $I_0$ by identity and $I_3$ to $I_0$ by identity. 

%Let $\on{ev}: W \otimes W^* \rightarrow \bf 1$ be the evaluation map. 

We define the annihilation operator $\mc C_{\on{ev}}^{\flat, I_1 \sqcup I_2}$ (annihilating paws $I_1 \sqcup I_2$) to be the composition of morphisms of complexes over $(X \sm N)^I \times (X \sm N)^I$:
\begin{equation}    \label{equation-annihilation-op-general}
\begin{aligned}
\restr{ \mc H_{I_1 \sqcup I_2 \sqcup I_3, W \boxtimes W^* \boxtimes W } }{ \Delta^{I_1 \sqcup I_2}\big( (X \sm N)^I \big) \times (X \sm N)^{I_3}  } \underset{\sim}{\xrightarrow{\chi_{I_1 \sqcup I_2}}}
\mc H_{I_0 \sqcup I_3, (W \otimes W^*) \boxtimes W} \\
\xrightarrow{\mc H(\on{ev} \boxtimes \Id_W )} \mc H_{I_0 \sqcup I_3, {\bf 1} \boxtimes W} \isom  E_{(X \sm N)^I} \boxtimes \mc H_{I_3, W}
\end{aligned}
\end{equation}
where $\chi_{I_1 \sqcup I_2}$ is the fusion isomorphism associated to the map
$$I_1 \sqcup I_2 \sqcup I_3 \twoheadrightarrow I_0 \sqcup I_3$$
sending $I_1$ to $I_0$ by identity, $I_2$ to $I_0$ by identity and $I_3$ to $I_3$ by identity.

\begin{lem}   \label{lem-Zorro-I-paws}
%$\mc C_{\on{ev}}^{\flat, \{1,2\}} \circ \mc C_{\delta}^{\sharp, \{2, 3\}} = \Id$.
The composition of morphisms of complexes over $(X \sm N)^I$:
\begin{equation}
\mc H_{I_1, W} \xrightarrow{\mc C_{\delta}^{\sharp, I_2 \sqcup I_3}} 
 \restr{ \mc H_{I_1 \sqcup I_2 \sqcup I_3, W \boxtimes W^* \boxtimes W} }{  \Delta^{I_1 \sqcup I_2 \sqcup I_3}\big( (X \sm N)^I \big) }
 \xrightarrow{\mc C_{\on{ev}}^{\flat, I_1 \sqcup I_2}} \mc H_{I_3, W} 
\end{equation}
is the identity.
\end{lem}

\quad

\begin{thm}   \label{thm-H-I-W-indsmooth}
$\mc H_{I, W}^j$ is ind-smooth over $(X \sm N)^I$. 
%In other words, for any geometric point $\ov x$ of $(X \sm N)^I$ and any specialization map 
%$\mf{sp}_{\ov x}: \ov{\eta_I} \rightarrow \ov{x}$,
%the induced morphism
%\begin{equation}  
%\mf{sp}_{\ov x}^*: \restr{ \mc H_{I, W}^j  }{ \ov{x}  } \rightarrow \restr{ \mc H_{I, W}^j  }{ \ov{\eta_I}  }
%\end{equation}
%is an isomorphism.
\end{thm}

\dem
Step (1): 
%For any family $(\ov{v_i})_{i \in I}$ of geometric points of $X \sm N$ and any family of specialization maps $\mf{\sp}_{v_i}: \ov{\eta} \rightarrow \ov{v_i}$, 
For every $i \in I$, for any geometric point $\ov{v_i} = \on{Spec} \Fqbar$ over a closed point of $X \sm N$ and any specialization map $\mf{sp}_{i}: \ov{\eta} \rightarrow \ov{v_i}$, 
we denote by $(X \sm N)_{(\ov{v_i})}$ the strict henselisation of $X \sm N$ at $\ov{v_i}$. Since $X \sm N$ is smooth of dimension $1$, $(X \sm N)_{(\ov{v_i})}$ is a trait. 

Applying Proposition \ref{prop-condition-of-smooth-multi-traits} to 
\begin{itemize}
\item $S_i = (X \sm N)_{\ov{v_i}}$, $\ov s_i = s_i = \ov{v_i}$, $\ov{\delta_i} = \ov{\eta}$, $\delta_i$ the image of $\ov{\eta}$ in $S_i$,

\item $\mc F = \mc H_{I, W}^j$, which is a pseudo-product sheaf over $\times_{i \in I} S_i$ by Proposition \ref{prop-H-ov-eta-I-1-ov-s-I-2-is-constant}, 

\item $\mc G = \mc H_{I_1 \sqcup I_2 \sqcup I_3, W \boxtimes W^* \boxtimes W}^j$, which is a pseudo-product sheaf over $(\times_{i \in I} S_i)^{\{1, 2, 3\}}$ by Proposition \ref{prop-H-ov-eta-I-1-ov-s-I-2-is-constant}, 

\item $\mc C^{\sharp}$ the creation morphism $\mc C_{\delta}^{\sharp, I_2 \sqcup I_3}$

\item $\mc C^{\flat}$ the annihilation morphism $\mc C_{\on{ev}}^{\flat, I_1 \sqcup I_2}$.
\end{itemize}
By Lemma \ref{lem-Zorro-I-paws}, the hypothesis of Proposition \ref{prop-condition-of-smooth-trait} is satisfied. 
We deduce that $\mc H_{I, W}^j$ is constant over $\times_{i \in I} S_i$.

\quad

Step (2): 
We denote by $(X \sm N)^I_{ (\times_{i \in I} \ov{v_i}) }$ the strictly henselisation of $(X \sm N)^I$ at $\times_{i \in I} \ov{v_i}$.
For every $j \in I$, the projections to $j$-th factor $(X \sm N)^I \rightarrow X \sm N$ and $\times_{i \in I} \ov{v_i} \rightarrow \ov{v_j}$ induce a morphism 
$$f_j: (X \sm N)^I_{ (\times_{i \in I} \ov{v_i}) } \rightarrow (X \sm N)_{(\ov{v_j})} = S_j.$$
We deduce a morphism
$$(f_j)_{j \in I}: (X \sm N)^I_{ (\times_{i \in I} \ov{v_i}) } \rightarrow \times_{j \in I} S_j.$$

Note that 
$$(X \sm N)^I_{ (\times_{i \in I} \ov{v_i}) } = (\times_{i \in I} S_i)_{  (\times_{i \in I} \ov{v_i})  }$$
Thus any specialization map $\mf{sp}: \ov{\eta_I} \rightarrow \times_{i \in I} \ov{v_i}$ in $(X \sm N)^I$ 
%we have $$\ov{\eta_I}   \rightarrow   (X \sm N)^I_{ (\times_{i \in I} \ov{v_i}) } \rightarrow \times_{i \in I} S_i.$$
is also a specialization map in $\times_{i \in I} S_i$.
By Step (1), $\mc H_{I, W}$ is constant over $\times_{i \in I} S_i$. Thus the induced morphism
\begin{equation}   \label{equation-H-I-W-times-v-i-to-eta-I}
\mf{sp}^*: \restr{  \mc H_{I, W}^j  }{ \times_{i \in I}\ov{v_i}  } \rightarrow  \restr{  \mc H_{I, W}^j  }{ \ov{\eta_I}  }
\end{equation}
is an isomorphism.

\quad

Step (3): let $\ov x$ be a geometric point of $(X \sm N)^I$ and $\mf{sp}_{\ov x}: \ov{\eta_I} \rightarrow \ov x$ a specialization map. We want to prove that the induced morphism $\mf{sp}_{\ov x}^*: \restr{  \mc H_{I, W}^j  }{ \ov{x}  } \rightarrow \restr{  \mc H_{I, W}^j  }{ \ov{\eta_I}  }$ is an isomorphism.

Surjectivity: there exists $(\ov{v_i})_{i \in I}$ with $\ov{v_i} = \on{Spec} \Fqbar$ a geometric point over a closed point of $X \sm N$ and a specialization map $\ov{x} \rightarrow \times_{i \in I} \ov{v_i}$. We have the induced morphisms
$$\restr{  \mc H_{I, W}^j  }{ \times_{i \in I}\ov{v_i} } \rightarrow \restr{  \mc H_{I, W}^j  }{ \ov{x}  } \xrightarrow{\mf{sp}_{\ov x}^*} \restr{  \mc H_{I, W}^j  }{ \ov{\eta_I}  } $$
By Step (2), the composition is an isomorphism. Thus $\mf{sp}_{\ov x}^*$ is surjective.

Injectivity: let $c \in \restr{  \mc H_{I, W}^j  }{ \ov{x}  } $ such that $\mf{sp}_{\ov x}^*(c)=0$. 
There exists $\mu$ and $\wt c \in \restr{  \mc H_{I, W}^{j, \, \leq \mu}  }{ \ov{x}  } $ such that $c$ is the image of $\wt c$ in $\restr{  \mc H_{I, W}^j  }{ \ov{x}  } $. 
Let $x$ be the image of $\ov x$ in $(X \sm N)^I$ and $\ov{ \{x\} }$ the Zariski closure of $x$ in $(X \sm N)^I$.
There exists $\Omega \subset \ov{ \{x\} }$ such that $\restr{  \mc H_{I, W}^{j, \, \leq \mu}  }{ \Omega  } $ is smooth. 
There exists $(\ov{u_i})_{i \in I}$ a family of $\Fqbar$-points over closed points of $X \sm N$ such that $\times_{i \in I} \ov{u_i}$ is a geometric point of $\Omega$. For any specialization map $\alpha: \ov{x} \rightarrow \times_{i \in I} \ov{u_i}$ in $\Omega$, the induced morphism $\alpha^*: \restr{  \mc H_{I, W}^{j, \, \leq \mu}  }{ \times_{i \in I}\ov{u_i} } \rightarrow \restr{  \mc H_{I, W}^{j, \, \leq \mu}  }{ \ov{x}  }$ is an isomorphism. So there exists $\wt c' \in \restr{  \mc H_{I, W}^{j, \, \leq \mu}  }{ \times_{i \in I}\ov{u_i} } $ such that $\alpha^*(\wt c')=\wt c$. Let $c'$ be the image of $\wt c'$ in $\restr{  \mc H_{I, W}^{j}  }{ \times_{i \in I}\ov{u_i} } $.

%There exists an étale neighbourhood $\Omega$ of $\ov x$ in $\ov{ \{x\} }$ and a section $\wt c \in \Gamma(\Omega, \mc H_{I, W}^j)$ such that the image of $\wt c$ by the restriction $$ \Gamma(\Omega, \mc H_{I, W}^j) \rightarrow  \restr{  \mc H_{I, W}^j  }{ \ov{x}  }$$ is $c$. 
%Note that there exists $(\ov{u_i})_{i \in I}$ a family of $\Fqbar$-points over closed points of $X \sm N$ such that $\times_{i \in I} \ov{u_i}$ is a geometric point of $\Omega$. Let $c'$ be the image of $\wt c$ by the restriction $$ \Gamma(\Omega, \mc H_{I, W}^j) \rightarrow  \restr{  \mc H_{I, W}^j  }{ \times_{i \in I}\ov{u_i} }$$
%For any specialization map $\alpha: \ov{x} \rightarrow \times_{i \in I} \ov{u_i}$ in $\Omega$, the following diagram is commutative
%$$\xymatrix{
%& \Gamma(\Omega, \mc H_{I, W}^j)  \ar[rd]^{\text{restriction}}  \ar[ld]  _{\text{restriction}} \\
%\restr{  \mc H_{I, W}^j  }{ \times_{i \in I}\ov{u_i} }   \ar[rr]^{\alpha^*}
%& & \restr{  \mc H_{I, W}^j  }{ \ov{x}  } 
%}$$

Now consider the morphisms
$$\restr{  \mc H_{I, W}^j  }{ \times_{i \in I}\ov{u_i} } \xrightarrow{\alpha^*} \restr{  \mc H_{I, W}^j  }{ \ov{x}  } \xrightarrow{\mf{sp}_{\ov x}^*} \restr{  \mc H_{I, W}^j  }{ \ov{\eta_I}  } $$
$$c' \mapsto c \mapsto 0$$
By step (2), the composition is an isomorphism. 
We deduce that 
%the image of $c'$ in $\restr{  \mc H_{I, W}^j  }{ \ov{\eta_I}  }$, which is equal to $\mf{sp}_x^*(c)$, 
$c' =0$. Hence $c=0$. We deduce that $\mf{sp}_{\ov x}^*$ is injective.
\cqfd

\quad

\section{Action of $\on{Weil}(X \sm N, \ov{\eta})^I$ on cohomology of stacks of shtukas}

\sssec{}
Let $v = \on{Spec} k(v)$ be a place of $X \sm N$. Fix an embedding $\ov F \subset \ov{F_v}$. Thus we have an inclusion $\on{Weil}(\ov{F_v} / F_v) \subset \on{Weil}(\ov F / F)$. 

Let $\ov{k(v)}$ be the algebraic closure of $k(v)$ in $\ov{F_v}$. 
Let $\ov v = \on{Spec} \ov{k(v)}$ be the geometric point over $v$.
Let $\ms I_v = \on{Ker}( \on{Weil}(\ov{F_v} / F_v) \rightarrow \on{Weil}( \ov{k(v)}  / k(v) )  )$ be the inertia group at $v$.

\sssec{}
Let $S = (X \sm N)_{(\ov v)}$ be the strict henselisation of $X \sm N$ at $\ov{v}$. By step (1) of the proof of Theorem \ref{thm-H-I-W-indsmooth}, the restriction $\restr{\mc H_{I, W}^j}{ S^I }$ is a constant sheaf over $S^I$. As a consequence, the action of $(\ms I_v)^I$ on $\restr{\mc H_{I, W}^j}{ \ov{\eta_I} }$ is trivial.

\sssec{}
The group $\on{Weil}(X \sm N, \ov{\eta})$ is the quotient of $\on{Weil}(\eta, \ov{\eta})$ by the subgroup generated by the $\ms I_v$ for all places $v \subset X \sm N$ and their conjugates. 

For any $v$ and any embedding $\ov F \subset \ov{F_v}$, the action of $(\ms I_v)^I$ on $\restr{\mc H_{I, W}^j}{ \ov{\eta_I} }$ is trivial. We deduce 

\begin{prop}   \label{prop-Weil-eta-factors-through-Weil-X}
The action of $\on{Weil}(\eta, \ov{\eta})^I$ on $\restr{ \mc H_{I, W}^j }{ \ov{\eta_I} }$ (defined in Proposition \ref{prop-FWeil-factors-though-Weil-I}) factors through $\on{Weil}(X \sm N, \ov{\eta})^I$.
\end{prop}

\quad

\section{The case of non necessarily split groups}   \label{section-non-split}

In this section, we use the context of \cite{vincent} Section 12. 

\sssec{}
Let $G$ be a geometrically connected smooth reductive group over $F$. As in \cite{vincent} Section 12.1, let $U$ be the maximal open subscheme of $X$ such that $G$ extends to a smooth reductive group scheme over $U$. 
We choose a parahoric integral model of $G$ on every point of $X \sm U$. %We glue these integral models on $U$ and the formal 
Then by gluing we obtain a smooth group scheme over $X$ that we still denote by $G$. Now $G$ is a smooth group scheme over $X$, reductive over $U$, of parahoric type at every point of $X \sm U$. All fibers of $G$ are geometrically connected.

We denote by $\wh N:=|N| \cup (X \sm U)$.

\sssec{}   \label{subsection-HN-non-split}
We use the Harder-Narasimhan truncation for $GL_r$ as in \cite{vincent} Section 12.1. We choose a vector bundle $\mc V$ of rank $r$ over $X$ equiped with a trivialization of $\on{det}(\mc V)$ and an embedding $\rho: G^{\mr{ad}} \rightarrow SL(\mc V)$. We have $\rho_*: \Bun_{G^{\mr{ad}}} \rightarrow \Bun_{GL_r}^0$.
For any  dominant coweight $\mu$ of $GL_r$, we define $\Bun_{G^{\mr{ad}}}^{\leq \mu}$ to be $(\rho_*)^{-1}( \Bun_{GL_r}^{0, \, \leq \mu} )$. We define $\Bun_{G}^{\leq \mu}$ as inverse image of $\Bun_{G^{\mr{ad}}}^{\leq \mu}$. The inductive systems induced by different choices of $\mc V$ and $\rho$ are compatible. 

\sssec{}
As in \cite{vincent} Section 12.1, let $\wt F$ be the finite Galois extension of $F$ such that $\on{Gal}(\wt F / F)$ is the image of $\on{Gal}(\ov F / F)$ in the group of automorphisms of the Dynkin diagram of $G$.
Let the $L$-group ${}^LG$ be the semi-direct product $\wh G \rtimes \on{Gal}(\wt F / F)$, where the semi-direct product is for the action of $\on{Gal}(\wt F / F)$ on $\wh G$ preserving a pinning.

%We view ${}^LG$ as an algebraic group over $\mb{Q}_{\ell}$. For all extension $E \supset \mb{Q}_{\ell}$, we denote by ${}^LG_E$ the extension of scalars to $E$.

\sssec{}   \label{subsection-def-H-I-W-non-split}
As in \cite{vincent} Section 12.3, we suppose $E$ large enough such that all irreducible representations of ${}^LG$ are defined over $E$. %{\color{blue}What does this mean???}

We denote by $\on{Rep}_{E}(({}^LG)^I)$ the category of finite dimensional $E$-linear representation of $({}^LG)^I$. 
Let $I$ be a finite set and $W \in \on{Rep}_{E}(({}^LG)^I)$. 
As in \cite{vincent} Section 12.3, we have a variant of the geometric Satake equivalence ($loc.cit.$ Théorème 12.16). 
We define the stack of $G$-shtukas $\Cht_{G, N, I, W} / \Xi$ over $(X \sm \wh N)^I$, the canonical perverse sheaf $\mc F_{G, N, I, W}^{\Xi}$ over $\Cht_{G, N, I, W} / \Xi$
and the morphism of paws $\mf p_G: \Cht_{G, N, I, W} / \Xi \rightarrow (X \sm \wh N)^I$.
For any dominant coweight $\mu$ of $GL_r$,
we define the truncated stack of shtukas $\Cht_{G, N, I, W}^{\leq \mu} / \Xi$ where the truncation follows from \ref{subsection-HN-non-split}.  For any $j \in \Z$, we define the sheaf of degree $j$ cohomology with compact support
$$\mc H_{G, N, I, W}^{j, \, \leq \mu}:= R^j (\mf p_G)_! \restr{ \mc F_{G, N, I, W}^{\Xi} }{ \Cht_{G, N, I, W}^{\leq \mu} / \Xi} $$ It is a constructible $E$-sheaf over $(X \sm \wh N)^I$.
We define 
$$\mc H_{G, N, I, W}^j:= \varinjlim _{\mu} \mc H_{G, N, I, W}^{j, \, \leq \mu}$$
in the abstract category of inductive limits of constructible $E$-sheaves over $(X \sm \wh N)^I$.

As in \cite{vincent} Section 12.3, the sheaf $\mc H_{G, N, I, W}^j$ is equiped with an action of the partial Frobenius morphisms and an action of the Hecke algebra.

\sssec{}
The results in Sections 1-5 still hold if we replace everywhere $X \sm N$ by $X \sm \wh N$ and replace everywhere $\wh G$ by ${}^LG$. 
We also replace the Eichler-Shimura relations in \cite{vincent} Sections 6-7 by the more general Eichler-Shimura relations in \cite{vincent} Section 12.3.3. Here are the slight modifications.

%In Lemma \ref{lem-ES-relation}, \ref{subsection-def-Hecke-local} and Lemma \ref{lem-S-equal-T}, we suppose that $v_i$ and $v$ are closed point of $X \sm \wh N$ where $G$ is split. 

%\sssec{}
%For any closed point $v$ of $X \sm \wh N$, let ${}^L G_v$ be the local $L$-group defined in \cite{vincent} Section 12.1. It is a subgroup of ${}^L G$.

\begin{lem}  (\cite{vincent} Section 12.3.3) \label{lem-ES-relation-non-split}
Let $(v_i)_{i \in I}$ be a family of closed points of $X \sm \wh N$.
Let $W = \boxtimes_{i \in I} W_i$ with $W_i \in \on{Rep}_E({}^L G)$. Then there exists $\kappa$, such that for any $\mu$ and any $i \in I$, we have
$$\sum_{\alpha =0}^{\on{dim}W_i} (-1)^{\alpha} S_{\wedge^{\on{dim}W_i - \alpha}W_i, v_i}(F_{\{i\}}^{\on{deg}(v_i)})^{\alpha} =0 \; \text{ in } \; \on{Hom}(\restr{\mc H_{I, W}^{j, \, \leq \mu} }{ \times_{i \in I} v_i}, \restr{\mc H_{I, W}^{j, \, \leq \mu + \kappa} }{ \times_{i \in I} v_i})$$
where $S_{\wedge^{\on{dim}W_i - \alpha}W_i, v_i}: \restr{\mc H_{I, W}^{j, \, \leq \mu} }{ \times_{i \in I} v_i} \rightarrow \restr{\mc H_{I, W}^{j, \, \leq \mu + \kappa} }{ \times_{i \in I} v_i}$ is defined in \cite{vincent} Section 12.3.3. 
\cqfd
\end{lem}

\begin{lem} \label{lem-S-equal-T-non-split} (\cite{vincent} Section 12.3.3)
The operator $S_{\wedge^{\on{dim}W_i - \alpha}W_i, v_i}$, which is a morphism of sheaves over $(X \sm \wh N)^I$, extends the action of the Hecke operator $T(h_{\wedge^{\on{dim}W_i - \alpha}W_i, v_i}) \in \ms H_{G, v_i}$, which is a morphism of sheaves over $(X \sm (\wh N \cup v_i))^I$ defined by Hecke correspondence.
\cqfd
\end{lem}

\begin{rem}
In \cite{vincent}, there are more general statements which use the local $L$-group. But Lemma \ref{lem-ES-relation-non-split} and Lemma \ref{lem-S-equal-T-non-split} are enough for us.
\end{rem}

\begin{lem}   \label{lem-H-is-union-of-sub-modules-non-split}
$\restr{\mc H_{G, N, I, W}^j}{\ov{\eta_I}}$ is an increasing union of $E$-vector subspaces $\mf M$ which are stable by $\on{FWeil}(\eta_I, \ov{\eta_I})$, and for which there exists a family $(v_i)_{i \in I}$ of closed points in $X \sm \wh N$ (depending on $\mf M$) such that $\mf M$ is stable under the action of $\otimes_{i \in I} \ms H_{G, v_i}$ and is of finite type as module over $\otimes_{i \in I} \ms H_{G, v_i}$.
\end{lem}
\dem
Note that the category $\on{Rep}_{E}(({}^LG)^I)$ is semisimple,
it is enough to prove the lemma for $W$ irreducible, which is of the form $W = \boxtimes_{i \in I} W_i$ with $W_i \in \on{Rep}_E({}^LG)$. 
Then the proof is the same as the proof of Lemma \ref{lem-H-is-union-of-sub-modules}. %except that we choose $(v_i)_{i \in I}$ such that $G$ is split at every $v_i$.
\cqfd

\begin{prop}   \label{prop-FWeil-factors-though-Weil-I-non-split}
The action of $\on{FWeil}(\eta_I, \ov{\eta_I})$ on $\restr{\mc H_{G, N, I, W}^j}{\ov{\eta_I}}$ factors through $\on{Weil}(\eta, \ov{\eta})^I$.
\end{prop}
\dem The same as Proposition \ref{prop-FWeil-factors-though-Weil-I}.
\cqfd

\begin{prop}  \label{prop-H-I-W-smooth-over-ov-eta-I-non-split}
$\restr{\mc H_{I, W}^j}{(\ov{\eta})^I}$ is smooth over $(\ov{\eta})^I$.
\end{prop}
\dem
%{\color{blue}Is $\on{Rep}_{E}(({}^LG)^I)$ semisimple???}
The proof is the same as the proof of Proposition \ref{prop-H-x-H-eta-I-isom}.% except that we choose $(y_i)_{i \in I}$ such that $G$ is split at every $y_i$.
\cqfd

\quad

The same argument as in Section 4 proves

\begin{thm}   \label{thm-H-I-W-smooth-non-split}
$\mc H_{G, N, I, W}^j$ is ind-smooth over $(X \sm \wh N)^I$. In other words, for any geometric point $\ov x$ of $(X \sm \wh N)^I$ and any specialization map $\mf{sp}_{\ov x}: \ov{\eta_I} \rightarrow \ov{x} $, the induced morphism
\begin{equation}    \label{equation-sp-v-*-I-general-non-split}
\mf{sp}_{\ov x}^*: \restr{ \mc H_{G, N, I, W}^j  }{ \ov{x}  } \rightarrow \restr{ \mc H_{G, N, I, W}^j  }{ \ov{\eta_I}  }
\end{equation}
is an isomorphism.
\cqfd
\end{thm}

The same argument as in Section 5 proves

\begin{prop}   \label{prop-Weil-eta-factors-through-Weil-X-non-split}
The action of $\on{Weil}(\eta, \ov{\eta})^I$ on $\restr{ \mc H_{G, N, I, W}^j }{ \ov{\eta_I} }$ (defined in Proposition \ref{prop-FWeil-factors-though-Weil-I-non-split}) factors through $\on{Weil}(X \sm \wh N, \ov{\eta})^I$.
\cqfd
\end{prop}

\quad

\section{Smooth cuspidal cohomology sheaf}

In this section, we assume that $G$ is split because the constant term morphisms are written only in the split case in \cite{coho-cusp}. 
%For the general case, \ref{subsection-HN-non-split} is not enough to construct the constant term morphisms, we would have to use the Harder-Narasimhan filtration for not necessarily split connected reductive groups (established in \cite{schieder}), which is more delicate.

\sssec{}
Let $P$ be a parabolic subgroup and $M$ be its Levi quotient.
In \cite{coho-cusp} Section 3.4, we defined the cohomology sheaf of stack of $M$-shtukas: $\mc H_{M, N, I, W}^{' \, j}$ over $(X \sm N)^I$. In $loc.cit.$ Section 3.5, we constructed the constant term morphism of sheaves over $\eta_I$
\begin{equation}  \label{equation-CT-morphism-eta-I}
\mc C_{G, N}^{P, \, j}: \restr{\mc H_{G, N, I, W}^j}{\eta_I} \rightarrow \restr{\mc H_{M, N, I, W}^{' \, j}}{\eta_I}
\end{equation}

Similarly to Theorem \ref{thm-H-I-W-indsmooth}, we prove
\begin{prop}   \label{prop-sp-x-*-M-isom}
For any geometric point $\ov x$ of $(X \sm N)^I$ and any specialization map 
$\mf{sp}_{\ov x}: \ov{\eta_I} \rightarrow \ov{x}$, the induced morphism
\begin{equation}   
\mf{sp}_{\ov x}^*: \restr{ \mc H_{M, N, I, W}^{' \, j}  }{ \ov{x}  } \rightarrow \restr{ \mc H_{M, N, I, W}^{' \, j}  }{ \ov{\eta_I}  }
\end{equation}
is an isomorphism.
\end{prop}

By Theorem \ref{thm-H-I-W-indsmooth} and Proposition \ref{prop-sp-x-*-M-isom}, $\mc H_{G, N, I, W}^j$ and $\mc H_{M, N, I, W}^{' \, j}$ are ind-smooth. We extend morphism (\ref{equation-CT-morphism-eta-I}) to a morphism over $(X \sm N)^I$
\begin{equation}
\mc C_{G, N}^{P, \, j}: \mc H_{G, N, I, W}^j \rightarrow \mc H_{M, N, I, W}^{' \, j}
\end{equation}

\begin{defi}
We define the cuspidal cohomology sheaf over $(X \sm N)^I$ to be
$$\mc H_{G, N, I, W}^{j, \, \on{cusp}}:= \bigcap_{P \subsetneq G} \Ker \mc C_{G, N}^{P, \, j} $$
\end{defi}

\sssec{}
By definition, $\restr{\mc H_{G, N, I, W}^{j, \, \on{cusp}} }{ \ov{\eta_I} }$ is the cuspidal cohomology group $H_{G, N, I, W}^{j, \, \on{cusp}}$ defined in \cite{coho-cusp} Definition 3.5.13.

\begin{prop}
$\mc H_{G, N, I, W}^{j, \, \on{cusp}}$ is a smooth $E$-sheaf over $(X \sm N)^I$.
\end{prop}
\dem
We deduce from Theorem \ref{thm-H-I-W-indsmooth} and Proposition \ref{prop-sp-x-*-M-isom} that for any geometric point $\ov x$ of $(X \sm N)^I$ and any specialization map 
$\mf{sp}_{\ov x}: \ov{\eta_I} \rightarrow \ov{x}$, the induced morphism
\begin{equation}   
\mf{sp}_{\ov x}^*: \restr{ \mc H_{G, N, I, W}^{j, \, \on{cusp}}  }{ \ov{x}  } \rightarrow \restr{ \mc H_{G, N, I, W}^{j, \, \on{cusp}}  }{ \ov{\eta_I}  }
\end{equation}
is an isomorphism. So $\mc H_{G, N, I, W}^{j, \, \on{cusp}} $ is ind-smooth.

Moreover, by \cite{coho-cusp} Theorem 0.0.1, $\restr{ \mc H_{G, N, I, W}^{j, \, \on{cusp}}  }{ \ov{\eta_I}  }$ has finite dimension. Thus $\mc H_{G, N, I, W}^{j, \, \on{cusp}}$ is a constructible sheaf.

We deduce that $\mc H_{G, N, I, W}^{j, \, \on{cusp}} $ is smooth.
\cqfd

\quad

%\sssec{}
%In \cite{coho-inte} Section 1.2 and Section 2.3, we defined the cohomology sheaf of stack of $M$-shtukas: $\mc H_{M, N, I, W}^{' \, j, \, \mc O_E}$ over $(X \sm N)^I$ and constructed the constant term morphism of sheaves over $\eta_I$
%\begin{equation}  \label{equation-CT-morphism-eta-I-Oe}
%\mc C_{G, N}^{P, \, j, \, \mc O_E}: \restr{\mc H_{G, N, I, W}^{j, \, \mc O_E}}{\eta_I} \rightarrow \restr{\mc H_{M, N, I, W}^{' \, j, \mc O_E}}{\eta_I}
%\end{equation}
%
%As in the $E$-coefficients case, we prove that $\mc H_{M, N, I, W}^{' \, j, \, \mc O_E}$ is also "ind-smooth" and extend (\ref{equation-CT-morphism-eta-I-Oe}) to a morphism of sheaves over $(X \sm N)^I$
%\begin{equation}   
%\mc C_{G, N}^{P, \, j, \, \mc O_E}: \mc H_{G, N, I, W}^{j, \, \mc O_E} \rightarrow \mc H_{M, N, I, W}^{' \, j, \, \mc O_E}
%\end{equation}
%We define $\mc H_{G, N, I, W}^{j, \, \mc O_E, \, \on{cusp}}:= \bigcap_{P \subsetneq G} \Ker \mc C_{G, N}^{P, \, j}$.
%
%\begin{prop}
%$\mc H_{G, N, I, W}^{j, \, \mc O_E, \, \on{cusp}}$ is a smooth $\mc O_E$-sheaf over $(X \sm N)^I$.
%\end{prop}
%\dem 
%Since $\mc H_{G, N, I, W}^{j, \, \mc O_E}$ and $\mc H_{M, N, I, W}^{' \, j, \, \mc O_E}$ are "ind-smooth" and $\mc C_{G, N}^{P, \, j, \, \mc O_E}$ is a morphism over $(X \sm N)^I$, we deduce that $\mc H_{G, N, I, W}^{j, \, \mc O_E, \, \on{cusp}}$ is also "ind-smooth".
%
%In addition, by \cite{coho-inte} Theorem 0.0.1, $\restr{ \mc H_{G, N, I, W}^{j, \, \mc O_E, \, \on{cusp}} }{ \ov{\eta_I} }$ is an $\mc O_E$-module of finite type. We deduce the proposition.
%\cqfd 

\quad

\end{document}